\documentclass[handcarry]{aiaa-tc}

\usepackage{color}

\usepackage{tabu}
\usepackage{fourier} 
\usepackage{array}
\usepackage{makecell}

\usepackage{multirow}
\usepackage{placeins}
\usepackage{indentfirst}
\usepackage{multirow}

\usepackage{amssymb}
\usepackage{amsmath,latexsym}
\usepackage{url}
\usepackage{enumerate}
\usepackage{graphicx}

\usepackage{booktabs}
 \newcommand{\ra}[1]{\renewcommand{\arraystretch}{#1}}
 
\usepackage{mathtools}

\usepackage{subcaption}
\usepackage{graphicx}
\usepackage{amsmath}
\usepackage[version=4]{mhchem}
\usepackage{siunitx}
\usepackage{longtable,tabularx}
\usepackage{amsmath}

\usepackage{examplep}
\usepackage[driverfallback=dvipdfm,
bookmarks=true,    
bookmarksopen=true,
bookmarksopenlevel=\maxdimen,
bookmarksdepth=10
]{hyperref}
\hypersetup{
pdfcenterwindow=true, 
pdfpagelabels=true,
pagebackref=true,            
hypertexnames=true,
plainpages=false,
unicode=false,     
pdftoolbar=true,                
    pdfmenubar=true,          
    pdffitwindow=false,         
    pdfstartview={FitH},        
    pdftitle={On Iterative Convergence for Damped Viscous Scheme},
    pdfauthor={Hiroaki Nishikawa, Yoshitaka Nakashima, Norihiko Watanabe},    
    pdfsubject={On Iterative Convergence for Damped Viscous Scheme},       
    pdfcreator={Hiroaki Nishikawa},   
    pdfproducer={Hiroaki Nishikawa}, 
    baseurl={http://www.hiroakinishikawa.com},
    pdfkeywords={}, 
    pdfnewwindow=true,     
    colorlinks=true,       
    linkcolor=black,       
    citecolor=black,       
    filecolor=black,       
    urlcolor=black,        
  breaklinks=true,
  hyperfigures=true,
  backref=true,
breaklinks=false, 
pdfpagelabels,      
pagebackref,        
hypertexnames=true, 
plainpages=false,   
naturalnames        
}
\usepackage[all]{hypcap}

\usepackage{cleveref}

\numberwithin{equation}{section}
\numberwithin{subsubsection}{subsection}
\numberwithin{subsection}{section}

\usepackage[utf8]{inputenc}
\usepackage{dashbox}

\def\dyad{\! \otimes \!}

\title{
Resolving Confusion Over Third Order Accuracy of U-MUSCL}

 \author{
{Emmett Padway}\thanks{Research Scholar ({emmett.padway@nianet.org}), 100 Exploration Way, Hampton, VA 23666 USA, Member AIAA}  \, and
{Hiroaki Nishikawa}\thanks{Research Fellow
({hiro@nianet.org}), 100 Exploration Way, Hampton, VA 23666 USA, Associate Fellow AIAA} \\
  {\normalsize\itshape 
{National Institute of Aerospace}, Hampton, VA 23666, USA}%
}

\def\o6{\frac{1}{6}} 

\def\djk{\partial_{\!\! j \! k} \!}
\def\ddjk{\partial^2_{\!\! j \! k} \!}
\def\dddjk{\partial^3_{\!\! j \! k} \!}

\begin{document}


 \maketitle

\begin{abstract} 
In this paper, we discuss the U-MUSCL reconstruction scheme -- an unstructured-grid extension of Van Leer's $\kappa$-scheme -- proposed by Burg for the edge-based discretization [AIAA Paper 2005-4999]. This technique has been widely used in practical unstructured-grid fluid-dynamics 
solvers but with confusions: e.g., third-order accuracy with $\kappa=1/2$ or $\kappa=1/3$. This paper clarifies some of these 
confusions: e.g., the U-MUSCL 
scheme can be third-order accurate in the point-valued solution with $\kappa=1/3$ on regular grids for linear equations in all dimensions, it can 
be third-order accurate with $\kappa=1/2$ as the QUICK scheme in one dimension. 
It is shown that the U-MUSCL scheme cannot be third-order accurate for nonlinear equations, except a very special case of $\kappa=1/2$ on regular simplex-element grids, but it can be an accurate low-dissipation second-order scheme. It is also shown that U-MUSCL extrapolates a quadratic function exactly with $\kappa=1/2$ on arbitrary
grids provided the gradient is computed by a quadratic least-squares method. 
Two techniques are discussed, which transform the U-MUSCL scheme into being genuinely third-order accurate on a regular grid: an efficient flux-reconstruction  method and a special source term quadrature formula for $\kappa=1/2$. 
\end{abstract}
 
\section{Introduction}
\label{intro} 

Despite decades of progress in high-order methods motivated by a desire for increase in accuracy, practical unstructured-grid computational fluid dynamics (CFD) solvers are still dominated by lower-order methods that can be written in the flux-balance form, where a residual at a node/cell $j$ is defined with a single numerical flux per face:
\begin{eqnarray} \label{eq:res}
{\bf Res}_j = \sum_{k \in \{ k_j\} }  {\Phi}_{jk} A_{jk},
\label{residual00}
\end{eqnarray}
where $\{ k_j\} $ is a set of neighbors, $A_{jk}$ is a face area, and ${\Phi}_{jk}$ is a numerical flux. Typically, such a flux-balance method is second-order accurate at best on unstructured grids. Higher-order accuracy will require significantly more complex algorithms (e.g., high-order flux quadrature) and such are rarely employed in practical codes. Under certain constraints, however, simpler high-order methods are available; one such example is the third-order edge-based 
method, which achieves third-order accuracy on arbitrary grids but only with simplex elements [\citen{katz_sankaran:JSC_DOI,liu_nishikawa_aiaa2017-0738,nishikawa_liu_source_quadrature:jcp2017}]. Another example is the unstructured MUSCL (U-MUSCL) scheme [\citen{burg_umuscl:AIAA2005-4999}], which is claimed to be third-order accurate on regular grids. These schemes are simple and efficient since they are based on a single flux evaluation per edge/face as in Equation (\ref{residual00}), and can be easily incorporated in existing unstructured-grid codes by modifying the solution/flux reconstruction and the source discretization [\citen{nishikawa_liu_source_quadrature:jcp2017}]. The U-MUSCL scheme is particularly simple to implement and is widely employed in many practical codes [\citen{burg_umuscl:AIAA2005-4999,fun3d_website,MurayamaYamamoto_jofa2008,Hashimoto_etal:AIAA2014-0240,white_nishikawa_baurle_aiaa2019-0127,scFLOW:Aviation2020}]. However, there still exist confusions about the U-MUSCL scheme: it is said to be third-order accurate with $\kappa=1/2$ [\citen{burg_umuscl:AIAA2005-4999}] or $\kappa=1/3$ [\citen{cfl3d_website}], $\kappa=-1/6$ [\citen{yang_harris:AIAAJ2016}], for example. Also, as discussed recently in Ref.[\citen{Nishikawa_FakeAccuracy:2020}], in fact, the U-MUSCL scheme cannot reach high-order accuracy for nonlinear problems even in 
one dimension. Therefore, in many CFD simulations, third-order accuracy is not actually achieved even on a regular grid. In many cases, as we will discuss later, improvements are due to lower dissipation not third-order accuracy. This paper clarifies the confusions and discusses 
correct ways to achieve high-order accuracy on regular grids with the U-MUSCL reconstruction scheme. 

The U-MUSCL reconstruction scheme, widely employed in practical unstructured-grid CFD codes today, was first proposed by Burg for the edge-based discretization in Ref.[\citen{burg_umuscl:AIAA2005-4999}]
as an unstructured-grid extension of Van Leer's $\kappa$-scheme   
[\citen{VLeer_Ultimate_III:JCP1977,VAN_LEER_MUSCL_AERODYNAMIC:J1985}]. For a one-dimensional uniform grid, the $\kappa$-scheme leads to a third-order accurate finite-volume advection scheme with $\kappa=1/3$ [\citen{Nishikawa_3rdMUSCL:2020IJNMF,VanLeer_Nishikawa_Ultimate_MUSCL:2021}]. Because of this special property, $\kappa=1/3$ has been a popular choice in MUSCL finite-volume schemes. On the other hand, when Burg proposed its extension,  U-MUSCL, to unstructured grids in Ref.[\citen{burg_umuscl:AIAA2005-4999}], he showed that third-order accuracy was obtained with $\kappa=1/2$ in the node-centered edge-based method for a one-dimensional nonlinear system of conservation laws. Note that Burg used $\chi$ for the parameter but we use $\kappa$ throughout the paper. 
Burg's U-MUSCL scheme must be finite-volume as he integrated source terms, but the numerical solution must not have been
a cell-average as in MUSCL because then third-order accuracy would have been obtained with $\kappa=1/3$. The only possibility that he obtained third-order accuracy
is that the numerical solution is a point value, and thus the scheme corresponds to the QUICK scheme of Leonard [\citen{Leonard_QUICK_CMAME1979,Nishikawa_3rdQUICK:2020}]. He then demonstrated superior performance (not third-order accurate) with $\kappa=1/2$ for 
practical CFD simulations. Since then, there has been a confusion about the value of $\kappa$ for third-order accuracy; one of the main contributions of this 
paper is to clarify the choice of $\kappa$ for third-order accuracy.

It is worth pointing out that a similar unstructured-grid MUSCL scheme had already been proposed in Refs.[\citen{Dervieux_IJNMF1998,DebiezDerieux:CF2000,CamarriSalvettiKoobusDerieux:CF2004}]. This scheme has a parameter $\beta$ (instead of $\kappa$) and is called the $\beta$-scheme. Its basic form is equivalent to the U-MUSCL scheme with $\kappa=1-2\beta$. The $\beta$-scheme was, however, proposed as a finite-difference-like scheme with point-valued solutions at nodes, not as a finite-volume scheme. In fact, it was shown that the scheme achieves third-order accuracy on a regular grid with $\beta=1/3$ (corresponding to $\kappa=1/3$, not $\kappa=1/2$) for linear equations [\citen{DebiezDerieux:CF2000}]. For nonlinear equations, they claim that the scheme cannot be higher-order, citing a `proof' of the impossibility of third-order accurate MUSCL schemes in Ref.[\citen{WuWangSun_NonExistence:1998}]. 
While the proof in Ref.[\citen{WuWangSun_NonExistence:1998}] has been shown to be false [\citen{Nishikawa_3rdMUSCL:2020IJNMF,VanLeer_Nishikawa_Ultimate_MUSCL:2021}], it is correct that the $\beta$-scheme cannot be higher-order accurate for nonlinear equations because 
it is a finite-difference scheme with solution reconstruction, as discussed in Refs.[\citen{Nishikawa_FakeAccuracy:2020,VanLeer_Nishikawa_Ultimate_MUSCL:2021}] 
and will also be discussed in this paper. 

Recently, in a series of papers [\citen{Nishikawa_3rdMUSCL:2020IJNMF,Nishikawa_3rdQUICK:2020,Nishikawa_FakeAccuracy:2020,Nishikawa_FSR_arXiv:2020}], we provided clarifications of the U-MUSCL type schemes. Ref.[\citen{Nishikawa_3rdMUSCL:2020IJNMF}] provides a detailed truncation error analysis of a finite-volume scheme with cell-averaged solutions (i.e., the MUSCL scheme) for a nonlinear conservation law and proves that third-order accuracy is achieved with $\kappa=1/3$. Ref.[\citen{Nishikawa_3rdQUICK:2020}] demonstrates that a finite-volume scheme with point-valued solutions (i.e., the QUICK scheme) is third-order accurate with $\kappa=1/2$.  Ref.[\citen{Nishikawa_FakeAccuracy:2020}] shows that the U-MUSCL scheme achieves third-order accuracy with $\kappa=1/3$ on a regular grid in multi-dimensions as a conservative finite-difference scheme but only for linear equations; a particular class of exact solutions effectively linearizes nonlinear equations and allows third-order accuracy to be observed. In this paper, we will follow these previous papers and discuss how the U-MUSCL scheme can be third-order accurate in one dimension and how it cannot be third-order in higher dimensions and for general nonlinear equations. 

Although not genuinely third-order accurate, the U-MUSCL scheme does improve solutions with $\kappa=1/2$ or $\kappa=1/3$ as demonstrated in many test cases (see, e.g., Refs [\citen{burg_umuscl:AIAA2005-4999,MurayamaYamamoto_jofa2008}]). However, such improvements are largely due to reduced dissipation by 
a value of $\kappa$ close to $1$ (zero dissipation), not by third-order accuracy. As discussed, for example, in Ref.[\citen{nishikawa_liu_aiaa2018-4166}], a third-order
scheme eliminates a leading dispersive error of a second-order scheme. Therefore, it is not dissipation but dispersion that is expected to be improved by a third-order
scheme. Genuine third-order accuracy would be desirable thus for problems requiring 
low dispersion errors such as aeroacoustics, and also for high-fidelity simulations using highly refined grids because second-order schemes are significantly 
less efficient on such grids. To this end, we discuss two techniques to construct U-MUSCL-type schemes that can be genuinely third-order 
accurate on regular grids in multi-dimensions. One is the flux reconstruction as demonstrated with an efficient chain-rule-based formula in Ref.[\citen{Nishikawa_FakeAccuracy:2020}]. 
 The other is a special source quadrature approach that eliminates the leading second-order error for $\kappa=1/2$. To the best of the 
 authors' knowledge, this 
 particular third-order $\kappa=1/2$ scheme is a newly discovered scheme. The chain-rule-based flux-reconstruction scheme can be extended to fourth- 
 and/or fifth-order, but in this paper, we will focus on third-order accuracy, which is sufficient to resolve confusions. Higher-order schemes will be discussed in 
 a subsequent paper. 

In addition, we provide the following clarifications: (1) the jump across a face between the left and right states reconstructed by U-MUSCL vanishes for quadratic functions on arbitrary grids for any value of $\kappa$, (2) the U-MUSCL solution reconstruction is exact for quadratic functions with $\kappa=1/2$ on arbitrary grids provided the gradients are computed by a quadratic least-squares (LSQ) method, and (3) the U-MUSCL reconstruction scheme can be directly applicable to cell-centered methods, 
but second-order accuracy will be lost if the face centroid is not exactly located at halfway between two adjacent centroids. 

The paper is organized as follows. Section 2 describes the edge-based discretization with the U-MUSCL scheme for the Euler equations. Section 3 discusses the clarification of the U-MUSCL scheme. Section 4 describes a genuinely third-order scheme based on efficient flux reconstruction. Section 5 presents the special third-order U-MUSCL scheme with $\kappa=1/2$. Section 6 discusses the quadratic exactness of the U-MUSCL reconstruction formula and vanishing solution jumps with quadratically exact gradients. Section 7 discusses the use of U-MUSCL in cell-centered methods, just for the
sake of completeness. Section 8 presents numerical results. Finally, in Section 9, the paper concludes with remarks.

\section{Edge-Based Discretization with U-MUSCL Reconstruction Scheme}
\label{fv_discretization}
 
\subsection{Edge-based discretization with U-MUSCL}

Consider the Euler equations: 
\begin{eqnarray}
\partial_t {\bf u}  + \mbox{div}  {\cal F} = {\bf s},
\label{diff_form}
\end{eqnarray}
where ${\bf s}$ is a source (or forcing) term vector, and 
 \begin{eqnarray}
  {\bf u} =  \left[  \begin{array}{c} 
               \rho       \\ [1ex]
               \rho {\bf v}     \\ [1ex]
               \rho E     
              \end{array} \right], \quad
  {\cal F} =  
 \left[  \begin{array}{c} 
               \rho {\bf v}      \\  [1ex]
               \rho  {\bf v} \dyad  {\bf v}  + p  {\bf I }  \\ [1ex]
               \rho  {\bf v} H 
              \end{array} \right] ,
              \label{Euler_Conservative}
\end{eqnarray}
where $\rho$ is the density, $p$ is the pressure, ${\bf v}$ is the velocity vector,
$ \dyad$ denotes the dyadic product, ${\bf I}$ is the identity matrix, and $E = H-p / \rho  = (\gamma-1) p /\rho   + {\bf v}^2 /2$ is the specific total energy ($\gamma=1.4$).

\begin{figure}[th!]
\begin{center}
\begin{minipage}[b]{0.7\textwidth}
\begin{center}
        \includegraphics[width=0.5\textwidth]{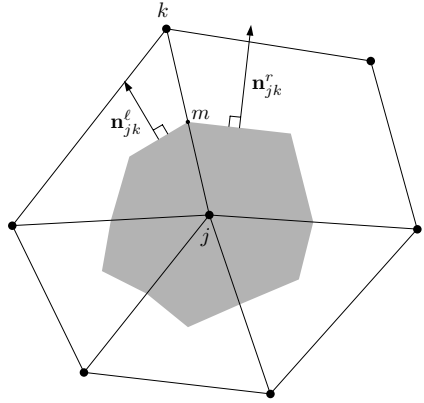}
\caption{Stencil for the edge-based discretization. The directed-area vector across the nodes $j$ and $k$
is denoted by ${\bf n}_{jk}$ and defined by ${\bf n}_{jk} = {\bf n}_{jk}^\ell+ {\bf n}_{jk}^r$.}
\label{fig:eb_stencil}
\end{center}
\end{minipage}
\end{center}
\end{figure}


Following Burg [\citen{burg_umuscl:AIAA2005-4999}], we store the solution values at nodes and discretize the system by the edge-based discretization method over the median dual volume (see Figure \ref{fig:eb_stencil}):
\begin{eqnarray}
  \frac{d {\bf u}_j}{dt} + \frac{1}{V_j}  {\bf Res}_j = {\bf s}_{j}, 
\label{eb_discretization}
\end{eqnarray}
where $V_j$ denotes the measure of the median dual volume around the node $j$ and the residual ${\bf Res}_j$ is defined by the edge-based discretization 
[\citen{andersonbonhaus:CandF1994,Haselbacher_Blazek_AIAAJ2000,Boris_Jim_NIA2007-08,liu_nishikawa_aiaa2016-3969}]: 
\begin{eqnarray}
{\bf Res}_j = \sum_{k \in \{ k_j\} }  {\Phi}_{jk} A_{jk},
\label{residual}
\end{eqnarray}
where $\{ k_j\} $ is a set of edge-connected neighbor nodes, $ {\Phi}_{jk} = {\Phi}_{jk} ( {\bf w}_L, {\bf w}_R, {\bf n}_{jk} )$ is a numerical flux evaluated at the edge midpoint in the direction of the directed-area vector ${\bf n}_{jk}$, ${\bf w}$ denotes a vector of primitive variables ${\bf w} = (\rho, u,v,p)$, $ {\bf w}_L$ and $ {\bf w}_R$ are 
primitive variables reconstructed from the nodes $j$ and $k$, respectively, and $A_{jk}$ is the magnitude of the directed-area vector $A_{jk}=|{\bf n}_{jk}|$. The numerical flux is computed by the Roe flux [\citen{Roe_JCP_1981}]: 
\begin{eqnarray}
 {\Phi}_{jk} =  \frac{1}{2} \left[   {\bf f}({\bf w}_L, \hat{\bf n}_{jk})   +  {\bf f}({\bf w}_R, \hat{\bf n}_{jk})  \right]  - \frac{1}{2} \hat{\bf D}_n \left[  {\bf u}({\bf w}_R) -  {\bf u}({\bf w}_L) \right],
\label{cfd_numerical_flux_2d_euler}
\end{eqnarray}
 where $ {\bf f} = {\cal F} \cdot \hat{\bf n}_{jk}$ is the flux projected along the directed-area vector $\hat{\bf n}_{jk} = {\bf n}_{jk} / |{\bf n}_{jk}|   $, and $\hat{\bf D}_n = | \partial {\bf f} / \partial {\bf u}  |$ is the dissipation term evaluated with the Roe averages [\citen{Roe_JCP_1981}]. We will also employ the Rusanov flux [\citen{Rusanov_Solver}] in accuracy verification tests in two dimensions, where $\hat{\bf D}_n$ is replaced by a scalar value defined by the maximum eigenvalue evaluated with the Roe averages. The left and right states are computed by the U-MUSCL reconstruction scheme, as proposed by Burg [\citen{burg_umuscl:AIAA2005-4999}], 
\begin{eqnarray}
{\bf w}_L  &=&  {\bf w}_{j} + \frac{\kappa}{2} ( {\bf w}_k - {\bf w}_j  ) + \frac{1-\kappa}{2}  \, \overline{\nabla} {\bf w}_j \cdot ( {\bf x}_{k} - {\bf x}_j ), \label{umuscl_of_burg_original_form}  \\[2ex]
{\bf w}_R  &=&  {\bf w}_{k} - \frac{\kappa}{2} ( {\bf w}_k -  {\bf w}_j  ) - \frac{1-\kappa}{2} \, \overline{\nabla} {\bf w}_k \cdot ( {\bf x}_{k} - {\bf x}_j ), 
\label{umuscl_of_burg_original_form_wR}
\end{eqnarray}
where $\overline{\nabla} {\bf w}_j$ and $\overline{\nabla} {\bf w}_k$ are the solution gradients numerically computed by a linear LSQ method with edge-connected neighbors, and $\kappa$ is a user-defined parameter.

\subsection{Equivalent forms}
An interesting equivalent form of the U-MUSCL scheme is given by 
\begin{eqnarray}
{\bf w}_L &=&  \kappa \frac{ {\bf w}_{j}  + {\bf w}_{k} }{2} + (1-\kappa) \left[  {\bf w}_{j} +  \frac{1}{2} \overline{\nabla} {\bf w}_j \cdot ( {\bf x}_{k} - {\bf x}_j )   \right], 
\label{umuscl_of_burg_original_form_wL_ver2}
\\ [2ex]
{\bf w}_R &=&  \kappa \frac{ {\bf w}_{k}  + {\bf w}_{j} }{2} + (1-\kappa) \left[  {\bf w}_{k} -  \frac{1}{2} \overline{\nabla} {\bf w}_k \cdot ( {\bf x}_{k} - {\bf x}_j )   \right],
\label{umuscl_of_burg_original_form_wR_ver2}
\end{eqnarray}
which show that the U-MUSCL scheme is a weighted average of the solution average over the edge and the linear extrapolation using the gradient  [\citen{nishikawa_LP_UMUSCL:JCP2020}].

The U-MUSCL scheme also has the following alternative form: 
\begin{eqnarray}
{\bf w}_L &=& {\bf w}_{j} +
 \frac{1}{4} \left[
   (1-\kappa) {\bf \Delta}^-
+ (1+\kappa) {\bf \Delta}^+ 
\right], 
\\ [2ex]
{\bf w}_R &=& {\bf w}_{k} - 
 \frac{1}{4} \left[   
    (1-\kappa) {\bf \Delta}^+
+ (1+\kappa) {\bf \Delta}^- 
\right], 
\end{eqnarray}
where 
\begin{eqnarray}
{\bf \Delta}^-= 2 \overline{\nabla} {\bf w}_j \cdot ( {\bf x}_{k} - {\bf x}_j ) - ( {\bf w}_{k} - {\bf w}_{j}),
 \quad
{\bf \Delta}^+  = {\bf w}_{k} - {\bf w}_{j},
  \label{kappa_diff_plus_minus_w}
\end{eqnarray}
which, by definition, reduces to Van Leer's $\kappa$-scheme [\citen{VLeer_Ultimate_III:JCP1977,VAN_LEER_MUSCL_AERODYNAMIC:J1985}] on a regular quadrilateral/hexahedral grid with the gradient $\overline{\nabla} {\bf w}_j$ reduced to the central difference gradient. It is noted that the $\beta$-scheme in Ref.[\citen{DebiezDerieux:CF2000}] with a parameter $\beta$ is equivalent to the $\kappa$-scheme with $\kappa=1-2\beta$ but extended to unstructured grids in a manner specific to a node-centered finite-volume discretization, which requires additional neighbor data and is not considered in this paper.

\subsection{U-MUSCL of Yang and Harris}

Yang and Harris [\citen{yang_harris:AIAAJ2016,yang_harris:CCP2018}] extended the U-MUSCL reconstruction scheme by adding a curvature term ${\bf C}_j/2$, 
\begin{eqnarray}
 {\bf w}_L  = {\bf w}_{j} 
 + \frac{\kappa}{2} ( {\bf w}_k - {\bf w}_j  ) + \frac{1-\kappa}{2}  \, \overline{\nabla}{\bf w}_j \cdot ( {\bf x}_{k} - {\bf x}_j )  + \frac{ {\bf C}_j }{2}  ,  
 \label{umuscl_yh_wL}
\end{eqnarray}
where ${\bf C}_j$ is defined, with an additional parameter $\kappa_3$, as 
\begin{eqnarray}
{\bf C}_j  =
\frac{\kappa_3}{4} ( \overline{\nabla}{\bf w}_k - \overline{\nabla}{\bf w}_j  ) \cdot ( {\bf x}_{k} - {\bf x}_j ) + \frac{1-\kappa_3}{4}  \,  ( {\bf x}_{k} - {\bf x}_j )^{\! T} \, \overline{\nabla}^2 {\bf w}_j \cdot ( {\bf x}_{k} - {\bf x}_j ),
\end{eqnarray}
where the superscript $T$ indicates the transpose. 
They set $\kappa=-1/6$ 
and $\kappa_3=-4/3$, and compute the second derivatives, $\overline{\nabla}^2 {\bf w}_j$, by applying a linear LSQ method successively. In Refs.[\citen{yang_harris:AIAAJ2016,yang_harris:CCP2018}], they claim that the scheme is third-order accurate with $\kappa=-1/6$ when $\kappa_3 = 0$: 
\begin{eqnarray}
{\bf w}_L = {\bf w}_{j} + \frac{\kappa}{2} ( {\bf w}_k - {\bf w}_j  ) + \frac{1-\kappa}{2}  \, \overline{\nabla} {\bf w}_j \cdot ( {\bf x}_{k} - {\bf x}_j )   + \frac{1}{2}  \,  ( {\bf x}_{k} - {\bf x}_j )^{\! T} \, \overline{\nabla}^2 {\bf w}_j \cdot ( {\bf x}_{k} - {\bf x}_j ).
\end{eqnarray}
As we will show later, this scheme is essentially equivalent to the original U-MUSCL scheme (\ref{umuscl_of_burg_original_form}) with $\kappa=1/3$ in one dimension. 
It can be third-order accurate only for linear problems [\citen{VanLeer_Nishikawa_Ultimate_MUSCL:2021,Nishikawa_FakeAccuracy:2020}]. 
Also, it requires additional cost and a larger stencil because it involves second-derivatives. This scheme will be further discussed in a subsequent paper
for fourth-order accuracy. 
 
\section{Clarification of U-MUSCL : third-order accurate with $\kappa=1/2$, $\kappa=1/3$, or $\kappa=-1/6$}
\label{truth_umuscl}

In this section, we will clarify the U-MUSCL scheme and discuss how and with what value of $\kappa$ it can be third-order accurate in one dimension. Note that third-order accuracy in one dimension does not necessarily imply third-order accuracy in higher dimensions; but the one-dimensional U-MUSCL scheme is the origin of some of the confusions and thus worth discussing.  In doing so, we will have to make a clear distinction between cell-averaged and point-valued solutions and source terms, whose difference is of second-order and thus matters for third- and higher-order accuracy.  
In short, we will show that the U-MUSCL scheme can be a third-order finite-volume scheme with point-valued solutions ($\kappa=1/2$) or with cell-averaged solutions ($\kappa=1/3$) for general nonlinear equations, or it can be a third-order finite-difference scheme ($\kappa=1/3$) for linear equations. Furthermore, $\kappa=-1/6$ of U-MUSCL of Yang and Harris is equivalent to the third-order finite-difference scheme based on solution reconstruction with $\kappa=1/3$; therefore, it is third-order only for linear 
equations [\citen{VanLeer_Nishikawa_Ultimate_MUSCL:2021,Nishikawa_FakeAccuracy:2020}].

\subsection{U-MUSCL scheme in one dimension}

For simplicity, we consider a scalar nonlinear conservation law:
\begin{eqnarray}
\partial_t u +  \partial_x f = s,
\label{diff_form_oned_cl}
\end{eqnarray}
where $f$ is a nonlinear function of $u$, $f=f(u)$. The U-MUSCL scheme in a one-dimensional grid of uniform spacing $h$, i.e., $x_{i+1}-x_i=h$, is given by
\begin{eqnarray} 
\frac{d {u}_i }{d t} + \frac{ \phi_{i+1/2} - \phi_{i-1/2} }{ h }  = s_i,
\label{fv_form_1D}
\end{eqnarray}
where ${u}_i$ denotes a solution stored at $x=x_i$, and the numerical flux $\phi$ is evaluated by the upwind flux: 
\begin{eqnarray} 
\phi_{i+1/2} = \frac{1}{2} [  f(u_L) + f(u_R) ] - \frac{1}{2} \left|  \frac{\partial f}{\partial u} \right|_{j+1/2} ( u_R -u_L),
\label{oned_numerical_flux}
\end{eqnarray}
where the left and right states are evaluated by the U-MUSCL reconstruction scheme, which is equivalent to Van Leer's scheme in one dimension:
\begin{eqnarray}
{u}_L =  {u}_{i} +
 \frac{1}{4} \left[
   (1-\kappa) {\Delta} u^-
+ (1+\kappa) {\Delta} u^+
\right],  \quad
{u}_R  = {u}_{i+1} - 
 \frac{1}{4} \left[   
    (1-\kappa) {\Delta} u^+
+ (1+\kappa) {\Delta} u^-
\right], 
\label{kappa_uL_and_uR}
\end{eqnarray}
where 
\begin{eqnarray}
{\Delta} u^- =  {u}_{i} - {u}_{i-1},
 \quad
{\Delta} u^+  = {u}_{i+1} - {u}_{i}.
  \label{kappa_diff_plus_minus}
\end{eqnarray}
It is important to note that the U-MUSCL scheme is typically implemented as above without clear definitions of the numerical solution and the source term: point 
value or cell average. The distinction is not important for second-order accuracy but important for third-order accuracy because they differ by $O(h^2)$. In 
what follows, we will discuss three cases where the U-MUSCL scheme achieves third-order accuracy, depending on the definitions of numerical solution 
$u_i$ and the source term evaluation $s_i$.

\subsection{Finite-volume with cell-averaged solution (MUSCL): third-order accurate with $\kappa=1/3$}
\label{oned_FVC}

Suppose the numerical solution is a cell-averaged solution and the source term is cell-averaged, so that the U-MUSCL scheme (\ref{fv_form_1D}) actually implements
\begin{eqnarray} 
\frac{d \overline{u}_i }{d t} + \frac{ \phi_{i+1/2} - \phi_{i-1/2} }{ h }  = \overline{s}_i,
\label{fv_form_1D_ca}
\end{eqnarray}
where   
\begin{eqnarray} 
  \overline{u}_i  = \frac{1}{h} \int_{x_j - h/2}^{x_j+h/2} u \, dx , \quad
  \overline{s}_i  = \frac{1}{h} \int_{x_j - h/2}^{x_j+h/2}   s  \, dx .
\end{eqnarray}
Then, this is a finite-volume scheme with cell-averaged solutions, known as MUSCL of Van 
Leer [\citen{VLeer_Ultimate_IV:JCP1977,VLeer_Ultimate_V:JCP1979}], approximating the exact integral form of the conservation law:
\begin{eqnarray}
  \frac{d  \overline{u}_i }{dt}  +  \overline{f_x} =   \overline{s}_i, 
  \label{integral_form_operator}
\end{eqnarray}
where
\begin{eqnarray}
\overline{f_x}  =   \frac{1}{h} \int_{x_i - h/2}^{x_i+h/2} f_x \, dx  =  \frac{1}{h}  [  f(u_{i+1/2}) - f(u_{i-1/2}) ] .
\end{eqnarray}
It is important to note that arbitrary high-order accuracy is possible in the MUSCL scheme because the integral form is exact. The accuracy of the discretization is then determined solely by the solution reconstruction, at least in one dimension: second-order accuracy with a linear reconstruction and third-order accurate with a quadratic reconstruction. Note also that the reconstruction is performed to obtain a point-valued solution at a face from cell-averaged numerical solutions stored at cells. The choice $\kappa=1/3$ does it exactly for quadratic functions (see Ref.[\citen{Nishikawa_3rdMUSCL:2020IJNMF}], for example). 

As derived in Ref.[\citen{Nishikawa_3rdMUSCL:2020IJNMF}], the truncation error of the scheme (\ref{fv_form_1D_ca}) is given by
\begin{eqnarray}
 \frac{d  \overline{u}_i }{dt}  + 
 \overline{f_x} +  \frac{3  \kappa -1}{12} \left[    f_{uu} u_x  u_{xx}   +  f_u  u_{xxx} 
\right] h^2
 + O(h^3) 
=   \overline{s}_i.
\label{fv_exact_form_te_not_yet001}
\end{eqnarray}
Comparing this with the target integral form  (\ref{integral_form_operator}), we find that the U-MUSCL scheme in one dimension achieves third-order accuracy in the cell-averaged 
solutions with $\kappa=1/3$ if the numerical solution is taken as a cell-averaged solution and the source term is cell averaged.

\subsection{Finite-volume with point-valued solution (QUICK): third-order accurate with $\kappa=1/2$}
\label{oned_FVP}

In the finite-volume scheme, the numerical solution does not have to be a cell average and may be taken as a point value at the cell center. A simple way to construct a third-order finite-volume scheme with point values is to convert the cell-averaged solution $\overline{u}_i$ in the time derivative
to the point value ${u}_i$ in the MUSCL finite-volume scheme (\ref{fv_form_1D_ca}). To do so, we express the cell-average $ \overline{u}_i $ in terms of the point values, $u_{i+1}$, $ u_i$, and $ u_{i-1}$, by cell-averaging the Taylor expansion:
\begin{eqnarray} 
  \overline{u}_i = \frac{1}{h} 
\int_{x_i-h/2}^{x_i+h/2}  \left[  u_i + (x-x_i) u_x + \frac{1}{2} (x-x_i)^2 u_{xx}   \cdots \right] dx 
  ,
\end{eqnarray}
and approximating the second derivative by the central difference approximation:
\begin{eqnarray} 
  \overline{u}_i = u_i + \frac{h^2}{24} \frac{ u_{i+1} - 2 u_i + u_{i-1} }{h^2} + O(h^4) =  u_i + \frac{1}{24} \left( u_{i+1} - 2 u_i + u_{i-1}  \right)+ O(h^4).
\end{eqnarray}
Substituting this into Equation (\ref{fv_form_1D_ca}), we obtain
\begin{eqnarray} 
\frac{d  {u}_i }{d t}  +    \frac{1}{24} \left( \frac{d u_{i+1}}{d t}  - 2 \frac{d u_i}{d t}  + \frac{d u_{i-1}}{d t}  \right) +\frac{ \phi_{i+1/2} - \phi_{i-1/2} }{ h } =   \overline{s}_i + O(h^4) .
\label{fv_form_1D_point_QUICK}
\end{eqnarray} 
This is a finite-volume scheme that updates the point-valued solution; the QUICK scheme of Leonard [\citen{Leonard_QUICK_CMAME1979}] is a good example.
Note that the time derivative terms are coupled with neighbors and thus a global linear system would need to be solved at each stage of an explicit 
time-stepping scheme. 
Other strategies are available for expressing the time derivative terms with point-valued solutions as discussed in Ref.[\citen{Nishikawa_3rdQUICK:2020}]. 
An example is shown in Ref.[\citen{Denaro:IJNMF1996}] for a general construction of finite-volume discretizations with point-valued solutions.

Because of the switch to the point-valued solution, we obtain a slightly different truncation error, as derived in Ref.[\citen{Nishikawa_3rdQUICK:2020}], 
\begin{eqnarray} 
\frac{d  {u}_i }{d t}  +    \frac{1}{24} \left( \frac{d u_{i+1}}{d t}  - 2 \frac{d u_i}{d t}  + \frac{d u_{i-1}}{d t}  \right) + 
   \overline{f_x}      +  \frac{2 \kappa -1}{8} \left[     f_{uu} u_x  u_{xx}   +   f_u  u_{xxx}  
\right] h^2   =    \overline{s}_i + O(h^3),
\label{QUICK_coupled}
\end{eqnarray}
or as an approximation to the integral form (\ref{integral_form_operator}), 
\begin{eqnarray} 
\frac{d  \overline{u}_i }{d t}   + 
   \overline{f_x}      +  \frac{2 \kappa -1}{8} \left[     f_{uu} u_x  u_{xx}   +   f_u  u_{xxx}  
\right] h^2   =    \overline{s}_i + O(h^3),
\label{QUICK_coupled2}
\end{eqnarray}
which shows that the scheme is third-order accurate at $\kappa=1/2$. 

Comparing this with the MUSCL scheme in the previous section, we conclude that the U-MUSCL scheme (\ref{fv_form_1D}) achieves third-order accuracy in one dimension with $\kappa=1/3$ if the numerical solution is a cell average, or with $\kappa=1/2$ if the numerical solution is a point value. 
In both cases, source terms need to be cell averaged. 
In Ref.[\citen{burg_umuscl:AIAA2005-4999}], Burg performed an accuracy verification test for a one-dimensional steady problem with a cell-averaged source term and apparently with point-valued solutions. Therefore, his scheme corresponds to the QUICK scheme [\citen{Nishikawa_3rdQUICK:2020}]. This is the reason that he obtained third-order accuracy with $\kappa=1/2$, not $\kappa=1/3$. If he had solved an unsteady problem, he would not have been able to obtain third-order accuracy without a consistent time-derivative treatment, e.g., as in Equation (\ref{fv_form_1D_point_QUICK}).

\subsection{Finite-difference scheme: third-order accurate with $\kappa=1/3$ for linear equations}
\label{1d_fd_fr}

If both the solution and the source term are point values as typical in a practical implementation of the U-MUSCL scheme {\color{black} (and equivalently the $\beta$-scheme [\citen{Dervieux_IJNMF1998,DebiezDerieux:CF2000,CamarriSalvettiKoobusDerieux:CF2004}])},
\begin{eqnarray} 
\frac{d {u}_i }{d t} + \frac{ \phi_{i+1/2} - \phi_{i-1/2} }{ h }  = s_i,
\label{fv_form_1D_FD}
\end{eqnarray}
where $s_i = s(x_i)$, then, this scheme is second-order accurate at best for nonlinear equations. As derived in Ref.[\citen{Nishikawa_3rdQUICK:2020}], the truncation error is given by
\begin{eqnarray} 
\frac{d {u}_i }{d t} +  f_x = s_i -   \frac{1}{24} \left[     f_{uuu} (u_x)^3 +  6 \kappa  f_{uu} u_x  u_{xx}   +  2 \left(
 3 \kappa - 1
  \right) f_u  u_{xxx}
\right] h^2 + O (h^3).
\label{umuscl_1d_te_2nd}
\end{eqnarray}
This shows that the third-order accuracy is achieved with $\kappa=1/3$ for linear equations, for which $ f_{uuu}= f_{uu}=0$; this fact has also been pointed out in Refs.[\citen{ZhangZhangShu2011,NLV6_INRIA_report:2008}]. Although there are high-order accuracy verification studies seen in the literature [\citen{yang_harris:AIAAJ2016,yang_harris:CCP2018,DementRuffin:aiaa2018-1305}], these high-order results are due to unexpected linearization of nonlinear equations by a particular class of exact solutions (i.e., a function with a small perturbation). Later, we will show an example of this phenomenon. See 
Ref.[\citen{Nishikawa_FakeAccuracy:2020}] for further details. 

Note that the scheme here cannot achieve third-order accuracy with $\kappa=1/2$. In Refs.[\citen{yang_harris:AIAAJ2016,yang_harris:CCP2018}], it is stated that the U-MUSCL scheme is third-order accurate with $\kappa=1/2$ as presented by Burg, but numerical results show only second-order accuracy. This is because the U-MUSCL scheme is treated as a finite-difference scheme as in Equation (\ref{fv_form_1D_FD}), not as a finite-volume scheme (\ref{fv_form_1D_ca}), and therefore it 
cannot be third-order accurate for general nonlinear
problems. See Refs.[\citen{VanLeer_Nishikawa_Ultimate_MUSCL:2021,Nishikawa_FakeAccuracy:2020}] for further details. 


\subsection{U-MUSCL of Yang and Harris: $\kappa=-1/6$}
\label{1d_umuscl_yh}

In one dimension, the U-MUSCL scheme of Yang and Harris (\ref{umuscl_yh_wL}) with $\kappa_3 = 0$ reduces to: 
 \begin{eqnarray}
{u}_L 
= 
{u}_{i} +
 \frac{1}{4} 
 \left[
  (1-\kappa) {\Delta} u^{-}
  + 
  (1+\kappa) {\Delta} u^{+}
\right] 
+
  \frac{1}{32}  \, 
(   u_{i+2} - 2 u_j + u_{i-2}    ).
\label{kappa_uL_and_uR_YH}
\end{eqnarray}
where, the extra term affects only the curvature term in the original U-MUSCL scheme.
They claim that the resulting scheme is third-order accurate with $\kappa=-1/6$ but it is essentially equivalent to the U-MUSCL of Burg 
with $\kappa=1/3$. To see this, let us write it in the quadratic reconstruction form: 
 \begin{eqnarray}
{u}_L 
= 
{u}_{i} +
\frac{1}{2} \frac{ u_{i+1} - u_{i-1} }{2h} \left( \frac{h}{2} \right)
+
\left[ 
\frac{3 \kappa}{2}
\frac{  u_{i+1} - 2 u_i + u_{i-1}  }{h^2}
+
\frac{3}{4}
\frac{  u_{i+2} - 2 u_i + u_{i-2}  }{(2h)^2} 
\right]
\left[
 \left( \frac{h}{2} \right)^2
 - 
  \frac{h^2}{12}
\right].
\label{kappa_uL_and_uR_YH2}
\end{eqnarray}
Without the extra term, the scheme becomes third-order accurate for linear equations with $\kappa=1/3$ as discussed in the previous section; therefore, the coefficient to the curvature term should be $3 \kappa/2=1/2$. Obviously, a different value of $\kappa$ is required
with the extra term. For third-order accuracy, it suffices to 
consider a quadratic solution. Let $u_{xx}$ be the constant second derivative of the quadratic solution. Then, we can write 
\begin{eqnarray}
{u}_L 
= 
{u}_{i} +
\frac{1}{2} \frac{ u_{i+1} - u_{i-1} }{2h} \left( \frac{h}{2} \right)
+
\left[ 
\frac{3 \kappa}{2} 
+
\frac{3}{4} 
\right]  u_{xx} 
\left[
 \left( \frac{h}{2} \right)^2
 - 
  \frac{h^2}{12}
\right],
\label{kappa_uL_and_uR_YH3}
\end{eqnarray}
and solving
\begin{eqnarray}  \frac{3 \kappa}{2} 
+
\frac{3}{4} 
= 
\frac{1}{2},
\end{eqnarray}
for $\kappa$, we obtain
\begin{eqnarray} 
 \kappa = -
\frac{1}{6} .
\end{eqnarray}
Therefore, the U-MUSCL scheme of Yang and Harris is essentially equivalent (although with a wider stencil) to the original U-MUSCL scheme with $\kappa=1/3$. It can be a third-order finite-volume scheme with cell-averaged solutions and a cell-averaged source term, or a third-order finite-difference scheme (for linear equations only). For this reason, we will not discuss their scheme any further and focus on the original U-MUSCL scheme in the rest of the paper.

\subsection{Remarks}
\label{umuscl_remarks}

As we have seen, the U-MUSCL scheme can be third-order accurate in one dimension in three different ways: (1) FVC, a finite-volume scheme with cell-averaged solutions and $\kappa=1/3$, (2) FVP, a finite-volume scheme with point-valued solutions and $\kappa=1/2$, (3) FD, a conservative finite-difference scheme with point-valued solutions, only for linear equations. 

Among these cases, the third-order finite-volume schemes, FVC and FVP, do not extend to higher dimensions with a single flux evaluation per face in the form  (\ref{eb_discretization}); higher-order flux quadrature will be required to achieve third-order accuracy as a finite-volume scheme. There exist techniques to achieve high-order without high-order quadrature [\citen{FVWENO:JSC2014,TamakiImamura:CF2017}], but they are currently applicable to Cartesian grids only (our interest is to develop an unstructured-grid scheme that reduces to a high-order scheme when a grid happens to be regular).  Therefore, the finite-difference scheme is the only way that the U-MUSCL scheme (\ref{eb_discretization}) can achieve third-order accuracy on a regular grid in higher dimensions, as it does not require high-order flux quadrature. But it can be high-order only for linear equations.

Therefore, all the existing U-MUSCL schemes [\citen{burg_umuscl:AIAA2005-4999,yang_harris:AIAAJ2016,yang_harris:CCP2018,DementRuffin:aiaa2018-1305}] are second-order accurate at best for nonlinear equations such as the Euler equations even on a regular grid. Nevertheless, these schemes have been demonstrated to 
provide improved resolution (see also Ref.[\citen{ZhangZhangShu2011}]); therefore they are useful second-order schemes for practical simulations. 
As we will show, these improvements are largely due to low dissipation with $\kappa=1/2$ in the case of the U-MUSCL scheme of Burg, or even smaller dissipation 
by an extended U-MUSCL scheme with second derivatives [\citen{yang_harris:AIAAJ2016,yang_harris:CCP2018}]. The latter will be discussed further in a subsequent 
paper. Yet, genuine third- or higher-order accuracy is sought for long-time convection problems requiring a low dispersion property (e.g., aeroacoustics) and 
for high-fidelity simulations requiring highly refined grids, for which high-order schemes are expected to be significantly more efficient than second-order schemes.
For completeness, we will discuss two techniques to achieve genuine third-order accuracy on regular grids.

\section{A Chain-Ruled Flux-Solution-Reconstruction Scheme with $\theta=1/3$ (CFSR3)}
\label{cfsr}

In this section, we discuss the U-MUSCL scheme made genuinely third-order accurate on regular grids by flux reconstruction. It is similar to the third-order edge-based method, but different in that it is third-order accurate only on regular grids with point-valued source and time-derivative terms whereas the third-order edge-based method is third-order on arbitrary simplex-element grids and requires special source quadrature formulas [\citen{nishikawa_liu_source_quadrature:jcp2017}].

\subsection{One dimension}

We consider the finite-difference-type U-MUSCL scheme (\ref{fv_form_1D_FD}) in Section \ref{1d_fd_fr}. This scheme cannot be third-order accurate for nonlinear equations because of the second-order error generated from the averaged flux term in the numerical flux  (\ref{oned_numerical_flux}), not the dissipation [\citen{Nishikawa_3rdMUSCL:2020IJNMF,Nishikawa_3rdQUICK:2020,Nishikawa_FakeAccuracy:2020}]. To eliminate this error, we need to reconstruct the flux directly, e.g., by applying the U-MUSCL flux reconstruction with a parameter $\theta$: 
\begin{eqnarray} 
\phi_{i+1/2} = \frac{1}{2} [  f_L + f_R] - \frac{1}{2} \left|  \frac{\partial f}{\partial u} \right|_{j+1/2} ( u_R -u_L),
\label{oned_numerical_flux_FR}
\end{eqnarray}
where $ f_L$ and $ f_R$ are given by
\begin{eqnarray}
{f}_L =  {f}_{i} +
 \frac{1}{4} \left[
   (1-\theta) {\Delta}^-
+ (1+\theta) {\Delta}^+ 
\right],  \quad
{f}_R  = {f}_{i+1} - 
 \frac{1}{4} \left[   
(1-\theta) {\Delta}^+
+ (1+\theta) {\Delta}^-
\right], 
\end{eqnarray}
\begin{eqnarray}
{\Delta}^-=  {f}_{i} - {f}_{i-1},
 \quad
{\Delta}^+  = {f}_{i+1} - {f}_{i}.
  \label{kappa_diff_plus_minus_FR}
\end{eqnarray}
This flux-reconstruction-based scheme has the following truncation error, 
\begin{eqnarray}
 \frac{d  {u}_i }{dt}  +  \frac{ 3 \theta - 1}{12} f_{xxx} h^2  + O(h^3) =  {s}_i,
\label{fd_oned_quickest_FR}
\end{eqnarray}
and therefore achieves third-order accuracy with $\theta=1/3$. Note that $u_L$ and $u_R$ are computed with the $\kappa$-reconstruction (\ref{kappa_uL_and_uR}) and the dissipation term does not contribute to the leading second-order error. Therefore, the choice of $\kappa$ is arbitrary; the scheme is third-order as long as $\theta=1/3$. In numerical experiments, we will set $\kappa=1/2$ unless otherwise stated. This scheme is referred to as the flux-solution-reconstruction (FSR) scheme 
[\citen{Nishikawa_FSR_arXiv:2020}]. It is third-order accurate for general nonlinear equations and, unlike FVC and FVP, extends to multi-dimensions.

Note that the FSR scheme is a conservative finite-difference scheme, approximating the differential form 
(\ref{diff_form_oned_cl}) at a point $x=x_i$, which may be the cell center. It is similar to the kappa-family finite-difference scheme of Van Leer [\citen{VLeer_Ultimate_III:JCP1977}], the QUICKEST scheme of Leonard [\citen{Leonard_QUICK_CMAME1979}], the schemes of Shu and Osher [\citen{Shu_Osher_Efficient_ENO_II_JCP1989}] and the NLV6 scheme [\citen{NLV6_INRIA_report:2008}].  However, the flux reconstruction can be very expensive in multi-dimensions. In the next section, we will discuss an efficient chain-rule-based flux reconstruction.

\subsection{Multi-dimensions}

Consider the U-MUSCL scheme on a regular quadrilateral grid, which can be expressed as
\begin{eqnarray}
\frac{d {\bf u}_{i,j}}{dt }  + \frac{ {\Phi}_{i,i+1} - {\Phi}_{i,i-1}  }{h_x} +  \frac{  {\Phi}_{j, j+1} - {\Phi}_{j,j-1}   }{h_y}= {\bf s}_{i,j},
\label{semi_discrete}
\end{eqnarray} 
where $(i,j)$ notation is used for the sake of convenience, $h_x = x_{i+1,j} - x_{i,j} $, $ h_y =  y_{i,j+1} - y_{i,j} $, $s_{i,j} = s(x_{i,j} , y_{i,j} )$, 
${\Phi}_{i,i+1}$ denotes the numerical flux in the $x$-direction across the face between the nodes $(i,j)$ and $(i+1,j)$, and similarly for 
${\Phi}_{i,i-1}$, and ${\Phi}_{j, j+1} $ and ${\Phi}_{j, j-1} $ are numerical fluxes in the $y$-direction. 
It is emphasized that our focus is on the unstructured-grid scheme and therefore it is not actually implemented with $(i,j)$ data structure, which is used just for convenience. 
Clearly, the scheme can be third-order as a finite-difference scheme with FSR applied in each coordinate direction. It can be implemented in an unstructured-grid code with the following flux function: 
\begin{eqnarray}
 {\Phi}_{jk} =  \frac{1}{2} \left[   {\bf f}_L   +  {\bf f}_R  \right]  - \frac{1}{2} \hat{\bf D}_n \left[  {\bf u}({\bf w}_R) -  {\bf u}({\bf w}_L) \right],
\label{cfd_numerical_flux_2d_euler_FR}
\end{eqnarray}
where the left and right fluxes are directly reconstructed,
\begin{eqnarray}
{\bf f}_L &=&  \theta \frac{ {\bf f}_{j}  + {\bf f}_{k} }{2} + (1-\theta) \left[  {\bf f}_{j} +  \frac{1}{2} \nabla {\bf f}_j \cdot ( {\bf x}_{k} - {\bf x}_j )   \right], 
\\ [2ex]
{\bf f}_R &=&  \theta \frac{ {\bf f}_{k}  + {\bf f}_{j} }{2} + (1-\theta) \left[  {\bf f}_{k} -  \frac{1}{2} \nabla {\bf f}_k \cdot ( {\bf x}_{k} - {\bf x}_j )   \right],
\end{eqnarray}
where $\theta = 1/3$, and ${\bf w}_L$ and ${\bf w}_R$ are computed by Equations
(\ref{umuscl_of_burg_original_form}) and (\ref{umuscl_of_burg_original_form_wR}), respectively. In contrast, the third-order edge-based method requires a similar flux reconstruction but $\theta$ must be zero [\citen{nishikawa_boundary_formula:JCP2015}]. As in one dimension, the choice of $\kappa$ for the solution reconstruction is arbitrary, and we set $\kappa=1/2$ unless otherwise stated. 

Although the flux reconstruction ensures genuine third-order accuracy, it is very expensive because the flux gradients need to be computed and stored for each flux: two and three flux vectors in two and three dimensions, respectively. The FSR scheme is, therefore, not a practical scheme; but it can be made very efficient by the chain rule [\citen{Nishikawa_FakeAccuracy:2020,nishikawa_hyperbolic_ns:AIAA2015,liu_nishikawa_aiaa2016-3969}]: 
\begin{eqnarray}
 \nabla {\bf f}_j =   \left(  \frac{\partial {\bf f}}{\partial {\bf w}} \right)_{\!\! j}  \nabla {\bf w}_j    
, \quad
 \nabla {\bf f}_k =   \left(  \frac{\partial {\bf f}}{\partial {\bf w}} \right)_{\!\! k}   \nabla {\bf w}_k  ,
 \label{chain_rule_grad}
\end{eqnarray}
where the flux Jacobian is given for the Euler equations, with the unit directed-area vector $\hat{\bf n}_{jk}$ and ${\bf v}$ defined as column vectors and the notation $u_n = {\bf v} \cdot  \hat{\bf n}_{jk}$, as
\begin{eqnarray}
 \frac{\partial {\bf f}}{\partial {\bf w}}
 =
  \left[  \begin{array}{ccc} 
u_n          &  \rho  \hat{\bf n}_{jk}^t   & 0    \\ [1ex]
u_n {\bf v} & \rho ( u_n {\bf I} + {\bf v} \dyad   \hat{\bf n}_{jk} )  &  \hat{\bf n}_{jk}    \\ [1ex]
u_n {\bf v}^2 /2&  \rho(H  \hat{\bf n}_{jk}^t + u_n {\bf v}^t ) & \gamma u_n /(\gamma-1)
              \end{array} \right].
\end{eqnarray}
In this way, we do not need to compute and store multiple fluxes and their gradients in two and three dimensions; it dramatically reduces the cost of flux reconstruction.
The U-MUSCL scheme with the efficient flux reconstruction is referred to as the third-order chain-rule-based flux-solution-reconstruction scheme (CFSR3) [\citen{Nishikawa_FSR_arXiv:2020}]. Extensions to higher-order accuracy are possible [\citen{Nishikawa_FSR_arXiv:2020}], which will be further 
discussed in a subsequent paper.

The CFSR scheme was first introduced in Ref.[\citen{Nishikawa_FSR_arXiv:2020}]. As an alternative, another efficient scheme called QFSR was also introduced in Ref.[\citen{Nishikawa_FSR_arXiv:2020}], where the flux is reconstructed in the form of a quadratic function of ${\bf w}$. It can also be used to achieve genuine third-order accuracy; however, we do not discuss QFSR here because CFSR is simpler. In the next section, we discuss a newly discovered technique to achieve genuine third-order accuracy without flux reconstruction on regular simplex-element grids. 

\section{A Special Third-Order U-MUSCL Scheme with  $\kappa=1/2$ (U-MUSCL-SSQ)}
\label{umuscl-s}

In this section, we present a unique way of achieving genuine third-order accuracy with the U-MUSCL reconstruction scheme based on vanishing residuals. It is similar to the third-order edge-based method, but different in that flux reconstruction is not necessary and it is third-order accurate only on regular simplex-element grids whereas the third-order edge-based method is third-order on arbitrary simplex-element grids and requires flux reconstruction.

\subsection{One dimension}

Consider the typical U-MUSCL scheme (\ref{fv_form_1D_FD}) and its truncation error (\ref{umuscl_1d_te_2nd}), it is noted that $\kappa=1/2$ reconstructs a quadratic solution exactly at a face from point-valued solutions stored at nodes [\citen{Nishikawa_3rdQUICK:2020}] (see also the next section). Then, the U-MUSCL scheme (\ref{fv_form_1D_FD}) is essentially equivalent to the second-order central finite-difference scheme; indeed, for $\kappa=1/2$, the truncation error (\ref{umuscl_1d_te_2nd}) becomes 
\begin{eqnarray} 
\frac{d {u}_i }{d t} +  f_x = s_i -   \frac{1}{24} \left[     f_{uuu} (u_x)^3 + 3  f_{uu} u_x  u_{xx}   +  f_u  u_{xxx}
\right] h^2 + O (h^3).
\label{umuscl_1d_te_2nd_3}
\end{eqnarray}
which simplifies, by $f_{xxx} =  f_{uuu} (u_x)^3 +  3 f_{uu} u_x  u_{xx}   +   f_u  u_{xxx}$, to
\begin{eqnarray} 
\frac{d {u}_i }{d t} +  f_x = s_i - \frac{1}{24} f_{xxx} h^2 + O (h^3).
\label{umuscl_1d_te_2nd_4}
\end{eqnarray} 
This is clearly a second-order scheme. However, if we were able to discretize the time derivative and source terms in such a way 
that the second-order truncation error is factored: 
 \begin{eqnarray} 
\frac{d {u}_i }{d t} +  f_x = s_i - \frac{1}{24} \left( \frac{d {u} }{d t}  + f_x - s \right)_{xx} h^2 + O (h^3),
\label{umuscl_1d_te_2nd_5}
\end{eqnarray}
then, since $\frac{d {u}_i }{d t} +  f_x - s_i =0$ for an exact solution, we would obtain a third-order scheme:  
 \begin{eqnarray} 
\frac{d {u}_i }{d t} +  f_x = s_i + O (h^3).
\label{umuscl_1d_te_2nd_6}
\end{eqnarray}
We have discovered such source discretization formula, which is given by
\begin{eqnarray}
s_i = 
s(x_j) + \frac{\kappa_s}{4} \left[ 
\,  s(x_{i+1}) -2  s(x_i) +  s(x_{i-1}) \, \right], 
\end{eqnarray}
where
\begin{eqnarray}
\kappa_s = \frac{1}{6},
\end{eqnarray}
and similarly for the time derivative term. It can be expressed with the U-MUSCL reconstruction formula as 
 \begin{eqnarray} 
s_i = \frac{ s_{i+1/2} + s_{i-1/2} }{2},
\label{source_special}
\end{eqnarray}
where
 \begin{eqnarray}
 s_{i-1/2} &=& \kappa_s \frac{ s(x_{i-1}) +  s(x_i) }{2} 
 + (1 - \kappa_s) \left[ s(x_i) - s_x(x_i) \frac{h}{2} \right], 
 \\ [2ex]
 s_{i+1/2} &=&  \kappa_s \frac{ s(x_{i+1}) +  s(x_i) }{2} 
 + (1 - \kappa_s) \left[ s(x_i) + s_x(x_i) \frac{h}{2} \right],
\\ [2ex]
  s_x(x_i) &=& \frac{  s(x_{i+1}) - s(x_{i-1} )   }{2 h}. 
\end{eqnarray}
Note that the above source discretization is compact with a three-point stencil. This scheme does not seem to have been known; it is here referred to as U-MUSCL-SSQ. Remarkably, the error-cancellation property extends to {\color{black} triangular and tetrahedral} grids as we will discuss in the next section. 

\subsection{Multi-dimensions}
 
  Consider, for simplicity, a scalar nonlinear equation with a source term, $\partial_x f(u) + \partial_y g(u) = s$, and its edge-based discretization with an extended special source discretization:
\begin{eqnarray}
Res_j = \frac{1}{V_j} \sum_{k \in \{ k_j\} }  {\phi}_{jk} A_{jk} 
- \tilde{s}_j, 
\quad
\tilde{s}_j = \frac{1}{V_j} \sum_{k \in \{ k_j\} }  {\psi}_{jk} V_{jk}, 
\label{residual_imuscl_s}
\end{eqnarray}
where, with $f_n = (f,g) \cdot \hat{\bf n}_{jk}$ and $D_n = \partial f_n / \partial u$, 
\begin{eqnarray}
 {\phi}_{jk} =  \frac{1}{2} \left[   {f_n}(u_L)   +  {f_n}(u_R)  \right]  - \frac{1}{2} \hat{D}_n ( u_R - {u}_L ), 
\end{eqnarray}
 \begin{eqnarray}
 V_{jk} = \frac{1}{4} (  {\bf x}_k - {\bf x}_j  )\cdot {\bf n}_{jk} 
 ,
 \quad
  {\psi}_{jk} = 
   \kappa_s \frac{ {s}_j + {s}_k }{2} 
 + (1 - \kappa_s) \left[ {s}_j +
 \frac{1}{2} \overline{\nabla} {s}_j 
 \cdot (  {\bf x}_k - {\bf x}_j  ) \right]
.
\end{eqnarray}
Note that the gradients are computed by a linear LSQ method, which is sufficient for an exact quadratic extrapolation on regular grids. 

{\color{black} For regular quadrilateral grids, this scheme reduces to the one-dimensional scheme applied in each coordinate direction and is second-order accurate in general. 
To see this, we expand the scheme on a Cartesian grid with spacing $h$,
\begin{eqnarray} 
\frac{d {u}_i }{d t} +  f_x + g_y = s_i 
- \frac{1}{24} 
\left[
   \left(   f_x - \frac{s}{2}   \right)_{xx}  
+  \left(   g_y - \frac{s}{2}  \right)_{yy}
\right] h^2 + O (h^3).
\label{umuscl_1d_te_2nd_5_2d}
\end{eqnarray}
The second-order errors do not vanish for $\partial_x f(u) + \partial_y g(u) = s$
unless $f_x - \frac{s}{2} = 0$ and $g_y - \frac{s}{2}$ = 0. }

  \begin{figure}[t] 
    \centering
      \begin{subfigure}[t]{0.4\textwidth}
  \includegraphics[width=0.99\textwidth,trim=0 0 0 0 ,clip]{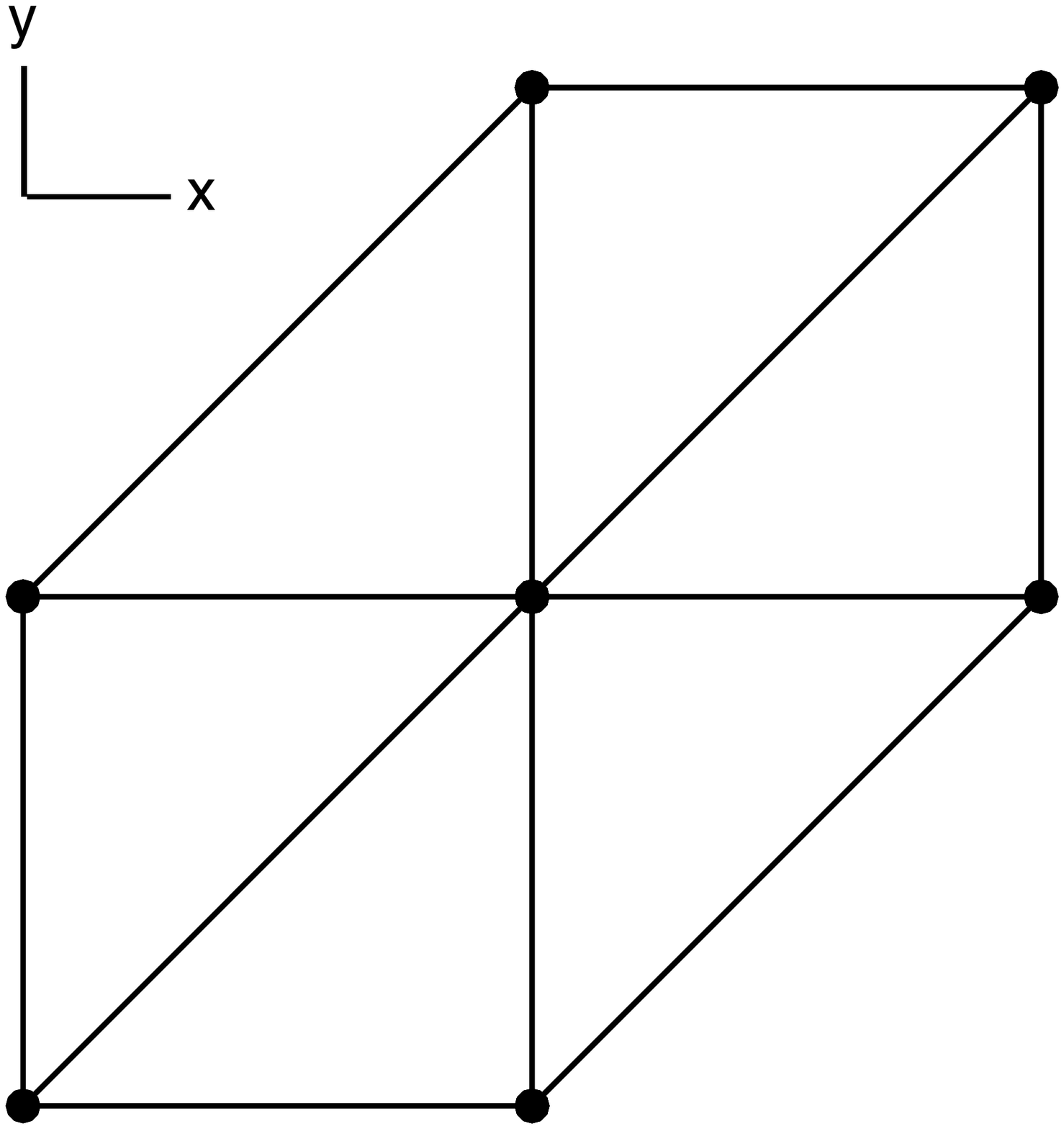}
  \caption[]{Regular triangular grid.}
 \label{fig:grid_reg_tria}
      \end{subfigure} 
      \begin{subfigure}[t]{0.4\textwidth}
  \includegraphics[width=0.99\textwidth,trim=0 0 0 0 ,clip]{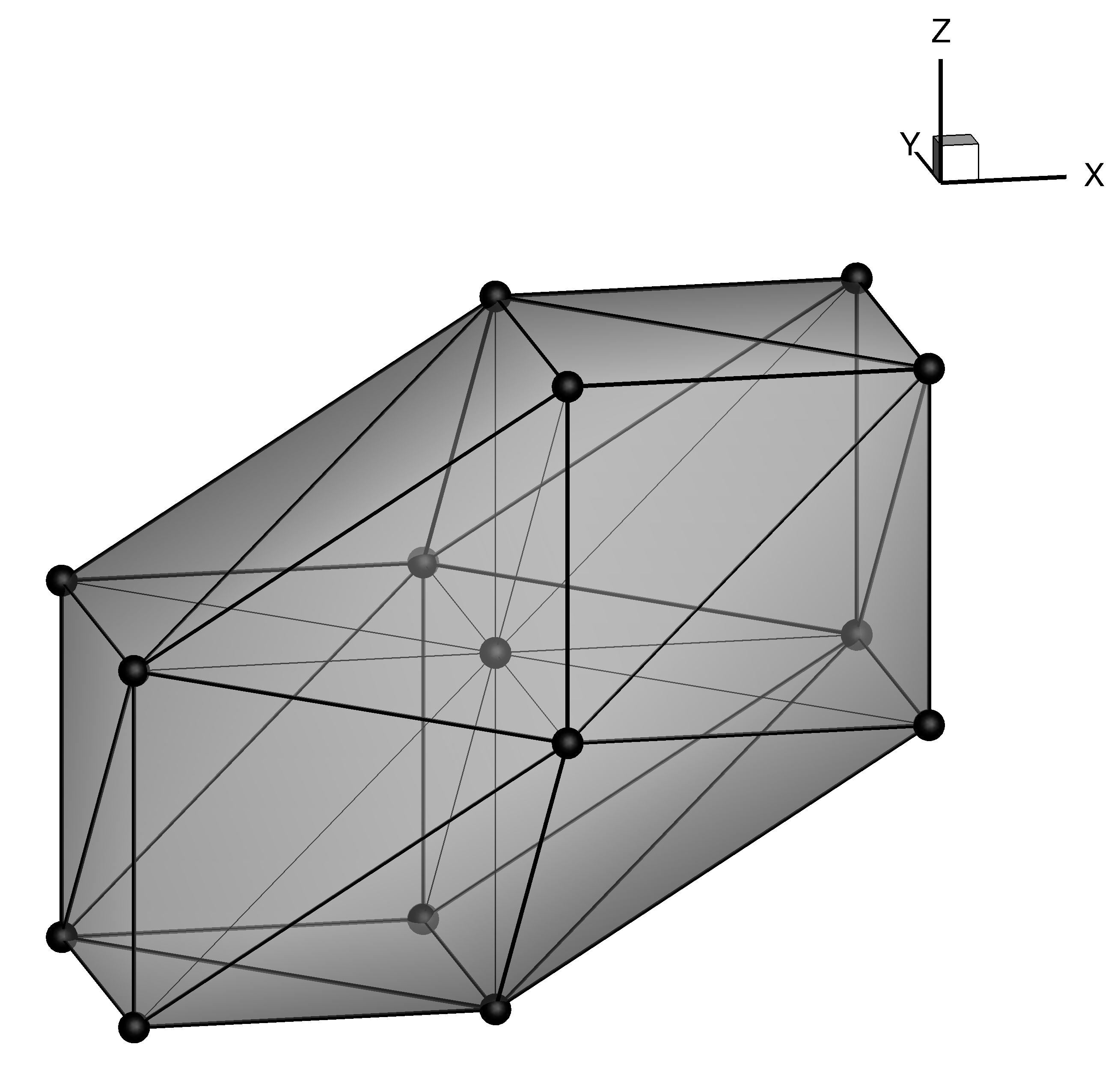}
  \caption[]{Regular tetrahedral grid.}
 \label{fig:grid_reg_tet}
      \end{subfigure}
      \caption{
\label{fig:regular_grids}%
Regular triangular and tetrahedral grids.} 
\end{figure}

On the other hand, this scheme achieves third-order accuracy for regular triangular grids with a particular value of $\kappa_s$. To see this, we substitute a smooth solution into the residual and expand it on 
a regular quadrilateral of spacing $h$ split into two right triangles as shown in Figure \ref{fig:grid_reg_tria} and obtain, for $\kappa=1/2$,
\begin{eqnarray}
Res_j &=& 
\partial_x f + \partial_y g - s 
+ \frac{1}{24} \left[ 
  \partial_{xx} (\partial_x f + \partial_y g) 
- \partial_{xy} (\partial_x f + \partial_y g) 
+ \partial_{yy} (\partial_x f + \partial_y g) 
]
\right] h^2
\nonumber \\ [2ex]
&-& \frac{\kappa_s}{6} \left[ 
  \partial_{xx} s 
- \partial_{xy} s 
+ \partial_{yy} s
\right] h^2
+ O(h^3), 
\end{eqnarray}
which suggests that we set $\kappa_s = 1/4$, so that we have
 \begin{eqnarray}
Res_j = 
r
+ \frac{1}{32} \left[ 
  \partial_{xx} r 
- \partial_{xy} r 
+ \partial_{yy} r
\right] h^2
+ O(h^3), \quad
r = \partial_x f + \partial_y g - s.
\end{eqnarray}
This gives the truncation error $TE_j$ when the smooth solution is an exact solution satisfying $r=0$:
  \begin{eqnarray}
TE_j = O(h^3). 
\end{eqnarray}
The result is valid for a time-dependent equation with the time derivative treated as a source term. It is valid also for a general regular triangular grid with each node having edge-connected nodes in a symmetric configuration, which can be shown in a similar manner to the proof provided in Ref.[\citen{nishikawa_liu_source_quadrature:jcp2017}]. See Appendix A for the proof. For example, for a triangular grid with equilateral triangles of side length $h$, we obtain
\begin{eqnarray}
Res_j = 
\partial_x f + \partial_y g - s 
+ \frac{1}{32} \left[ 
  \partial_{xx} (\partial_x f + \partial_y g) 
+ \partial_{yy} (\partial_x f + \partial_y g) 
\right] h^2
- \frac{\kappa_s}{8} \left[ 
  \partial_{xx} s 
+ \partial_{yy} s
\right] h^2
+ O(h^3), 
\end{eqnarray}
which becomes, again for $\kappa_s = 1/4$,
 \begin{eqnarray}
Res_j = 
r
+ \frac{1}{32} \left[ 
  \partial_{xx} r 
+ \partial_{yy} r
\right] h^2
+ O(h^3), \quad
r = \partial_x f + \partial_y g - s,
\end{eqnarray}
thus giving
  \begin{eqnarray}
TE_j = O(h^3). 
\end{eqnarray}

The special property extends to three dimensions. 
Consider a scalar conservation law $\mbox{div} {\bf f} = s$, where ${\bf f}$ is a flux vector, and its edge-based discretization on a tetrahedral grid:
\begin{eqnarray}
Res_j = \frac{1}{V_j} \sum_{k \in \{ k_j\} }  {\phi}_{jk} A_{jk} 
- \tilde{s}_j, 
\quad
\tilde{s}_j = \frac{1}{V_j} \sum_{k \in \{ k_j\} }  {\psi}_{jk} V_{jk}, 
\label{residual_imuscl_s}
\end{eqnarray}
where, with $f_n = {\bf f} \cdot \hat{\bf n}_{jk}$ and $D_n = \partial f_n / \partial u$, 
\begin{eqnarray}
 {\phi}_{jk} =  \frac{1}{2} \left[   {f_n}(u_L)   +  {f_n}(u_R)  \right]  - \frac{1}{2} \hat{D}_n ( u_R - {u}_L ), 
\end{eqnarray}
 \begin{eqnarray}
 V_{jk} = \frac{1}{6} (  {\bf x}_k - {\bf x}_j  )\cdot {\bf n}_{jk} 
 ,
 \quad
  {\psi}_{jk} = 
   \kappa_s \frac{ {s}_j + {s}_k }{2} 
 + (1 - \kappa_s) \left[ {s}_j +
 \frac{1}{2} \overline{\nabla} {s}_j 
 \cdot (  {\bf x}_k - {\bf x}_j  ) \right]
.
\end{eqnarray}
Note that these expressions are the same as those in two dimensions except the factor in the partial volume $V_{jk}$.
The directed-area vector ${\bf n}_{jk}$ is defined by the sum of dual triangular faces formed by the cell centroid, face centroids,
and the edge midpoint (see, e.g., Ref.[\citen{liu_nishikawa_aiaa2016-3969}]).
For a regular tetrahedral grid generated by subdividing a regular hexahedral grid by splitting a hexadedron into six tetrahedra, we 
obtain, for $\kappa=1/2$,
\begin{eqnarray}
Res_j &=& 
\mbox{div} \, {\bf f} - s 
+ \frac{1}{24} \left[ 
  \partial_{xx} (\mbox{div} \,{\bf f}) 
+ \partial_{yy} (\mbox{div} \,{\bf f}) 
+ \partial_{zz} (\mbox{div} \,{\bf f}) 
- \partial_{xy} (\mbox{div} \,{\bf f}) 
- \partial_{yz} (\mbox{div}\, {\bf f}) 
+ \partial_{zx} (\mbox{div} \,{\bf f}) 
]
\right] h^2
\\
&+& \frac{5\kappa_s}{36} \left[ 
  \partial_{xx} s 
+ \partial_{yy} s
+ \partial_{zz} s
- \partial_{xy} s 
- \partial_{yz} s 
+ \partial_{zx} s
\right] h^2
+ O(h^3), 
\end{eqnarray}
which suggests that we set $\kappa_s = 3/10$, so that we have
 \begin{eqnarray}
Res_j = 
r
+ \frac{1}{24} \left[ 
  \partial_{xx} r
+ \partial_{yy} r 
+ \partial_{zz} r 
- \partial_{xy} r 
- \partial_{yz} r 
+ \partial_{zx} r
\right] h^2
+ O(h^3), \quad
r = \mbox{div} \,{\bf f} - s,
\end{eqnarray}
thus giving
  \begin{eqnarray}
TE_j = O(h^3). 
\end{eqnarray}
A more general proof can be found in Appenix A. 

Unlike the third-order edge-based method, the U-MUSCL-SSQ scheme loses third-order accuracy on irregular grids for a reason mentioned in the next section. Nevertheless, as we will show later, it can be a very accurate second-order scheme compared with other second-order schemes.

\section{Exact quadratic reconstruction with $\kappa=1/2$ on Arbitrary Grids}
\label{truth_umuscl_2_exact_quadratic_kappa1o2}

As it was known (but a proof was never provided)  [\citen{nishikawa_hyperbolic_ns:AIAA2015}], the U-MUSCL reconstruction scheme by itself can be  high-order on arbitrary grids: it can be exact for quadratic functions if the nodal gradients are computed by an algorithm exact for quadratic functions, e.g., by a quadratic LSQ method. To prove this, let us consider a quadratic function $ u_{quadratic}(x,y) $ defined around a node $j$: 
\begin{eqnarray}
 u_{quadratic}({\bf x}) 
 &=&   {u}_j +  (x-x_j) \partial_x u_j + (y-y_j) \partial_y u_j +     \frac{1}{2}  (x-x_j)^2 \partial_{xx} u_j  \nonumber \\
 &+&  (x-x_j) (y-y_j) \partial_{xy} u_j +
    \frac{1}{2}  (y-y_j)^2 \partial_{yy} u_j,
\end{eqnarray}
where all the derivatives involved in the above expression are constants defined at $j$. Then, the function value at a neighbor node $k$ is exactly expressed by
\begin{eqnarray}
  u_k  = u_{quadratic}({\bf x}_k)  =  u_j +  \partial_{jk} u_j + \frac{1}{2}  \partial_{jk}^2 u_j, 
    \end{eqnarray}
where we have introduced the notation: 
\begin{eqnarray}
 \partial_{jk} \equiv ( {\bf x}_{k} - {\bf x}_j )  \cdot (\partial_x, \partial_y).
 \end{eqnarray}
Assume now that gradients are computed by a quadratic LSQ method and thus can be written exactly as
\begin{eqnarray}
\overline{\nabla} u_j = \nabla u_j , \quad 
\overline{\nabla} u_k = \nabla u_j + \partial_{jk}  \nabla u_j.
\end{eqnarray}
Then, substituting these exact expressions into the U-MUSCL reconstruction schemes, we find
\begin{eqnarray}
u_L 
&=& u_{j} + \frac{\kappa}{2} ( u_k- u_j  ) +   \frac{1-\kappa}{2}  \nabla u_j \cdot ( {\bf x}_{k} - {\bf x}_j )   \nonumber  \\ [2ex]
&=& u_{j} + \frac{\kappa}{2} (   \partial_{jk} u_j + \frac{1}{2}  \partial_{jk}^2 u_j   ) +  \frac{1-\kappa}{2}   \partial_{jk} u_j   \nonumber \\ [2ex]
&=& u_{j}  +  \frac{1}{2}   \partial_{jk} u_j + \frac{\kappa}{4}  \partial_{jk}^2 u_j  ,
\end{eqnarray}
and for $u_R$,
\begin{eqnarray}
u_R
&=& u_{k} + \frac{\kappa}{2} ( u_j- u_k  ) +   \frac{1-\kappa}{2}  \nabla u_k \cdot ( {\bf x}_{j} - {\bf x}_k )   \nonumber  \\ [2ex]
&=& (u_j +  \partial_{jk} u_j + \frac{1}{2}  \partial_{jk}^2 u_j) - \frac{\kappa}{2} (   \partial_{jk} u_j + \frac{1}{2}  \partial_{jk}^2 u_j   ) +  \frac{1-\kappa}{2} ( \nabla u_j + \partial_{jk}  \nabla u_j)  \cdot ( {\bf x}_{j} - {\bf x}_k )   \nonumber \\ [2ex]
&=& u_{j}  + \left( 1 - \frac{\kappa}{2}  -    \frac{1-\kappa}{2}   \right) \partial_{jk} u_j  
  + \left(   \frac{1}{2} - \frac{\kappa}{4} -          \frac{1-\kappa}{2}         \right) \partial_{jk}^2 u_j   \nonumber  \\ [2ex]
&=& u_{j}  +    \frac{1 }{2}  \partial_{jk} u_j    +   \frac{\kappa}{4} \partial_{jk}^2 u_j,
\end{eqnarray}
which have just proved two things. First, the jump at the face vanishes for quadratic functions on arbitrary grids for any $\kappa$:
\begin{eqnarray}
u_R -  u_L = \left( u_{j}  +  \frac{1}{2}   \partial_{jk} u_j + \frac{\kappa}{4}  \partial_{jk}^2 u_j  \right) -  \left( u_{j}  +  \frac{1}{2}   \partial_{jk} u_j + \frac{\kappa}{4}  \partial_{jk}^2 u_j  \right) =  0.
\end{eqnarray}
It is however important to note that it merely shows $u_L = u_R$ and the values $u_L$ and $u_R$ are not exact for a quadratic function. 
 Note also that if the gradients are computed by a linear LSQ method, then the jump can vanish for a quadratic function only with $\kappa=1$, which makes the jump identically zero for any function. Second, the U-MUSCL reconstruction scheme reconstructs a quadratic function exactly on an arbitrary grid with $\kappa=1/2$:
\begin{eqnarray}
u_L = u_R = u_{quadratic}\left( \frac{ {\bf x}_j + {\bf x}_k}{2} \right) =  u_{j}  +  \frac{1 }{2}  \partial_{jk} u_j  +   \frac{1}{8} \partial_{jk}^2 u_j.
\end{eqnarray}
Therefore, $\kappa=1/2$ is the only choice that satisfies both the zero jump condition and the exact quadratic reconstruction. 
Note that both of the two things discussed here are true for the reconstruction at halfway between two nodes; the same will be true even when it is applied to a cell-centered 
scheme as long as the reconstruction is performed at the midpoint of two adjacent centroids. 

It is important to note that $\kappa=1/2$ reconstructs the point-valued solution at the edge-midpoint from the point-valued solutions stored at nodes. This does not lead to a third-order scheme for two reasons. First, on a regular grid, it will result in the second-order central scheme as discussed in Section \ref{umuscl-s}.
It can be third-order accurate with the special source quadrature formula but only on simple-element grids. Remember also that a finite-volume scheme with point-valued solutions can be third-order accurate with $\kappa=1/2$ but it does not extend to higher dimensions. Second, on irregular grids, the edge-based flux quadrature formula is exact for quadratic fluxes (thus third-order) only with $\kappa=0$ and special source discretization formulas [\citen{nishikawa_liu_source_quadrature:jcp2017,nishikawa_boundary_formula:JCP2015}] and only for simplex-element grids. Nevertheless, 
numerical results indicate that $\kappa=1/2$ gives superior performance over other values as we will show later. 

Incidentally, the ability to quadratically interpolate the solution at the edge-midpoint is relevant to the third-order residual-distribution method originally developed by Caraeni and Fuchs [\citen{caraeni_fuchs:CSCC2000,Caraeni_Journal_2002}], which may be considered as a reconstruction-based third-order continuous Petrov-Galerkin finite-element method. In their method, they reconstruct a $P_2$ simplex element from a $P_1$ simplex element by interpolating the solution from two nodes of an edge to the edge-midpoint, and then calculate the element residual by integrating the governing equations of interest over the reconstructed $P_2$ element. Their interpolation formula for the edge-midpoint solution $u_m$ is equivalent to the arithmetic average of the U-MUSCL reconstruction schemes with $\kappa=1/2$:
\begin{eqnarray}
u_m =  \frac{ u_L + u_ R}{2}  =  \frac{u_j + u_k}{2} - \frac{1}{8} \left(   \overline{\nabla} {u}_k - \overline{\nabla} {u}_j  \right)  \cdot ( {\bf x}_{k} - {\bf x}_j ),
\end{eqnarray}
where the gradients are computed by a quadratic LSQ method (or by a linear method; then third-order accurate on smooth grids only). This technique was later used by Nishikawa, Rad, and Roe for developing high-order fluctuation-splitting schemes for Cauchy-Riemann, inviscid, viscous equations [\citen{nishikawa_rad_roe:AIAA2001,nishikawa_roe:ICCFD_2004,nishikawa:VKI_LS_2005,nishikawa_roe:ICCFD_2006,nishikawa_multigrid_cr:IJNMF}].

\section{U-MUSCL in Cell-Centered Methods}
\label{truth_umuscl_ccfv}

The U-MUSCL reconstruction scheme is very simple to implement and often employed in cell-centered finite-volume codes also. 
However, it loses second-order accuracy when the face centroid is not located exactly halfway between two adjacent centroids. Ref.[\citen{nishikawa_LP_UMUSCL:JCP2020}] demonstrates that the loss of second-order accuracy results in serious accuracy deterioration on unstructured grids, and proposes a suitable modification. A similar modification has been proposed also in Ref.[\citen{Tamaki_PhDThesis2018}] for preserving second-order accuracy on adaptive Cartesian grids. Burg mentioned this problem in his original paper [\citen{burg_umuscl:AIAA2005-4999}], but no modification techniques were proposed. 

To illuminate this problem, we consider the U-MUSCL reconstruction scheme in the form shown previously:
\begin{eqnarray}
{u}_L = \kappa  \frac{{u}_j + {u}_k}{2} +   (1-\kappa) \left[ u_j  + \overline{\nabla}{u}_j  \cdot ( {\bf x}_{f_c} -  {\bf x}_j ) \right],
 \label{umuscl_L2}
\end{eqnarray}
where ${\bf x}_{f_c}$ denotes a face centroid position. It is clear in this form that the U-MUSCL scheme is a linear combination of 
the linear interpolation (the first term) and the linear extrapolation (the second term). Then, we see that the linear interpolation will not be
exact for linear functions on unstructured grids if the face centroid does not coincide with the midpoint of a line segment connecting the centroids of the two adjacent cells, $j$ and $k$. As a result, second-order accuracy will be lost. As suggested in Ref.[\citen{nishikawa_LP_UMUSCL:JCP2020}], one way to recover second-order accuracy is to replace $ {u}_k$ by a linearly-extrapolated value $u_p$ at a location ${\bf x}_p$ such that the face centroid is the midpoint of the centroid of $j$ and the point ${\bf x}_p$ (i.e., ${\bf x}_p = 
2 {\bf x}_{f_c}-{\bf x}_j)$):  
 \begin{eqnarray}
{u}_L = \kappa  \frac{{u}_j + {u}_p}{2} +   (1-\kappa) \left[ u_j  + {\nabla}{u}_j  \cdot ( {\bf x}_{f_c} -  {\bf x}_j ) \right],
 \label{umuscl_L2_fixed}
\end{eqnarray}
where
\begin{eqnarray}
 u_p =  u_k  +  {\nabla} u_k  \cdot ( 
 2 {\bf x}_{f_c}-{\bf x}_j -  {\bf x}_k ).
\end{eqnarray}
This modification makes the reconstruction exact for linear functions on arbitrary grids. A similar modification needs to be applied to the right 
state $u_R$ also if necessary. See Ref.[\citen{nishikawa_LP_UMUSCL:JCP2020}] for further details and its impact on accuracy in unstructured grids.
  
A cell-centered finite-volume method with the U-MUSCL reconstruction scheme can achieve third-order accuracy in the same way as the edge-based method on regular quadrilateral grids. Therefore, as we have concluded in Section \ref{umuscl_remarks}, all the U-MUSCL-based cell-centered finite-volume codes are second-order accurate at best for any value of $\kappa$ for the nonlinear Euler equations, and it can be third-order accurate in point-valued solutions (not cell averages) with $\kappa=1/3$ on 
regular grids only for linear problems. 

One can achieve third-order accuracy on regular quadrilateral grids with CFSR (not with U-MUSCL-SSQ). It is emphasized again that third-order accuracy will be obtained in the point-valued solution, not in the cell-averaged solution. Moreover, time derivatives and source terms must be evaluated at the cell center (not integrated). Hence, the scheme must be treated as a finite-difference scheme with point-valued solutions, rather than a finite-volume scheme. It is also pointed out that third-order accuracy remains to be demonstrated for regular triangular grids, in which the cell-centered stencil is not the same everywhere (there are two different configurations) and therefore it is not clear if any error cancellation occurs. Studies on cell-centered schemes are left as future work.

\section{Results}
\label{results}


\subsection{Numerical tests in one dimensions}

We present numerical results for third-order schemes in one dimension. The finite-volume scheme with cell averaged solutions as discussed in Section \ref{oned_FVC} will be referred to as FVC. Similarly, FVP denotes the finite-volume scheme with point valued solutions in Section \ref{oned_FVP}. The U-MUSCL scheme in the form (\ref{fv_form_1D_FD}) is simply referred to as U-MUSCL. These schemes are summarized along with FSR and U-MUSCL-SSQ in Table \ref{Tab.classificaiton_1D}.

\begin{table}[t]
\ra{1.1}
\begin{center}
{\tabulinesep=0.5mm
\begin{tabu}{lccrrrrr}\hline\hline 
\multicolumn{1}{r}{  \multirow{2}{*}{Type}      }        &
\multicolumn{1}{r}{  \multirow{2}{*}{$\kappa$}      }        &
\multicolumn{1}{r}{  \multirow{2}{*}{$\theta$}      }        &
\multicolumn{1}{r}{ \multirow{2}{*}{Solution}       } &
\multicolumn{1}{r}{ \multirow{2}{*}{Reconstruction} } &
\multicolumn{1}{r}{ \multirow{2}{*}{$s$ and $\frac{\partial u}{\partial t}$} } &
\multicolumn{2}{r}{ \multirow{1}{*}{ Order of Accuracy} }      
\\
 \cmidrule(lr){7-8}  
&
&
&
&
&
&
\multicolumn{1}{l}{Linear}     &
\multicolumn{1}{l}{Nonlinear}  
  \\ \hline 
U-MUSCL: FVC & $ 1/3$  & -- & cell average & solution   & cell average  & $\boldsymbol {\color{red} O(h^3)} $ & $\boldsymbol {\color{red} O(h^3)} $  \\
U-MUSCL: FVP&  $1/2$  & --   & point value  & solution   &  cell average  & $\boldsymbol {\color{red} O(h^3)} $ & $\boldsymbol {\color{red} O(h^3)} $  \\
U-MUSCL: FD &  $1/3$  & --   & point value  & solution   &  point value   & $\boldsymbol {\color{red} O(h^3)} $ & $O(h^2)$  \\
FSR & $\mathbb{R}$ & $1/3$  & point value  & flux and solution  &  point value   & $\boldsymbol {\color{red} O(h^3)} $ & $\boldsymbol {\color{red} O(h^3)} $  \\
U-MUSCL-SSQ &     $ 1/2$   & --  & point value  & solution   &  quadrature (\ref{source_special})  & $\boldsymbol {\color{red} O(h^3)} $ & $\boldsymbol {\color{red} O(h^3)} $ \\
  \hline  \hline
\end{tabu}
}

\caption{Third-order schemes in the form of the U-MUSCL scheme in one dimension. $\mathbb{R}$ indicates any real value. Orders of accuracy are shown for linear and nonlinear equations. }

\label{Tab.classificaiton_1D}
\end{center}
\end{table}

\subsubsection{Accuracy verification with steady problems}

We first consider a steady state problem for Burgers' equation:
\begin{eqnarray}
 \partial_x f = s(x),
\end{eqnarray}
where $f=u^2/2$ and $s(x) =  a C^2 \exp(2 a x) $, so that the exact solution is given by 
\begin{eqnarray}
u(x)  = C \exp(a x),
\end{eqnarray}
where we set $C=1.57$ and $a=1.23$. We solve the problem over a unit domain in a series of uniform grids with $n$ points, where $n = 16, 32, 64,128$, which are taken to be the cell centers. The uniform grid spacing is denoted by $h = 1/(n-1)$. To focus on the accuracy in the interior, we specify the exact solution, point value 
\begin{eqnarray}
{u}_i =   C \exp(a x_i),
\end{eqnarray}
or cell average, 
\begin{eqnarray}
\overline{u}_i =    \frac{1}{h} \int_{x_i - h/2}^{x_i+h/2} u \, dx  =  \frac{C}{a h} \left[   \exp(    a(x_i + h/2 )  )  -   \exp(    a(x_i - h/2 )  )  \right],
\end{eqnarray}
 at the left boundary point and its neighbor, and similarly for the right boundary. The forcing term $s(x)$ is discretized either by a point evaluation $s(x_i)$
\begin{eqnarray}
 {s}_i = s(x_i) = a C^2 \exp(2 a x_i),
\end{eqnarray}
or a cell average:
\begin{eqnarray}
\overline{s}_i =   \frac{1}{h} \int_{x_i - h/2}^{x_i+h/2} s \, dx  =   \frac{C^2}{2 h} \left[   \exp(   2 a(x_i + h/2 )  )  -   \exp(  2  a(x_i - h/2 )  )  \right].
\end{eqnarray}

First, three schemes are considered: FVC(1/3), FVP(1/2), and FSR(1/3). Note again that FVP(1/2) is the scheme Burg used for accuracy verification of U-MUSCL for a one-dimensional problem in his paper [\citen{burg_umuscl:AIAA2005-4999}]; it is the QUICK scheme [\citen{Nishikawa_3rdQUICK:2020}]. We compute steady state solutions by iterating with the explicit SSP RK3 time-stepping scheme [\citen{SSP:SIAMReview2001}] until the residual is reduced by seven orders of magnitude.

  \begin{figure}[t] 
    \centering
      \begin{subfigure}[t]{0.24\textwidth}
  \includegraphics[width=0.99\textwidth,trim=0 0 0 0 ,clip]{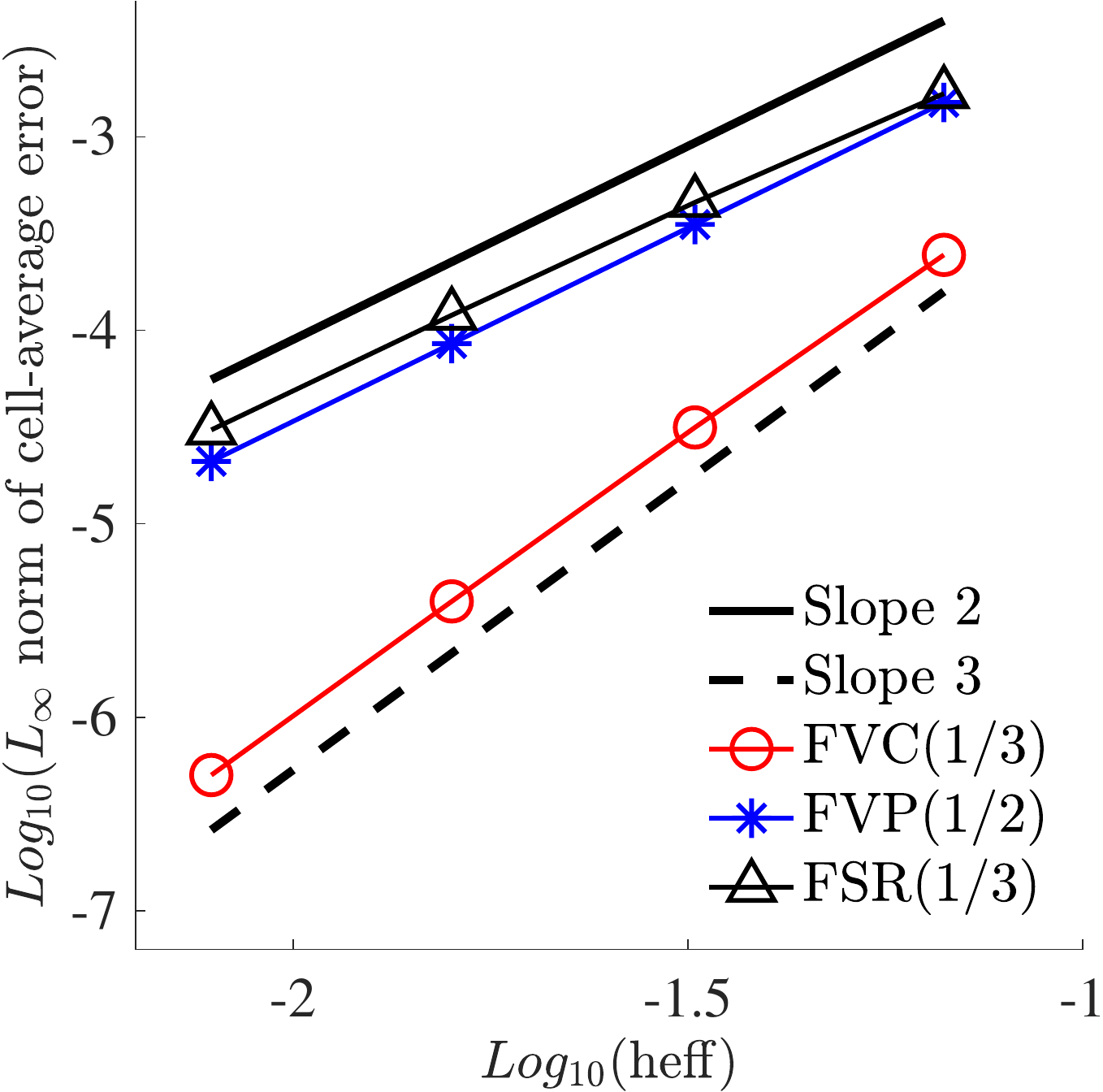}
  \caption[]{Cell average error for $\overline{s}_i$}
 \label{fig:oned_steady_cc}
      \end{subfigure} 
      \begin{subfigure}[t]{0.24\textwidth}
  \includegraphics[width=0.99\textwidth,trim=0 0 0 0 ,clip]{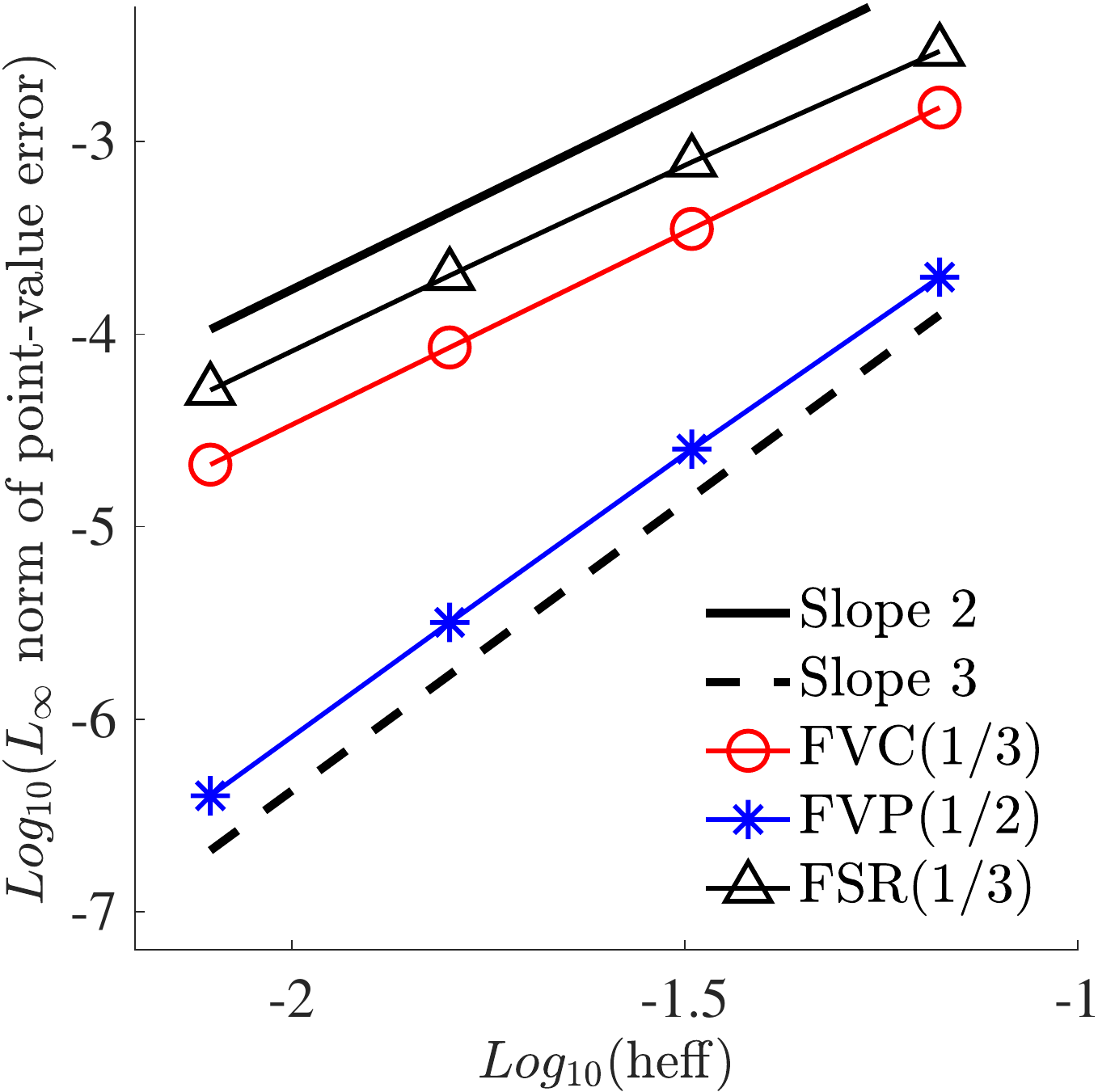}
  \caption[]{Point value error for $\overline{s}_i$}
 \label{fig:oned_steady_pc}
      \end{subfigure}
      \begin{subfigure}[t]{0.24\textwidth}
  \includegraphics[width=0.99\textwidth,trim=0 0 0 0 ,clip]{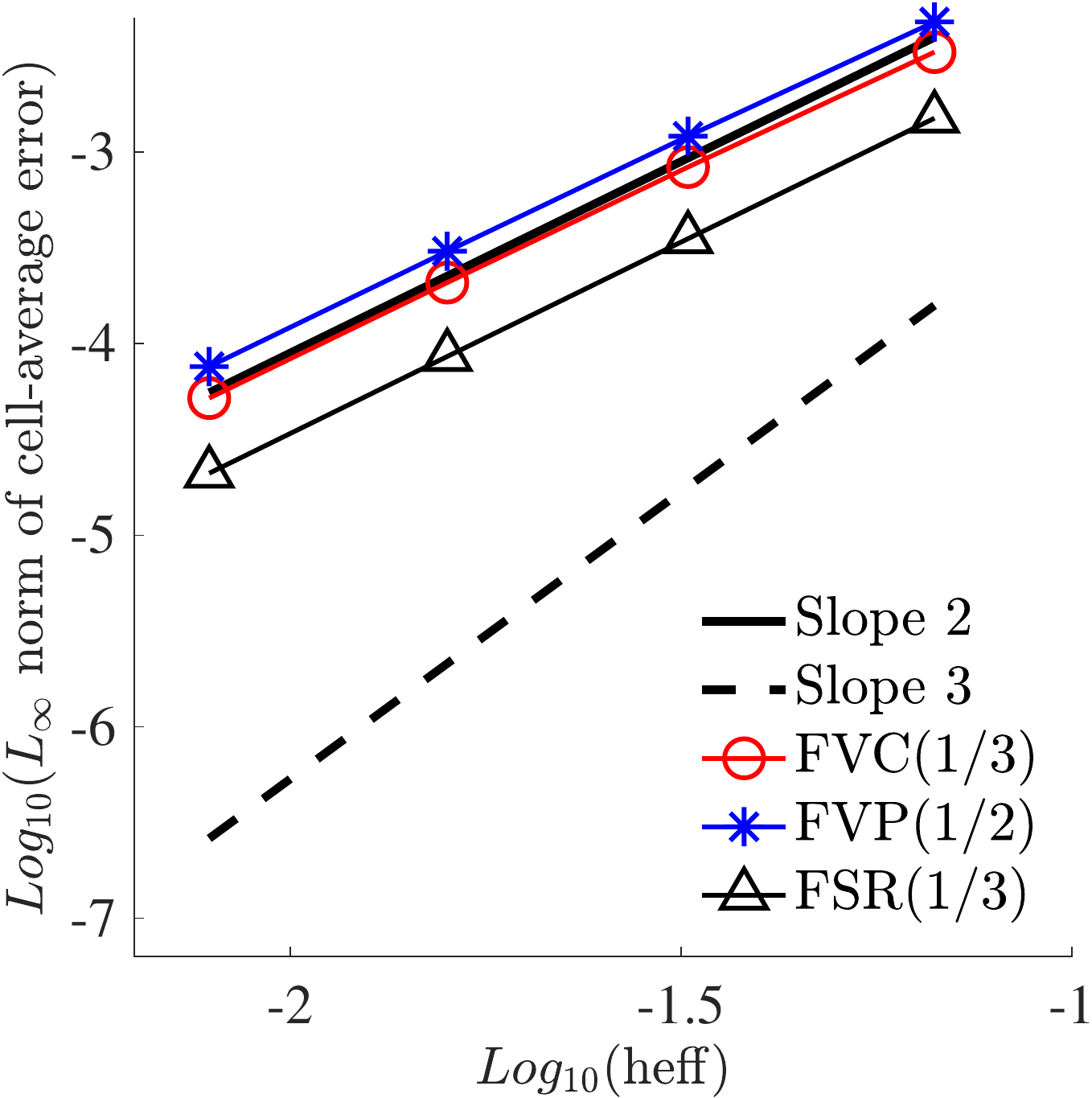}
  \caption[]{Cell average error for $ {s}_i$}
 \label{fig:oned_steady_cp}
      \end{subfigure} 
      \begin{subfigure}[t]{0.24\textwidth}
  \includegraphics[width=0.99\textwidth,trim=0 0 0 0 ,clip]{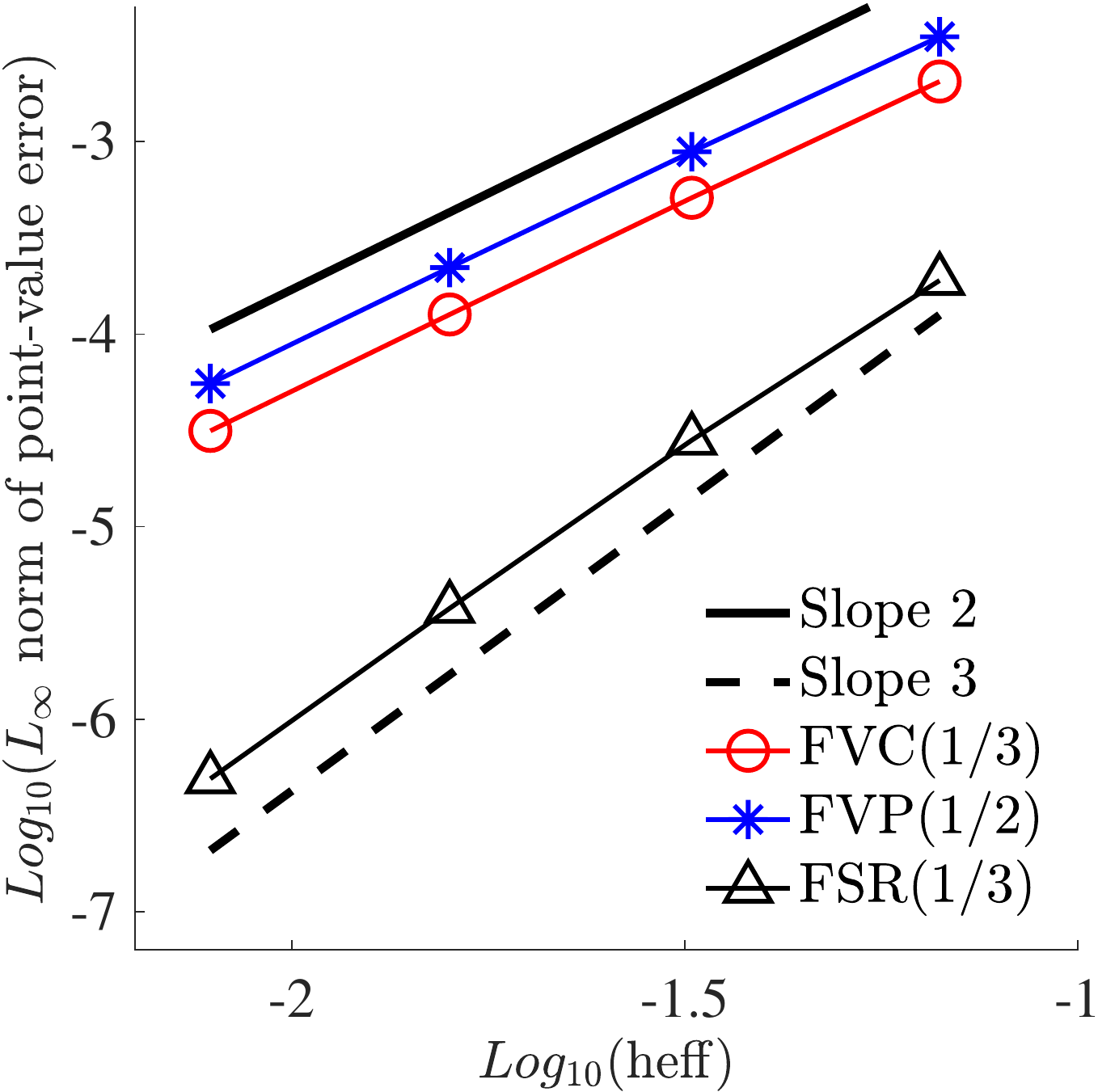}
  \caption[]{Point value error for ${\overline{}s}_i$}
 \label{fig:oned_steady_pp}
      \end{subfigure}
      \caption{
\label{fig:oned_steady}%
Error convergence results for the steady Burgers equation: $ \partial_x f = s$, where $f=u^2/2$, with the cell-averaged source term in (a) and (b), and with the point-value source term in (c) and (d)} 
\end{figure}

Results with the cell-averaged source are shown in Figures \ref{fig:oned_steady_cc} and \ref{fig:oned_steady_pc}. In Figure \ref{fig:oned_steady_cc}, the discretization errors computed against the exact cell-averaged solution are presented. As expected, FVC(1/3) is third-order accurate (MUSCL); and others are second-order 
accurate. Figure \ref{fig:oned_steady_pc} shows the discretization errors computed against the exact point-valued solution at the cell center. Here, as expected, only FVP(1/2) gives third-order accuracy (QUICK). Results with the point-valued source are shown in Figures \ref{fig:oned_steady_cp} and \ref{fig:oned_steady_pp}. In this case, FSR(1/3) is third-order accurate in the point-valued solution and others are second-order accurate as expected.  
 
  \begin{figure}[htbp!]
    \centering
      \begin{subfigure}[t]{0.45\textwidth}
  \includegraphics[width=0.85\textwidth,trim=0 0 0 0 ,clip]{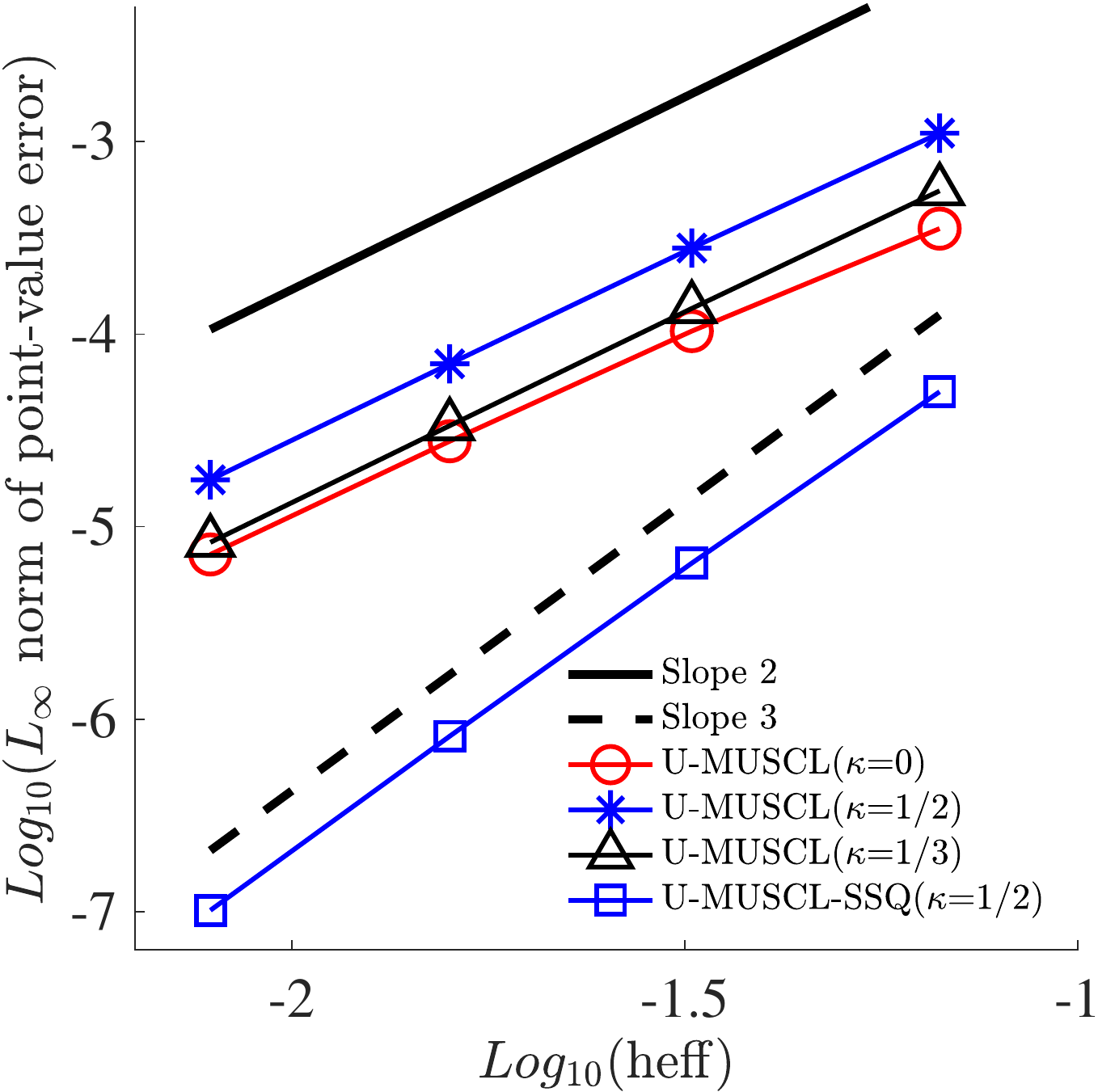}
  \caption[]{Burgers' equation.}
\label{fig:oned_steady_c_nonlinear}%
      \end{subfigure}
      \begin{subfigure}[t]{0.45\textwidth}
  \includegraphics[width=0.85\textwidth,trim=0 0 0 0, clip]{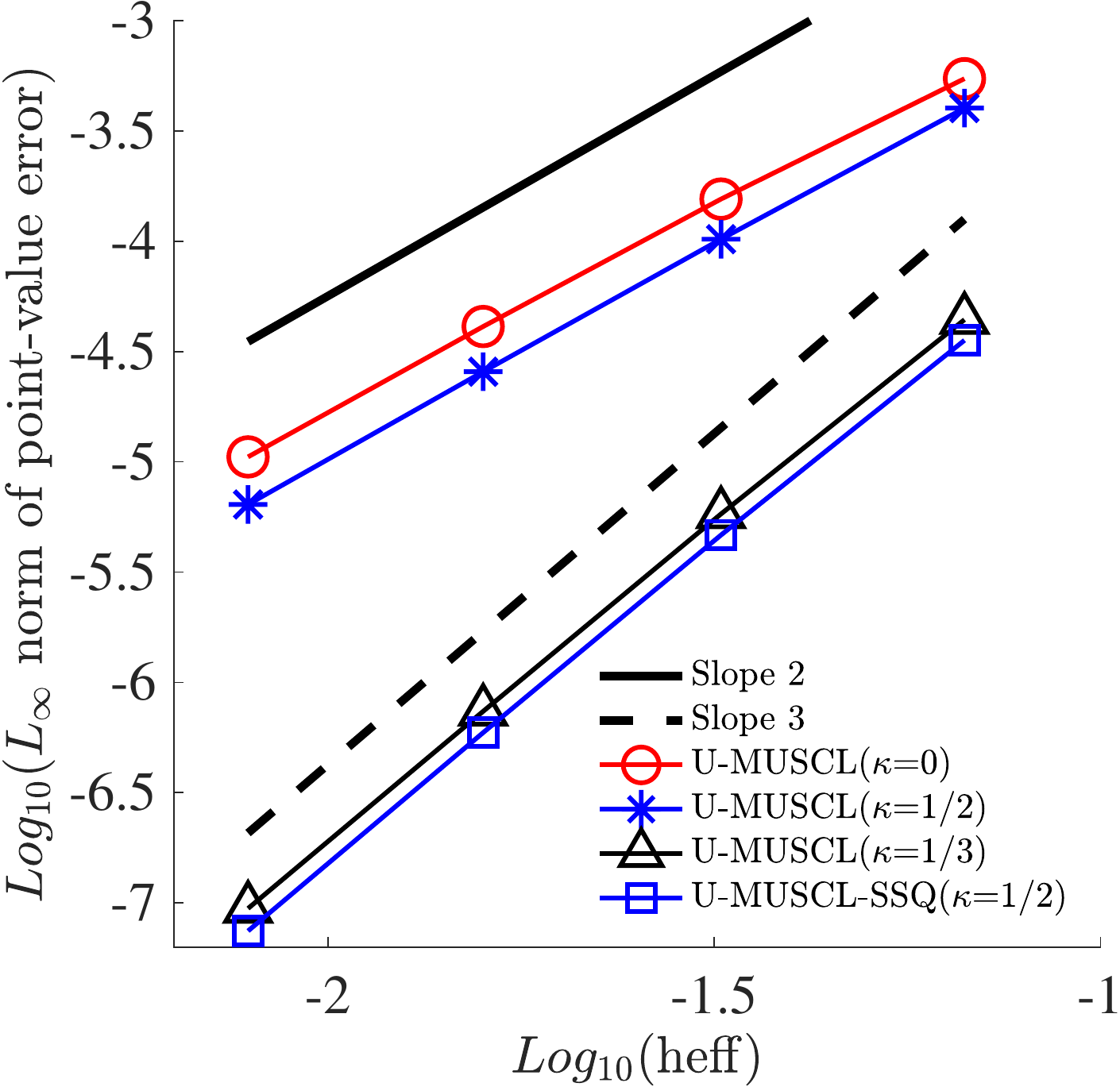}
  \caption[]{Linear advection}
\label{fig:oned_steady_p_linear}%
      \end{subfigure}
      \caption{
\label{fig:oned_steady_UMSUCL}%
Error convergence results of the U-MUSCL scheme (\ref{fv_form_1D_FD}) with various values of $\kappa$ for steady Burger's and linear advection equations.} 
\end{figure}

We then tested the U-MUSCL scheme (\ref{fv_form_1D_FD}) without flux reconstruction, which is the form typically used in practical solvers, for $\kappa=0$, $1/3$, $1/2$, and also the U-MUSCL-SSQ scheme, which is with $\kappa=1/2$ and the special source term (\ref{source_special}).  To demonstrate third-order accuracy of U-MUSCL for linear equations, we performed the test also for a linear advection equation with $f=u$. Results are shown in Figure \ref{fig:oned_steady_UMSUCL}. As expected, the U-MUSCL scheme with $\kappa=1/3$ (indicated by black triangles) is third-order accurate for the linear equation but second-order accurate for the nonlinear Burgers equation. On the other hand, the U-MUSCL-SSQ scheme is third-order accurate for both Burgers' and linear advection equations. 


  \begin{figure}[htbp!]
    \centering
      \begin{subfigure}[t]{0.45\textwidth}
  \includegraphics[width=0.85\textwidth,trim=0 0 0 0 ,clip]{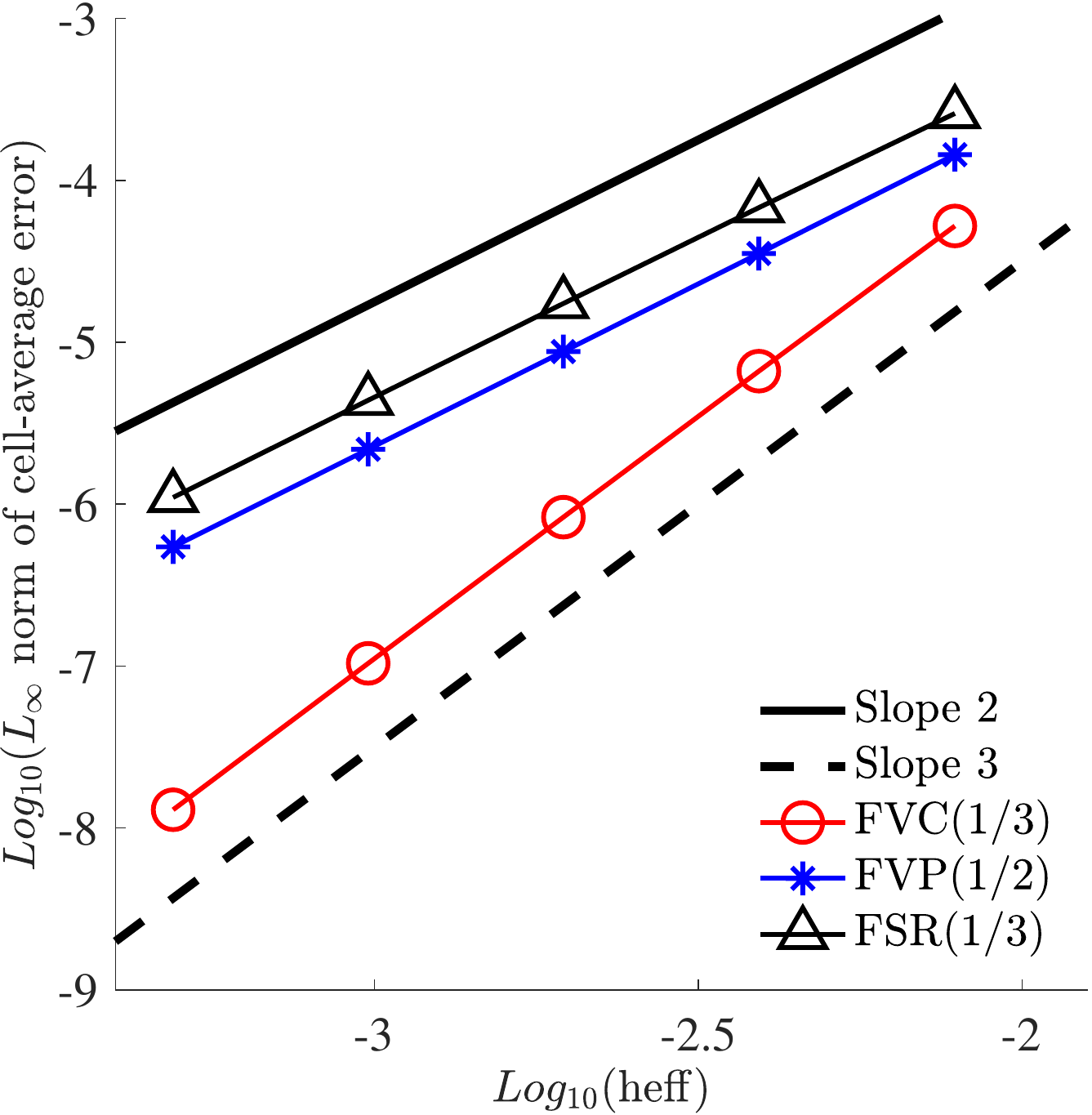}
  \caption[]{Cell average error.}
\label{fig:oned_unsteady_c}%
      \end{subfigure}
      \begin{subfigure}[t]{0.45\textwidth}
  \includegraphics[width=0.85\textwidth,trim=0 0 0 0, clip]{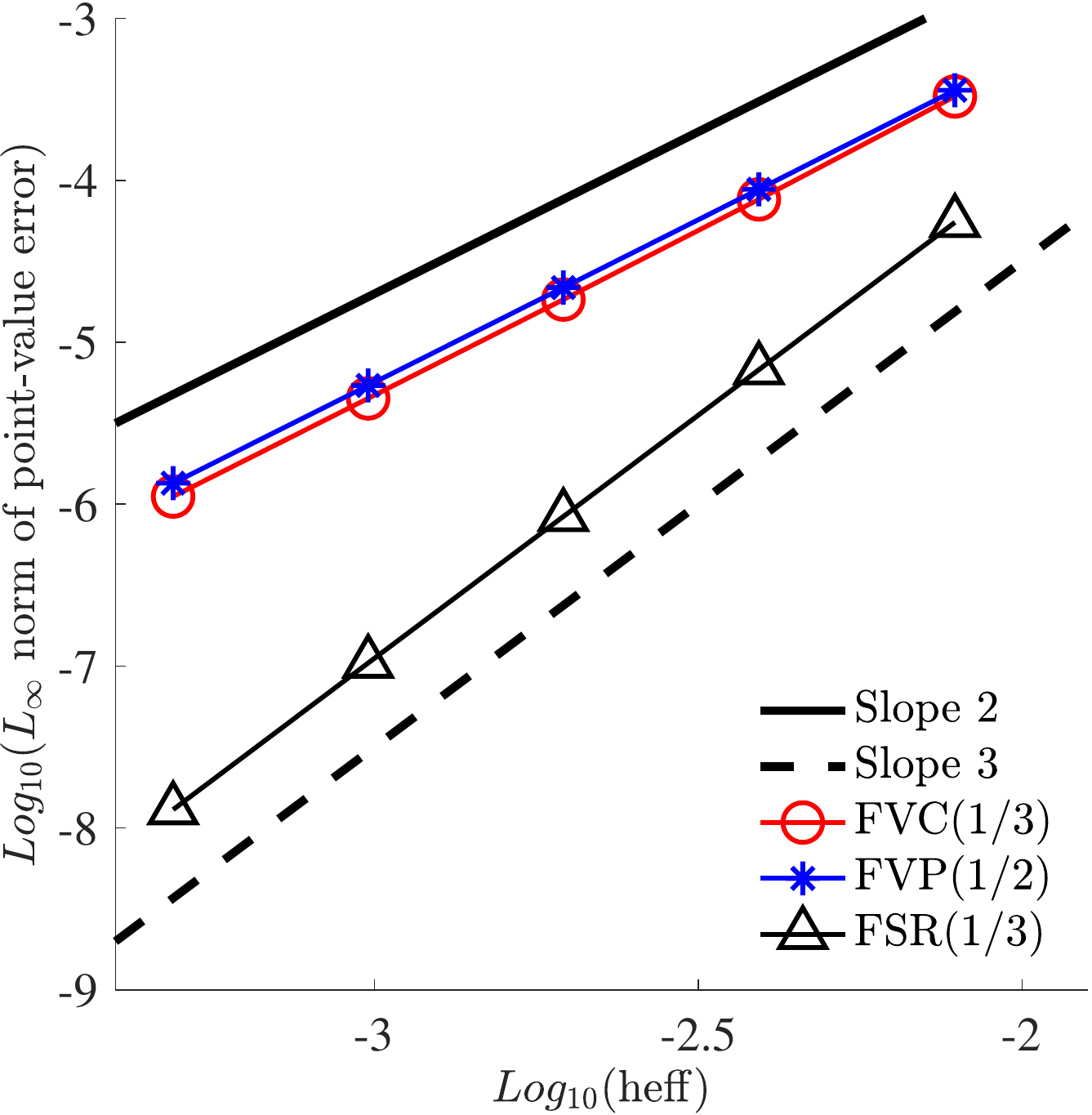}
  \caption[]{Point value error.}
\label{fig:oned_unsteady_p}%
      \end{subfigure}
      \caption{
\label{fig:oned_unsteady}%
Error convergence results for the unsteady Burgers equation: $\partial_t u +  \partial_x f = 0$, where $f=u^2/2$.} 
\end{figure}

 %

\subsubsection{Accuracy verification with unsteady problems}

Next, we consider an unsteady problem for Burgers' equation in the domain $x\in[0,1]$ for a series of uniform grids with $2^{m+6}$-$1$ cells, where $m=1,2,3,4,5$. The initial solution is set to be $\sin(2 \pi x)$. The time integration is performed with the explicit SSP RK3 scheme for 800 time steps with $\Delta t = 0.0001$ to the final time $t_{f} = 0.08$. The left and right boundaries are made periodic; thus there is no boundary. The exact solution (a moment before a shock is formed) is computed by iteratively solving $u_i = \sin(  2 \pi (x_i - u_i {t_f} ) )$ at the cell center $x=x_i$ and the exact cell average is computed by a three-point Gaussian quadrature formula with the exact cell center solution and the exact solutions at the left and right faces of the cell.

  \begin{figure}[htbp!]
    \centering
      \begin{subfigure}[t]{0.45\textwidth}
  \includegraphics[width=0.85\textwidth,trim=0 0 0 0 ,clip]{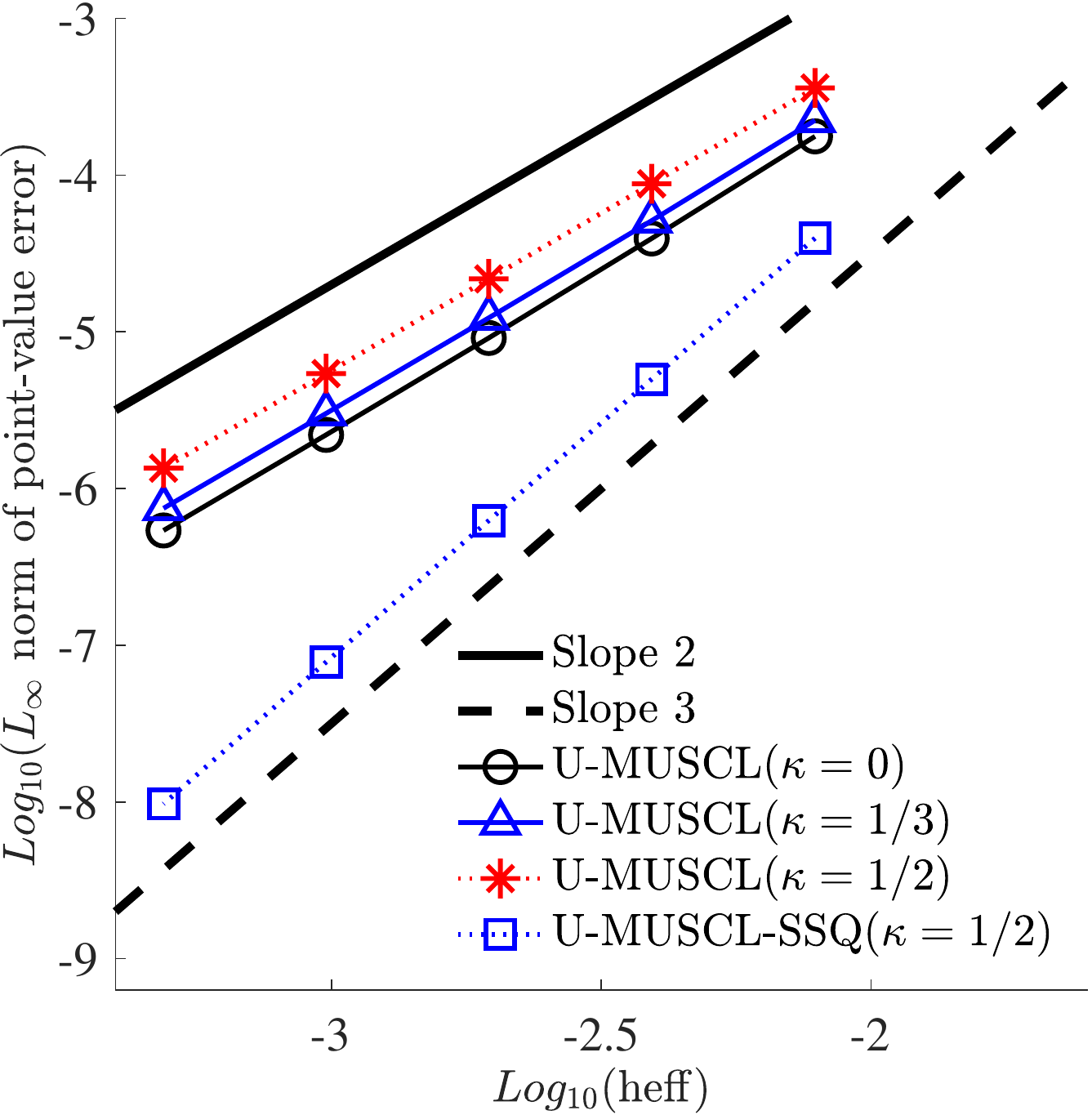}
  \caption[]{Burgers' equation.}
\label{fig:oned_unsteady_c_nonlinear}%
      \end{subfigure}
      \begin{subfigure}[t]{0.45\textwidth}
  \includegraphics[width=0.85\textwidth,trim=0 0 0 0, clip]{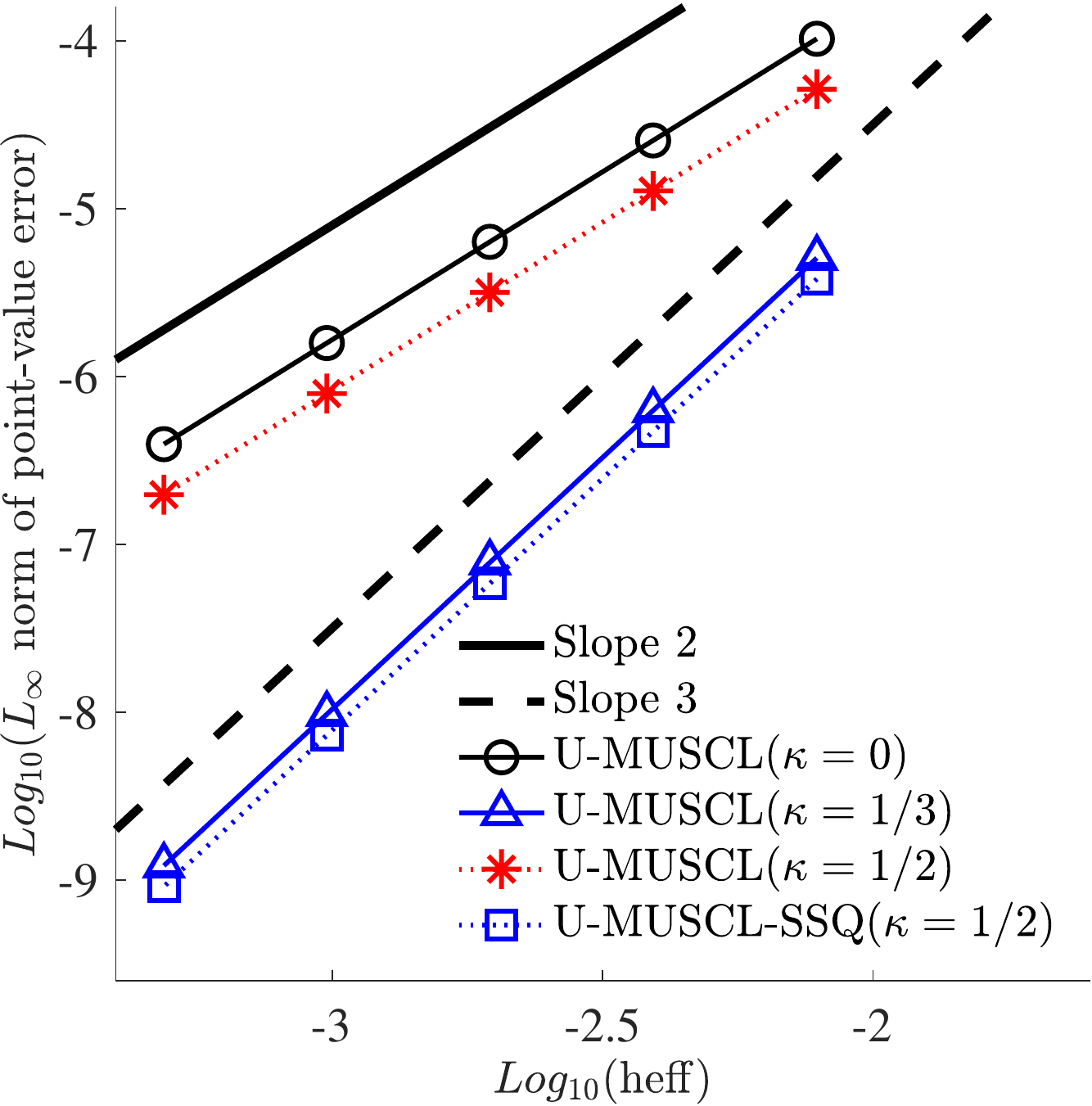}
  \caption[]{Linear advection}
\label{fig:oned_unsteady_p_linear}%
      \end{subfigure}
      \caption{
\label{fig:oned_unsteady_UMSUCL}%
Error convergence results of the U-MSUCL scheme (\ref{fv_form_1D_FD}) with various values of $\kappa$ for unsteady Burger's and linear advection equations.} 
\end{figure}

We tested the same three schemes but simply by integrating the form (\ref{fv_form_1D_ca}) with $s_i=0$, as typically done in the U-MUSCL scheme. However, the initial values are given by the cell-averages of $\sin(2 \pi x)$ for FVC, otherwise evaluated pointwise, $\sin(2 \pi x_i)$. 
Results are shown in Figure \ref{fig:oned_unsteady}. As expected, FVC(1/3) achieves third-order accuracy in the cell-averaged solution and FSR(1/3) achieves third-order accuracy in the point-valued solution. FVP(1/2) is  second-order accurate because the coupling terms in the time derivative as in Equation  (\ref{fv_form_1D_point_QUICK}) are ignored. This means that the U-MUSCL scheme in the typical form (i.e., without a consistent spatial discretization of the time derivative term or without the flux reconstruction) cannot be third-order accurate.

Finally, we performed the same unsteady calculations with the U-MUSCL scheme (\ref{fv_form_1D_FD}) with $\kappa=0$, $\kappa=1/3$, and $\kappa=1/2$, and U-MUSCL-SSQ for the nonlinear Burgers equation and a linear advection equation: $\partial_t u +  \partial_x f = 0$, where $f=u^2/2$ and $f=u$, respectively. The same exact solution as before was used in both cases. In the linear case, the initial sine wave is simply convected to the right. 

As expected and similarly to the steady case, the U-MUSCL scheme is second-order accurate when applied to the nonlinear Burgers equation as shown in Figure \ref{fig:oned_unsteady_c_nonlinear}, but achieves third-order accuracy with $\kappa=1/3$ for the linear advection equation as shown in Figure \ref{fig:oned_unsteady_p_linear}. On the other hand, U-MUSCL-SSQ is third-order accurate for both linear and nonlinear equations.

Not shown, but the FVP scheme gives third-order accuracy if the time derivatives are discretized in a compatible manner as demonstrated in Ref[\citen{Nishikawa_3rdQUICK:2020}]. To the best of the authors' knowledge, such compatible time-derivative discretizations, including that of U-MUSCL-SSQ, are never used with the U-MUSCL scheme in practical codes.

  \begin{figure}[htbp!]
    \centering
      \begin{subfigure}[t]{0.32\textwidth}
  \includegraphics[width=0.85\textwidth,trim=0 0 0 0 ,clip]{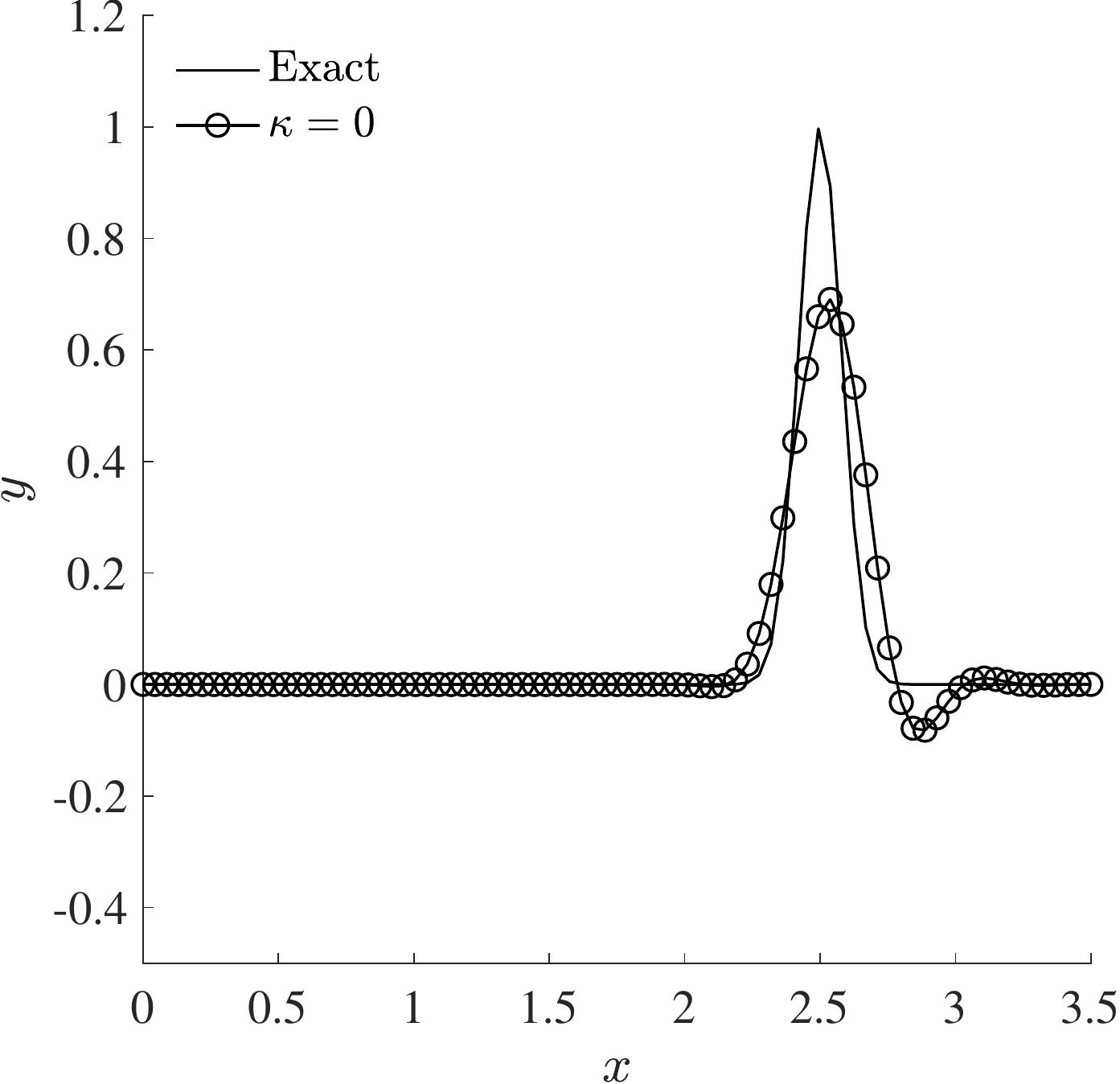}
  \caption[]{$\kappa=0$.}
\label{fig:oned_adv_kappa0p0}%
      \end{subfigure}
      \begin{subfigure}[t]{0.32\textwidth}
  \includegraphics[width=0.85\textwidth,trim=0 0 0 0 ,clip]{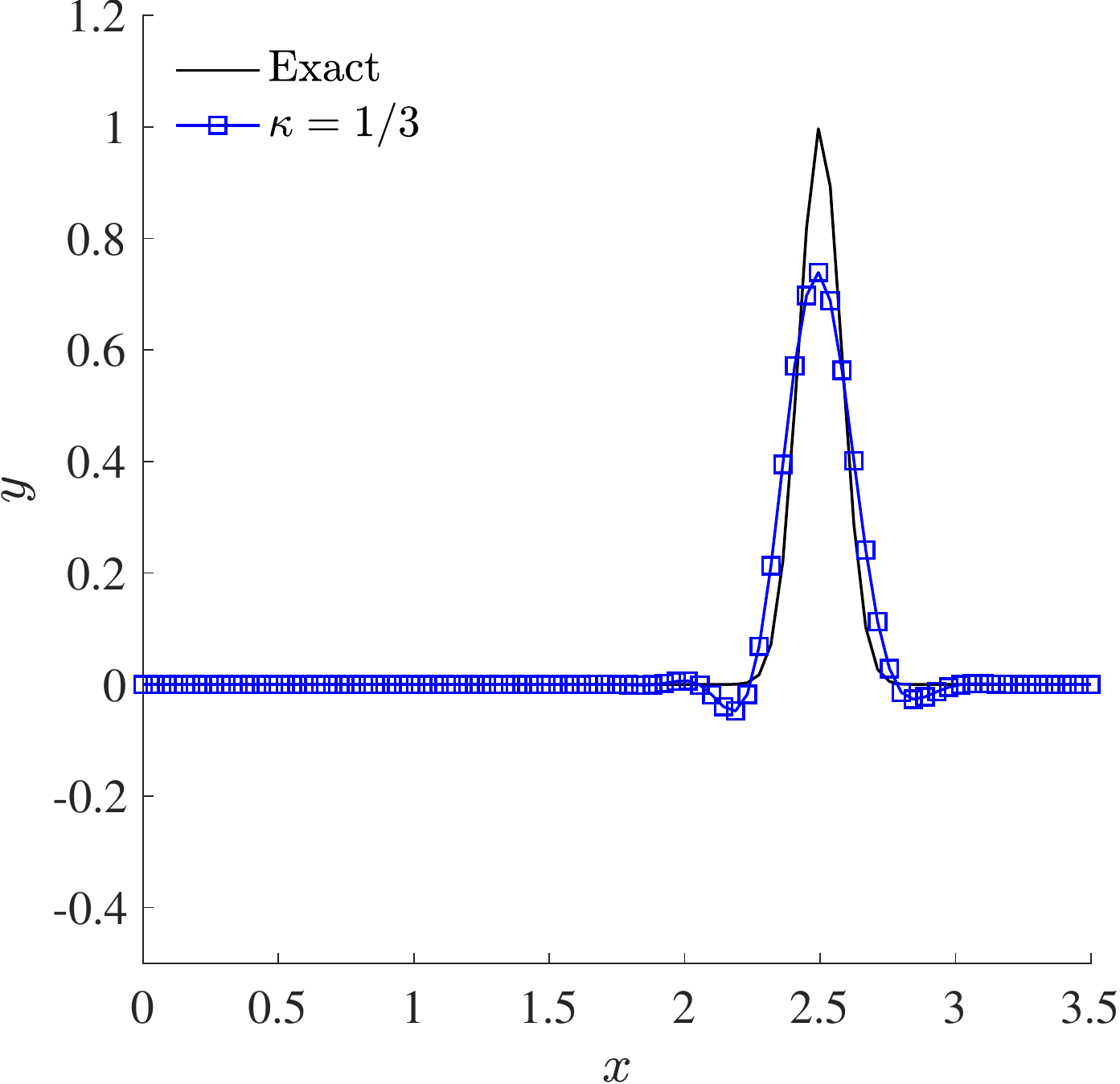}
  \caption[]{$\kappa=1/3$.}
\label{fig:oned_adv_kappa1o3}%
      \end{subfigure}
      \begin{subfigure}[t]{0.32\textwidth}
  \includegraphics[width=0.85\textwidth,trim=0 0 0 0 ,clip]{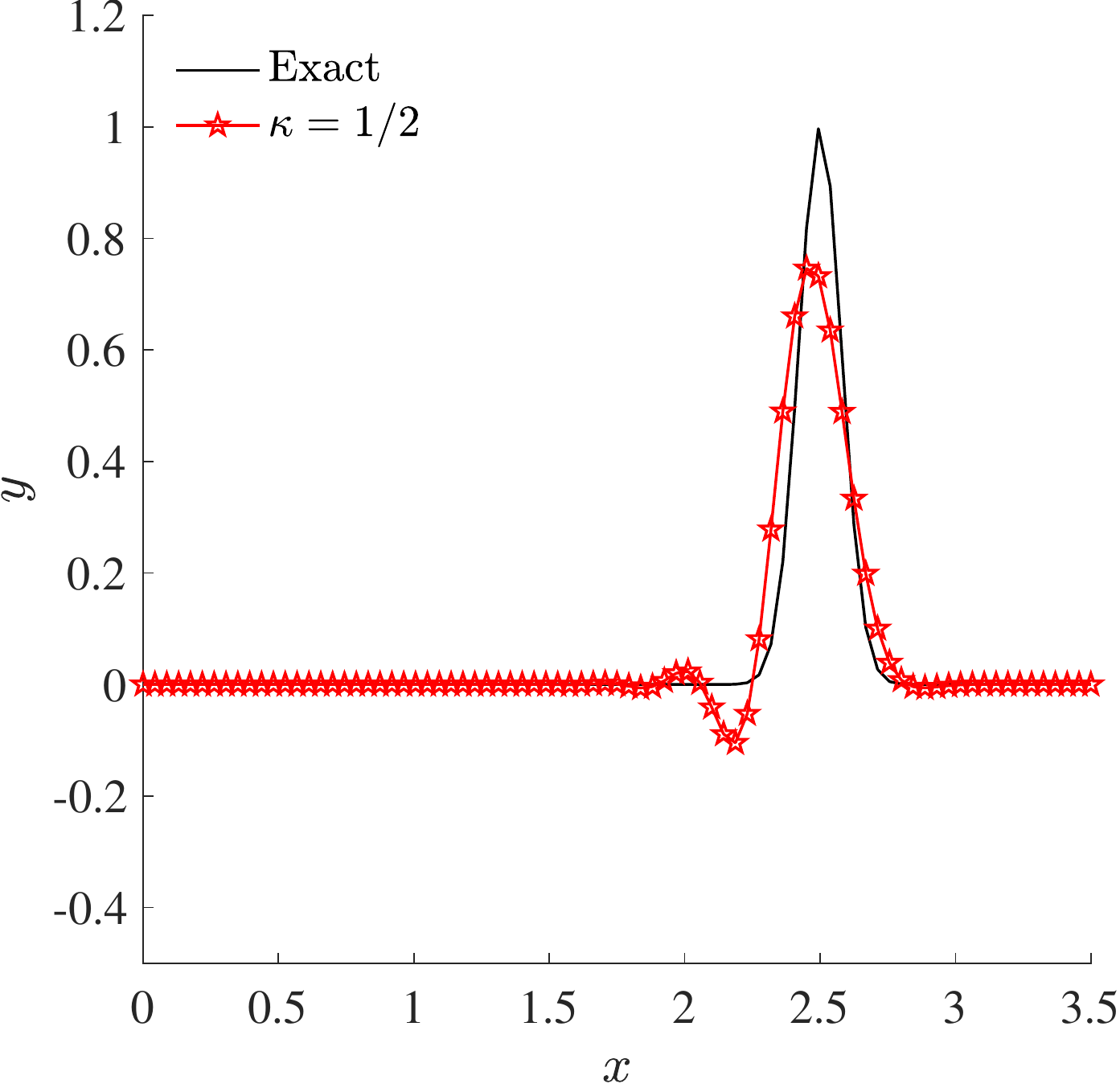}
  \caption[]{$\kappa=1/2$.}
\label{fig:oned_adv_kappa1o2}%
      \end{subfigure}
      \begin{subfigure}[t]{0.32\textwidth}
  \includegraphics[width=0.85\textwidth,trim=0 0 0 0 ,clip]{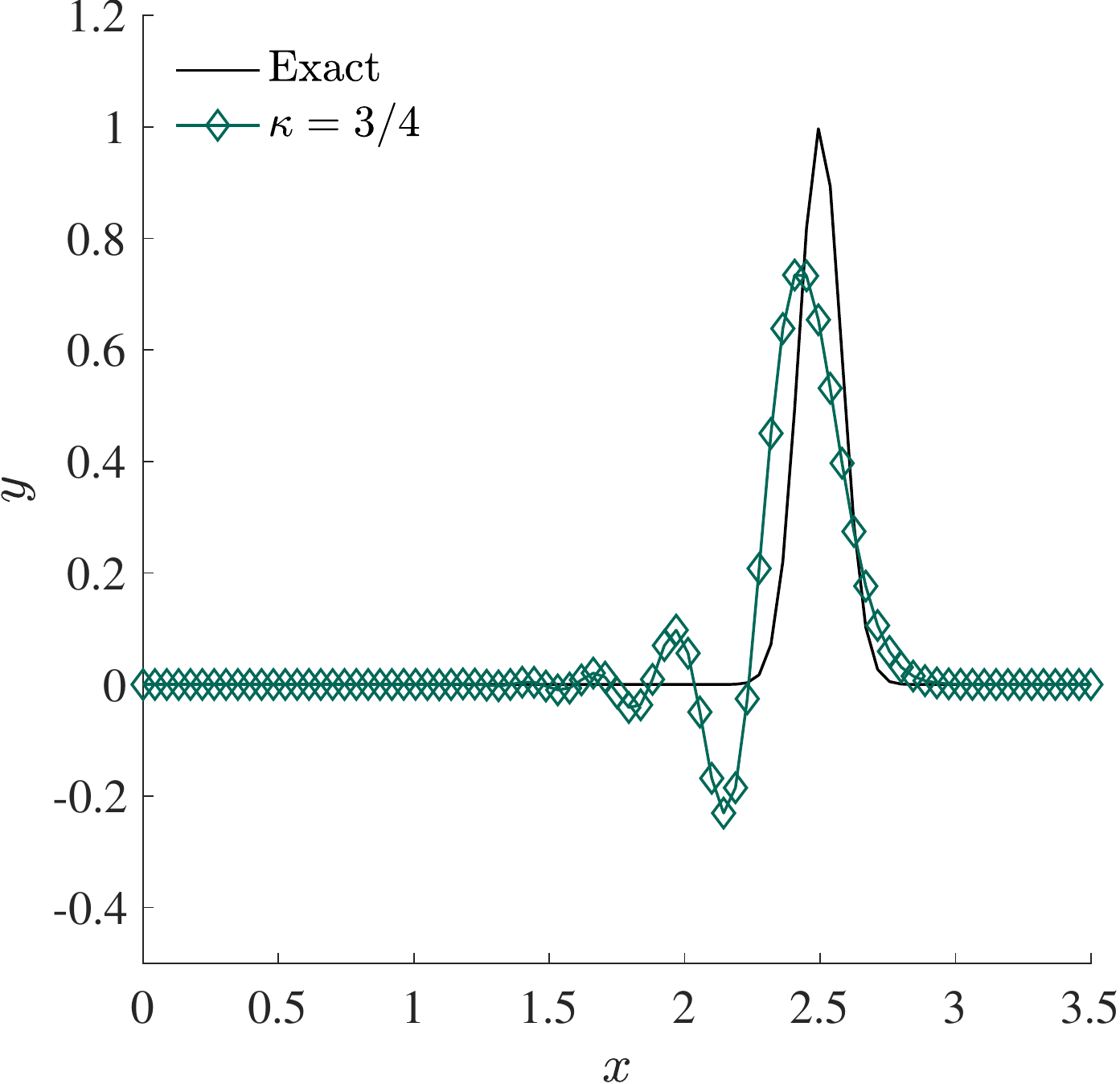}
  \caption[]{$\kappa=3/4$.}
\label{fig:oned_adv_kappa3o4}%
      \end{subfigure}
      \begin{subfigure}[t]{0.32\textwidth}
  \includegraphics[width=0.85\textwidth,trim=0 0 0 0 ,clip]{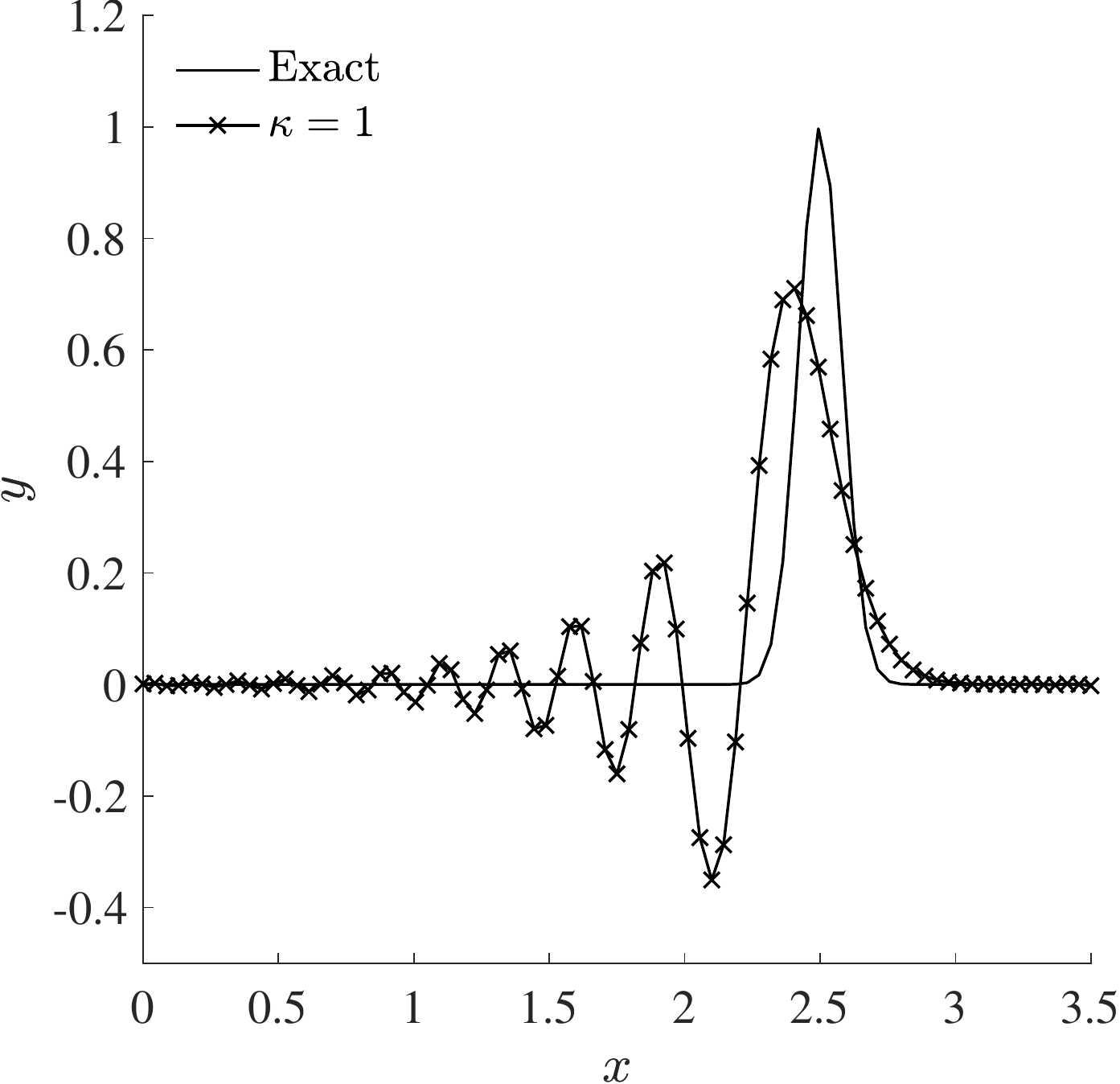}
  \caption[]{$\kappa=1$.}
\label{fig:oned_adv_kappa1p0}%
      \end{subfigure}
      \begin{subfigure}[t]{0.32\textwidth}
  \includegraphics[width=0.85\textwidth,trim=0 0 0 0 ,clip]{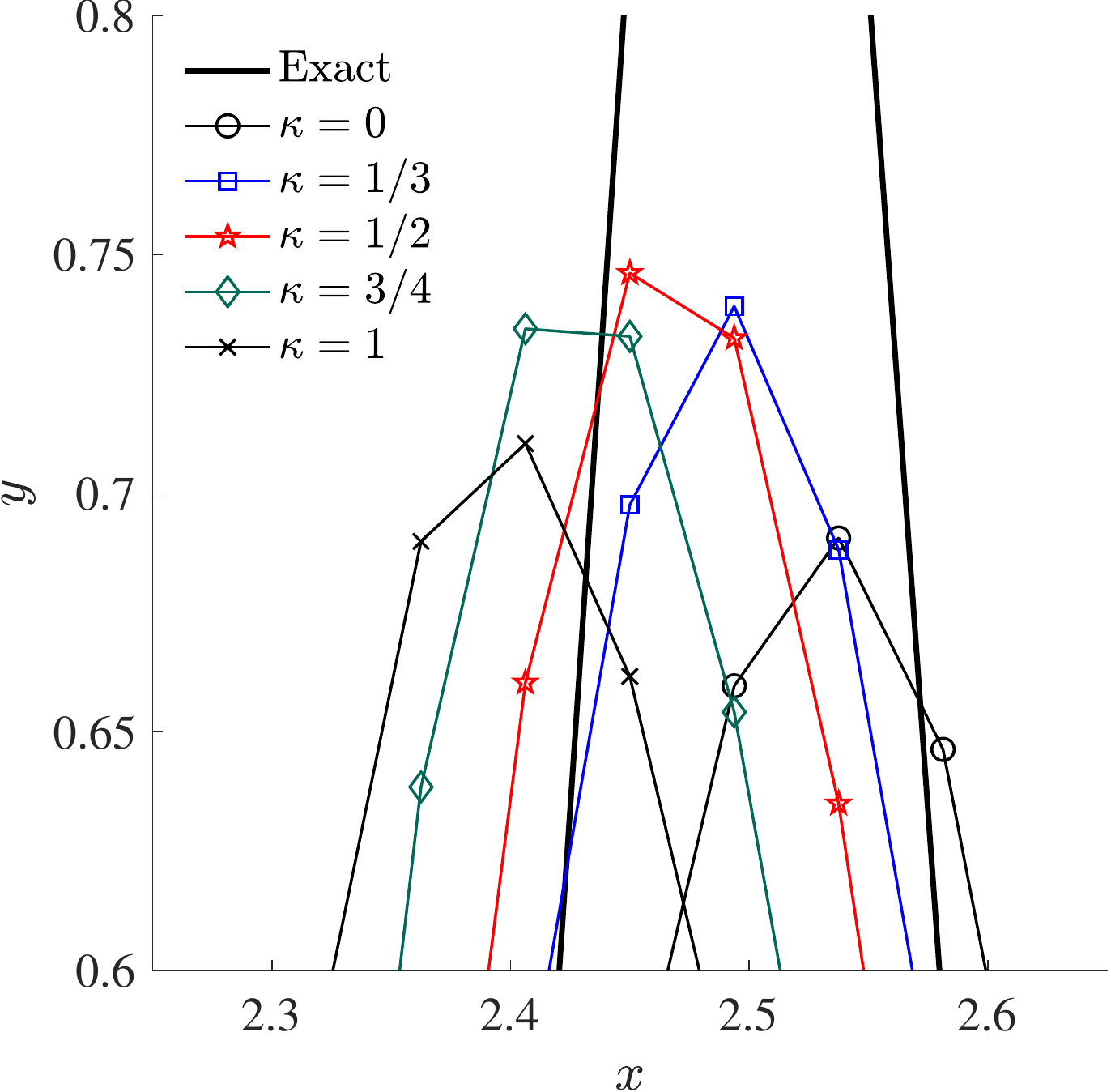}
  \caption[]{Comparison.}
\label{fig:oned_adv_kappa_comparison}%
      \end{subfigure}
      \caption{
\label{fig:oned_adv_kappa}%
Numerical solutions obtained with various values of $\kappa$ for the linear advection: $\partial_t u + \partial_x u = 0$ at the final time $t=2.0$ on a uniform grid with 81 nodes.} 
\end{figure}
  
\subsubsection{Dissipative and dispersive errors}
\label{oned_diss_disp}

We provide a qualitative comparison among various values of $\kappa$ for a linear advection equation $\partial_t u + a \partial_x u = 0$ with a Gaussian pulse: $u(x) = \exp \left( - 80 (x-1/2)^2 \right)$ traveling to the right at a speed $a=1$ to the final time $t=2.0$. The grid is a uniform grid with 81 nodes. The time integration is performed with the third-order SSP Runge-Kutta scheme as before with CFL$=0.1$. 

Results are shown in Figure \ref{fig:oned_adv_kappa}. We observe severe dispersive errors towards $\kappa=1$ (the central scheme); high-frequency modes travel at wrong speeds and lag behind. This is a typical behavior of second-order schemes whose leading error is dispersive. The most accurate solution is obtained with $\kappa=1/3$, which is third-order for linear problems and thus eliminates the leading dispersive error of  second-order schemes. Hence, it is not necessarily less dissipative than second-order schemes. To reduce dissipation, one has to employ a fourth-order scheme or a low-dissipation numerical flux. The latter has been demonstrated as a practical approach in Ref.[\citen{nishikawa_liu_aiaa2018-4166}]. 
It is important to note also that less dissipative schemes with $\kappa$ closer to $1$ does not necessarily better preserve the peak of the Gaussian pulse. See Figure \ref{fig:oned_adv_kappa_comparison}. The peak is better kept with $\kappa=1/2$ or $\kappa=1/3$, not with the zero dissipation scheme of $\kappa=1$. We will have a similar observation for a nonlinear problems later.

\subsection{Zero jump for any $\kappa$ and exact quadratic extrapolation with $\kappa=1/2$}
\label{quadratic_extrapolation}

To verify the zero solution jump for any $\kappa$ and the exact quadratic extrapolation with $\kappa=1/2$ when gradients are computed by a quadratic LSQ method as discussed in Section \ref{truth_umuscl_2_exact_quadratic_kappa1o2}, we consider a series of irregular triangular grids with 2304, 4096, 6400, 9216, 12544 nodes in a unit square domain (similar to the one shown in Figure \ref{fig:twod_verification_tria_irrg_grid}). We compute the maximum norms of the jump $u_R-u_L$ given by the U-MUSCL reconstruction scheme at each edge-midpoint and the error $|u_L - u_{exact} |+| u_R - u_{exact}|$, where $u_{exact}$ is a given function evaluated at the edge-midpoint. We consider two types of functions: a sine function, 
\begin{eqnarray}
u(x,y) = 1 + 0.2 \sin( 2.3 \pi x + 2.5 \pi y ), 
\end{eqnarray}
and a quadratic function,
\begin{eqnarray}
u(x,y) = 8.75 - 1.3 x + 3.7 y + 2.1 x^2 +0.3 x y - 7.5 y^2.
\end{eqnarray}
$L_\infty$ norms of the solution jump and the error are computed over all edges in a given grid. 

\begin{figure}[htbp!]
    \centering
      \begin{subfigure}[t]{0.235\textwidth}
  \includegraphics[width=0.99\textwidth,trim=0 0 0 0 ,clip]{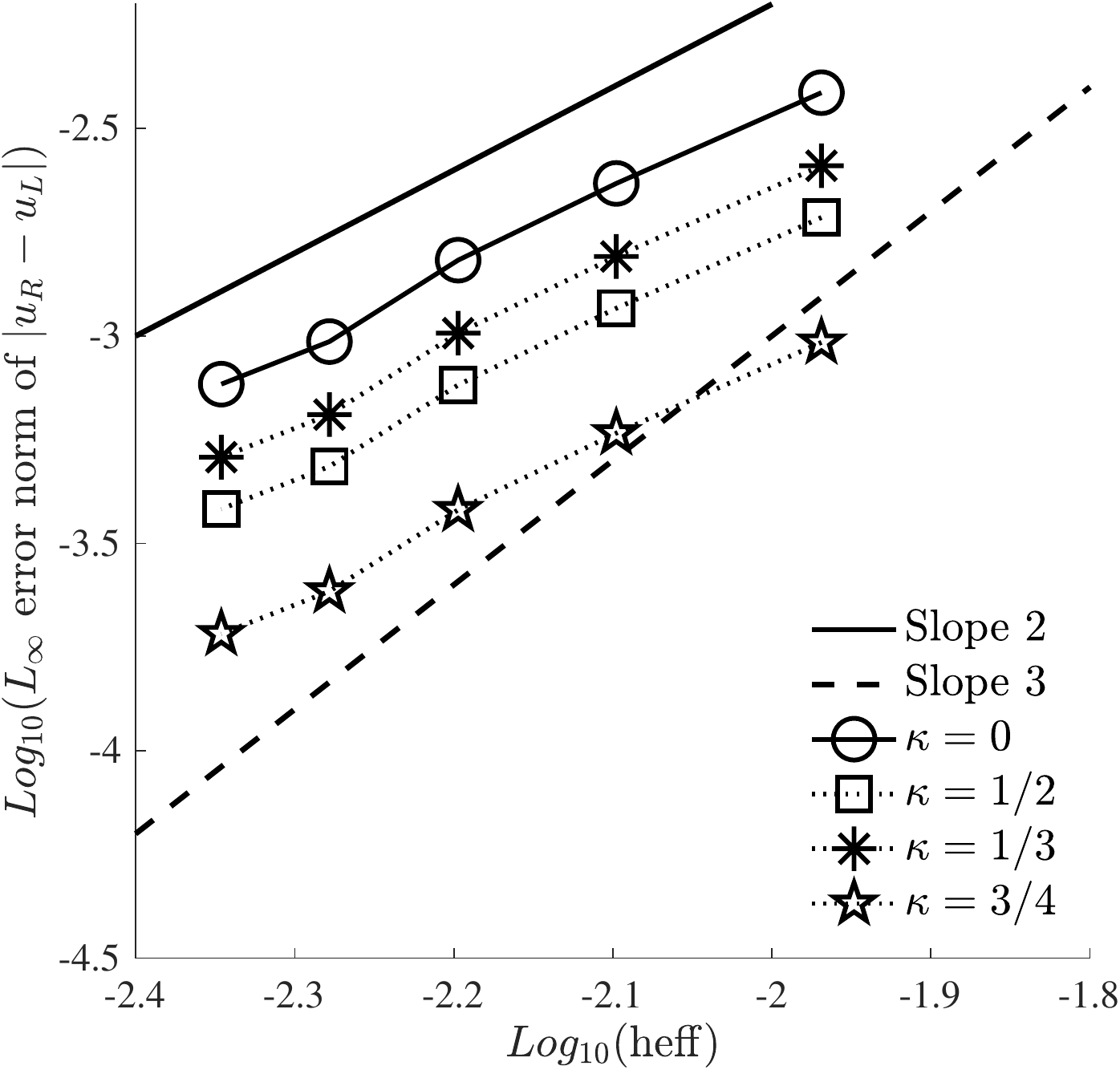}
  \caption[]{Jumps: a sine function.}
\label{fig:twod_test_jump_grad2x2}%
      \end{subfigure}
      \begin{subfigure}[t]{0.235\textwidth}
  \includegraphics[width=0.99\textwidth,trim=0 0 0 0 ,clip]{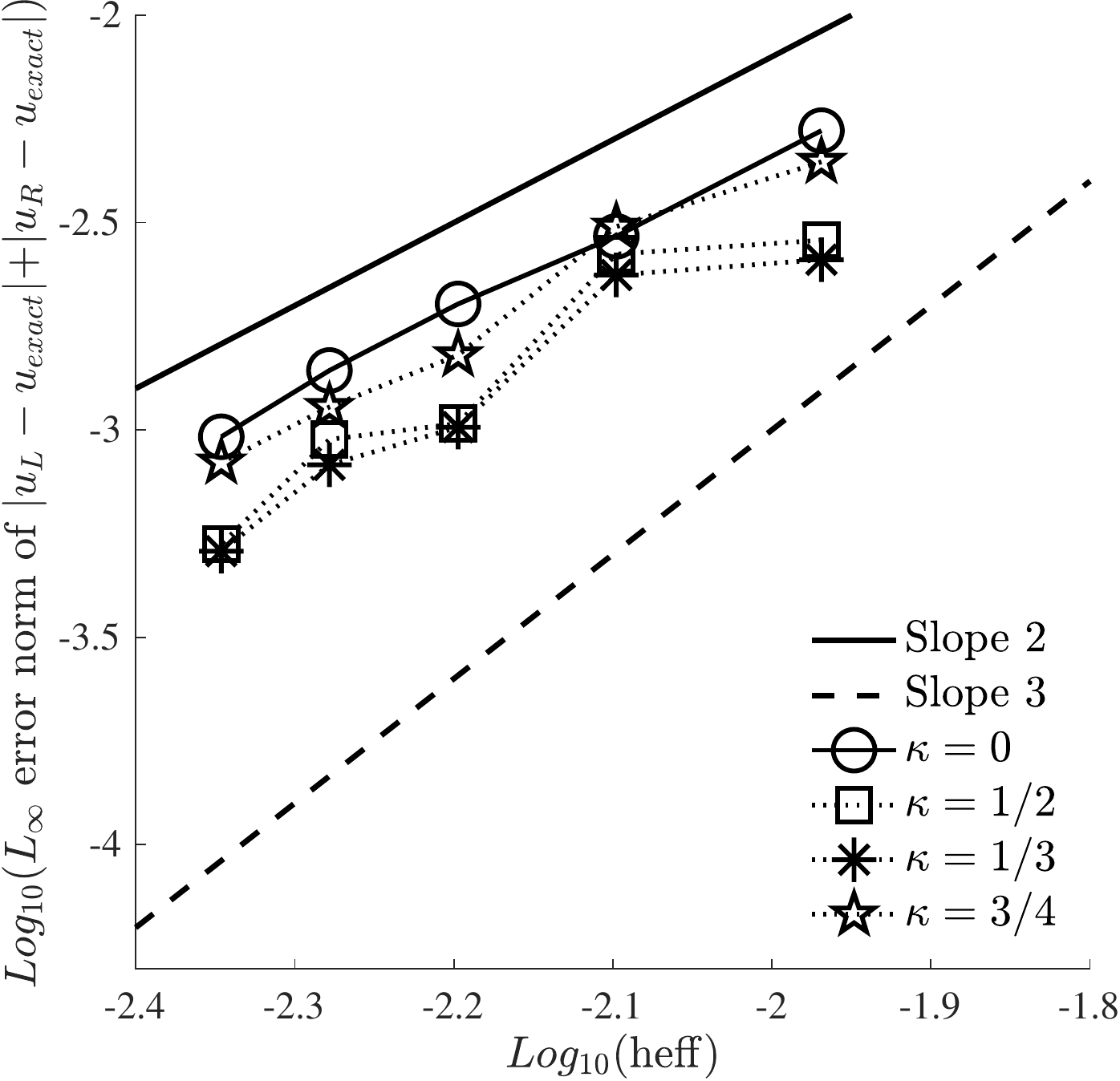}
  \caption[]{Errors: a sine function.}
\label{fig:twod_test_error_grad2x2}%
      \end{subfigure}
      \begin{subfigure}[t]{0.235\textwidth}
  \includegraphics[width=0.99\textwidth,trim=0 0 0 0 ,clip]{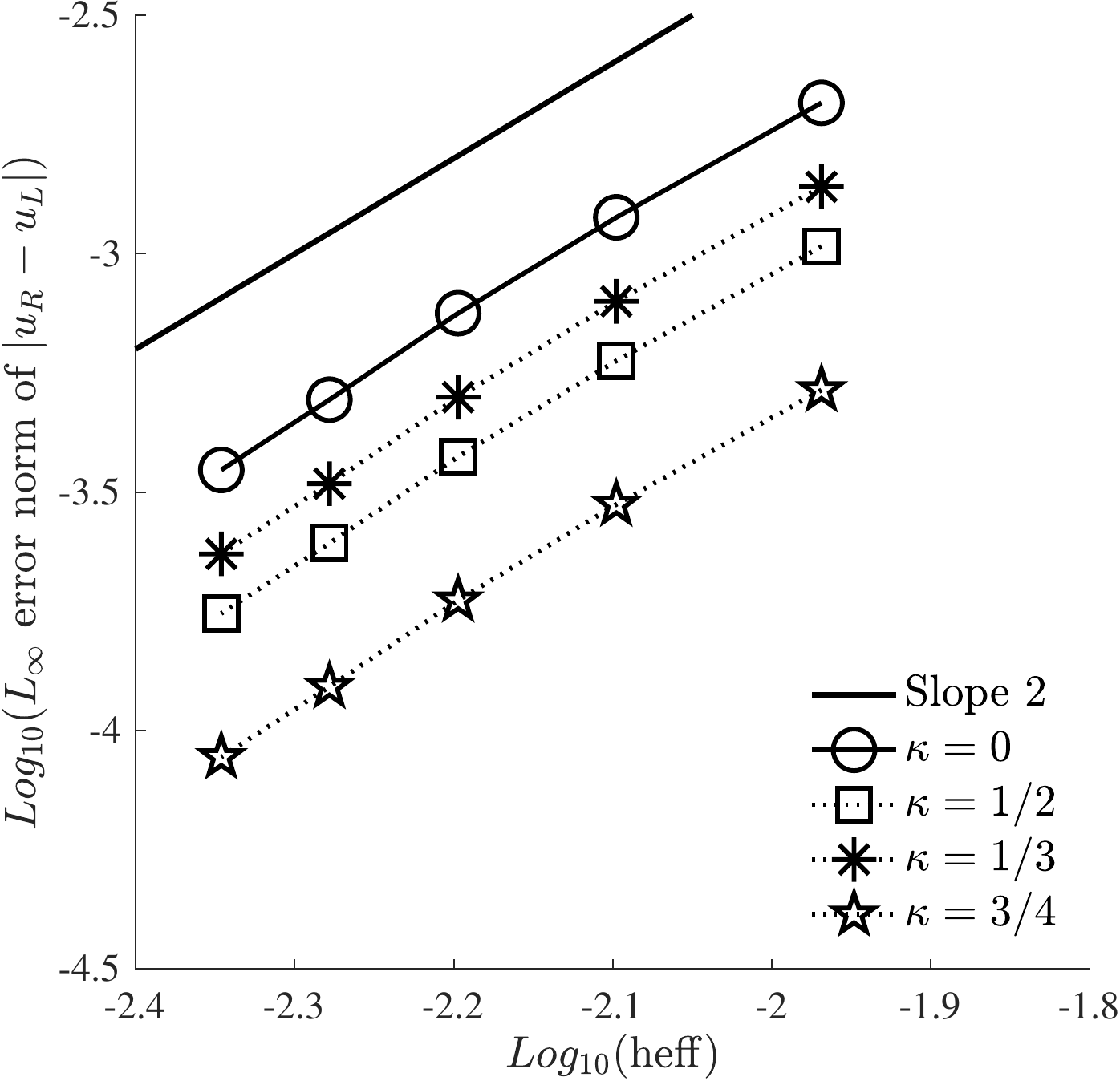}
  \caption[]{Jumps: quadratic.}
\label{fig:twod_test_jump_grad2x2_2}%
      \end{subfigure}
      \begin{subfigure}[t]{0.235\textwidth}
  \includegraphics[width=0.99\textwidth,trim=0 0 0 0 ,clip]{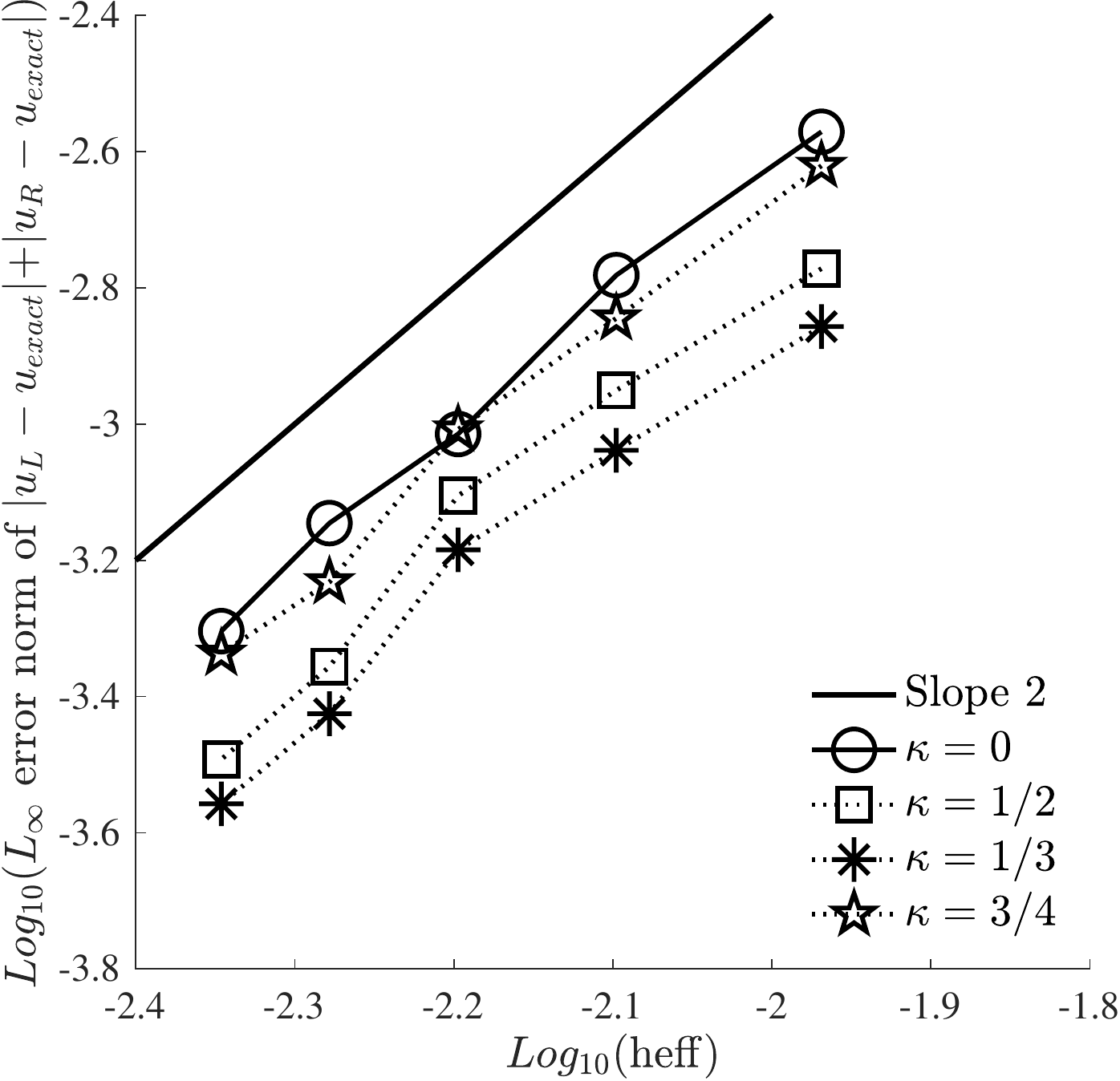}
  \caption[]{Errors: quadratic.}
\label{fig:twod_test_error_grad2x2_2}%
      \end{subfigure}
      \caption{
\label{fig:twod_test_jump_error_grad2x2}%
Solution jump and error at the edge midpoint with linear LSQ gradients.} 
\end{figure}

First, we employed a linear LSQ method and obtained the results shown in Figure \ref{fig:twod_test_jump_error_grad2x2}. For the sine function, both the jump and the error are of $O(h^2)$  as expected. The same results are obtained for the quadratic function. Note that the jump decreases as $\kappa$ increases. This is expected because the jump will vanish identically with $\kappa=1$. Then, if we compare $\kappa=1/3$ and $\kappa=1/2$, we see that the latter gives slightly smaller jumps (less dissipative). 

\begin{figure}[htbp!]
    \centering
      \begin{subfigure}[t]{0.235\textwidth}
  \includegraphics[width=0.99\textwidth,trim=0 0 0 0 ,clip]{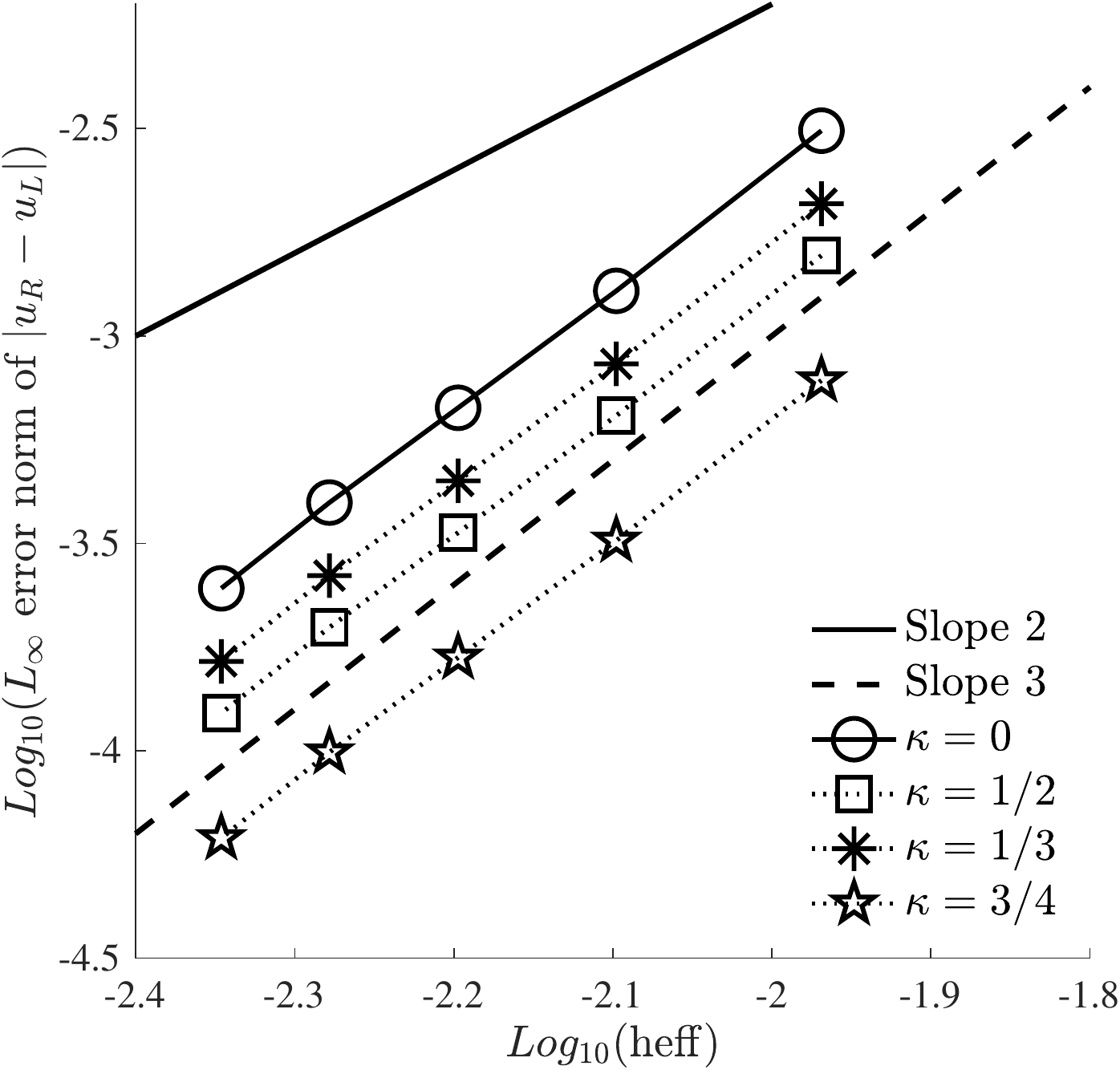}
  \caption[]{Jumps: a sine function.}
\label{fig:twod_test_jump_grad5x5}%
      \end{subfigure}
      \begin{subfigure}[t]{0.235\textwidth}
  \includegraphics[width=0.99\textwidth,trim=0 0 0 0 ,clip]{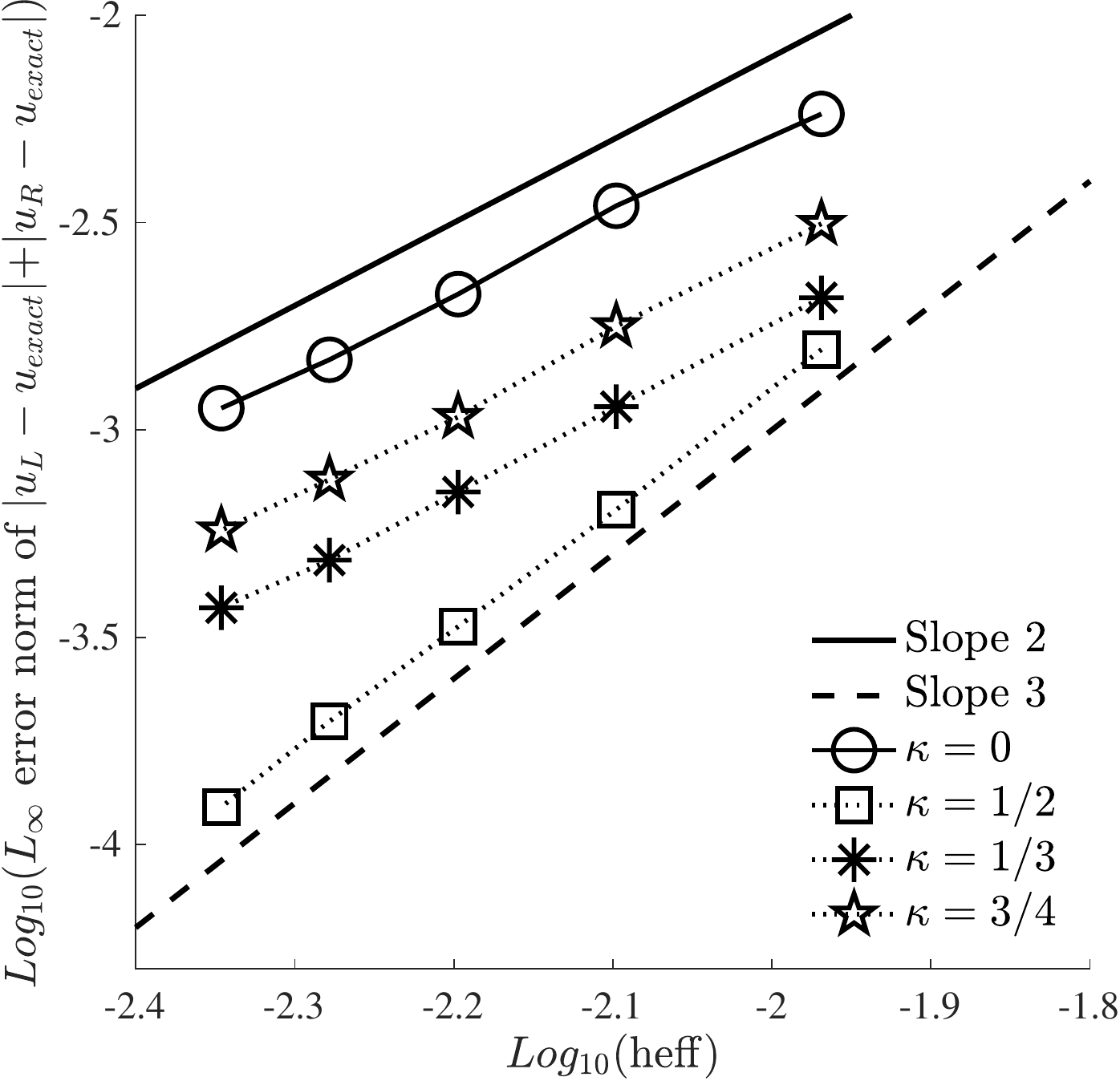}
  \caption[]{Errors: a sine function.}
\label{fig:twod_test_error_grad5x5}%
      \end{subfigure}
      \begin{subfigure}[t]{0.235\textwidth}
  \includegraphics[width=0.99\textwidth,trim=0 0 0 0 ,clip]{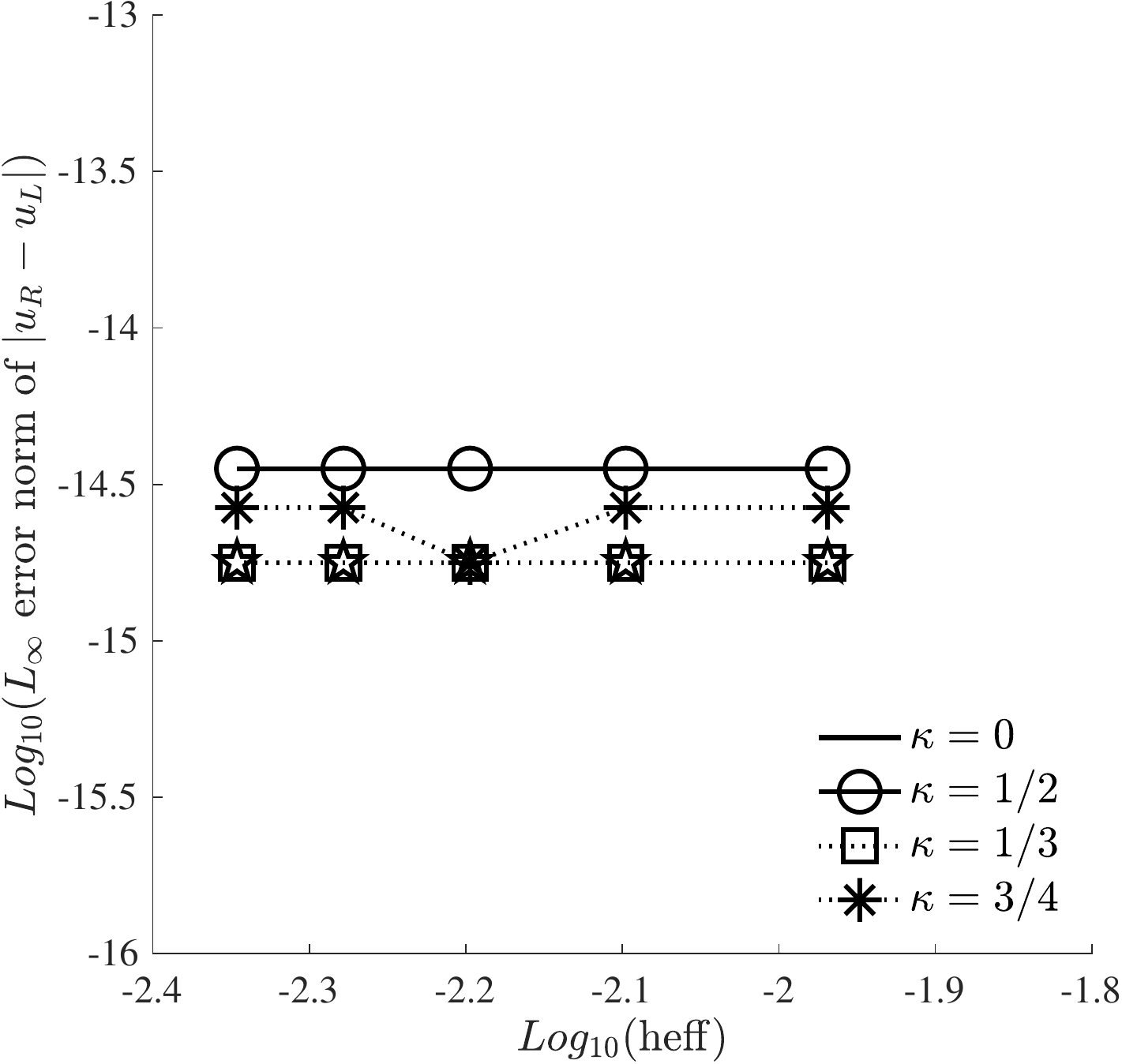}
  \caption[]{Jumps: quadratic.}
\label{fig:twod_test_jump_grad5x5_2}%
      \end{subfigure}
      \begin{subfigure}[t]{0.235\textwidth}
  \includegraphics[width=0.99\textwidth,trim=0 0 0 0 ,clip]{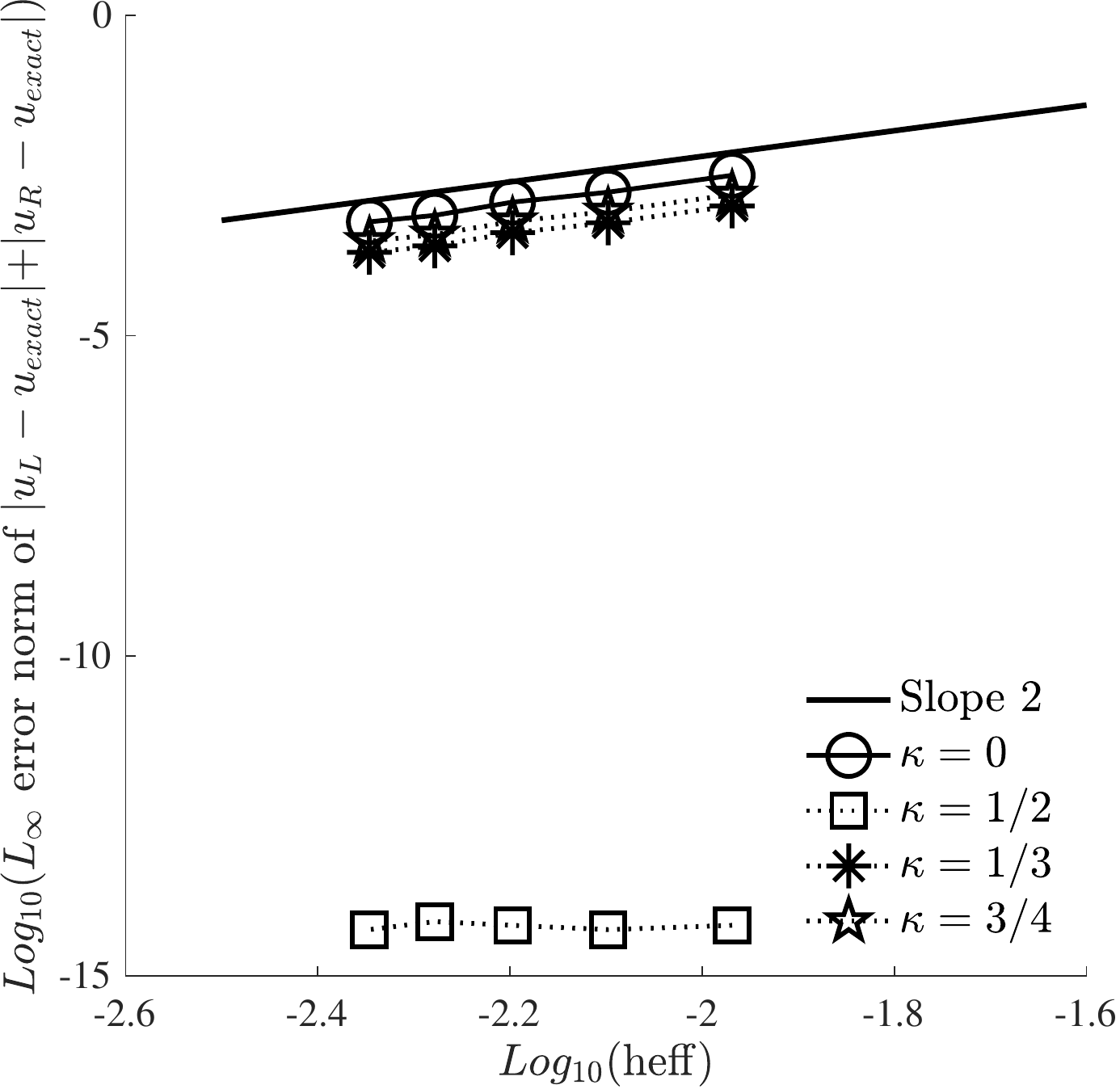}
  \caption[]{Errors: quadratic.}
\label{fig:twod_test_error_grad5x5_2}%
      \end{subfigure}
      \caption{
\label{fig:twod_test_jump_error_grad5x5}%
Solution jump and error at the edge midpoint with quadratic LSQ gradients.} 
\end{figure}

Next, we employed a quadratic LSQ method and obtained the results shown in Figure \ref{fig:twod_test_jump_error_grad5x5}. For the sine function, the jump is of $O(h^3)$ as shown in Figure \ref{fig:twod_test_jump_grad5x5}; this is expected because the jump vanishes for any $\kappa$ for a quadratic function, which is then numerically confirmed as shown in Figure \ref{fig:twod_test_jump_grad5x5_2}. On the other hand, the error is $O(h^2)$ for $\kappa=0$, $1/3$, and $3/4$ but $O(h^3)$ for $\kappa=1/2$ as shown in Figure \ref{fig:twod_test_error_grad5x5}; this is expected because the reconstruction is exact for quadratic functions with $\kappa=1/2$ as confirmed in Figure \ref{fig:twod_test_error_grad5x5_2}. Again, the jump is smaller with $\kappa=1/2$ than with $\kappa=1/3$.

\begin{table}[t]
\ra{1.1}
\begin{center}
{\tabulinesep=0.5mm
\begin{tabu}{lcclllll}\hline\hline 
\multicolumn{1}{l}{ \multirow{2}{*}{Type}      }        &
\multicolumn{1}{l}{ \multirow{2}{*}{$\kappa$}      }        &
\multicolumn{1}{l}{ \multirow{2}{*}{$\theta$}      }        &
\multicolumn{1}{l}{ \multirow{2}{*}{$s$, $\frac{\partial u}{\partial t}$} } & \multicolumn{1}{l}{ \multirow{2}{*}{Regular quadrilaterals}       } &
\multicolumn{1}{l}{ \multirow{2}{*}{Regular triangles}  } &
\multicolumn{1}{l}{ \multirow{2}{*}{Irregular grids}   }
\\
& & & & & &
\\ \hline 
  U-MUSCL & $\mathbb{R}$ & -- & point evaluation & $O(h^2)$  & $O(h^2)$  & $O(h^2)$    \\
    CFSR3 & $\mathbb{R}$ & $1/3$ & point evaluation & $\boldsymbol {\color{red} O(h^3)} $  & $\boldsymbol {\color{red} O(h^3)} $ & $O(h^2)$    \\
U-MUSCL-SSQ & $1/2$ & -- & quadrature (\ref{residual_imuscl_s}) &
 ${\color{black} O(h^2)} $, $\kappa_s = 1/6$  & $\boldsymbol {\color{red} O(h^3)} $, $\kappa_s = 1/4$ & $O(h^2)$ \\
  \hline  \hline
\end{tabu}
}

\caption{U-MUSCL-type schemes in two dimensions. $\mathbb{R}$ indicates any real value. The numerical solution is a point value in all schemes. The order of accuracy is for general nonlinear equations.}

\label{Tab.classificaiton_2D}
\end{center}
\end{table}

\subsection{Third-order accuracy verification in two dimensions}
\label{results_third_order}

In this section, we present accuracy verification results obtained for the Euler equations by the method of manufactured solutions on various types of grids. We will compare the U-MUSCL scheme, CFSR3, and U-MUSCL-SSQ as summarized in Table \ref{Tab.classificaiton_2D}. For the U-MUSCL scheme, we will further compare various values of $\kappa$. In all cases, the residual equations are solved by an implicit defect-correction solver with the exact Jacobian of the first-order residual (i.e., no LSQ gradients) relaxed by the Gauss-Seidel scheme. The Rusanov flux is employed here for robustness, which is effective especially for  $\kappa=3/4$. The iteration is taken to be converged when the residual is reduced by six orders of magnitude in the $L_1$ norm. To focus on accuracy of the interior scheme, we specify the exact solution,
{\color{black}
\begin{eqnarray}
\rho =  1.0  +  C \exp\left( \pi (   0.3 x + 0.3 y    ) \right), \quad
u = 0.15  + C \exp\left( \pi (   0.3 x + 0.3 y    ) \right), \\ 
v =  0.02  + C \exp\left( \pi (   0.3 x + 0.3 y    ) \right), \quad
p =  1.0  + C \exp\left( \pi (   0.3 x + 0.3 y    ) \right), 
\end{eqnarray}
}
where we set $C=0.3$, at boundary nodes, their neighbors, and the neighbors of the neighbors. The source term and its gradient are computed numerically and stored at each node. See Ref.[\citen{nishikawa_centroid:JCP2020}] for details
on the source term computation. The discretization error will be shown for the pressure (results are similar for the other variables). 

\subsubsection{Regular quadrilaterals}

\begin{figure}[htbp!]
    \centering
      \begin{subfigure}[t]{0.48\textwidth}
  \includegraphics[width=0.99\textwidth,trim=0 0 0 0 ,clip]{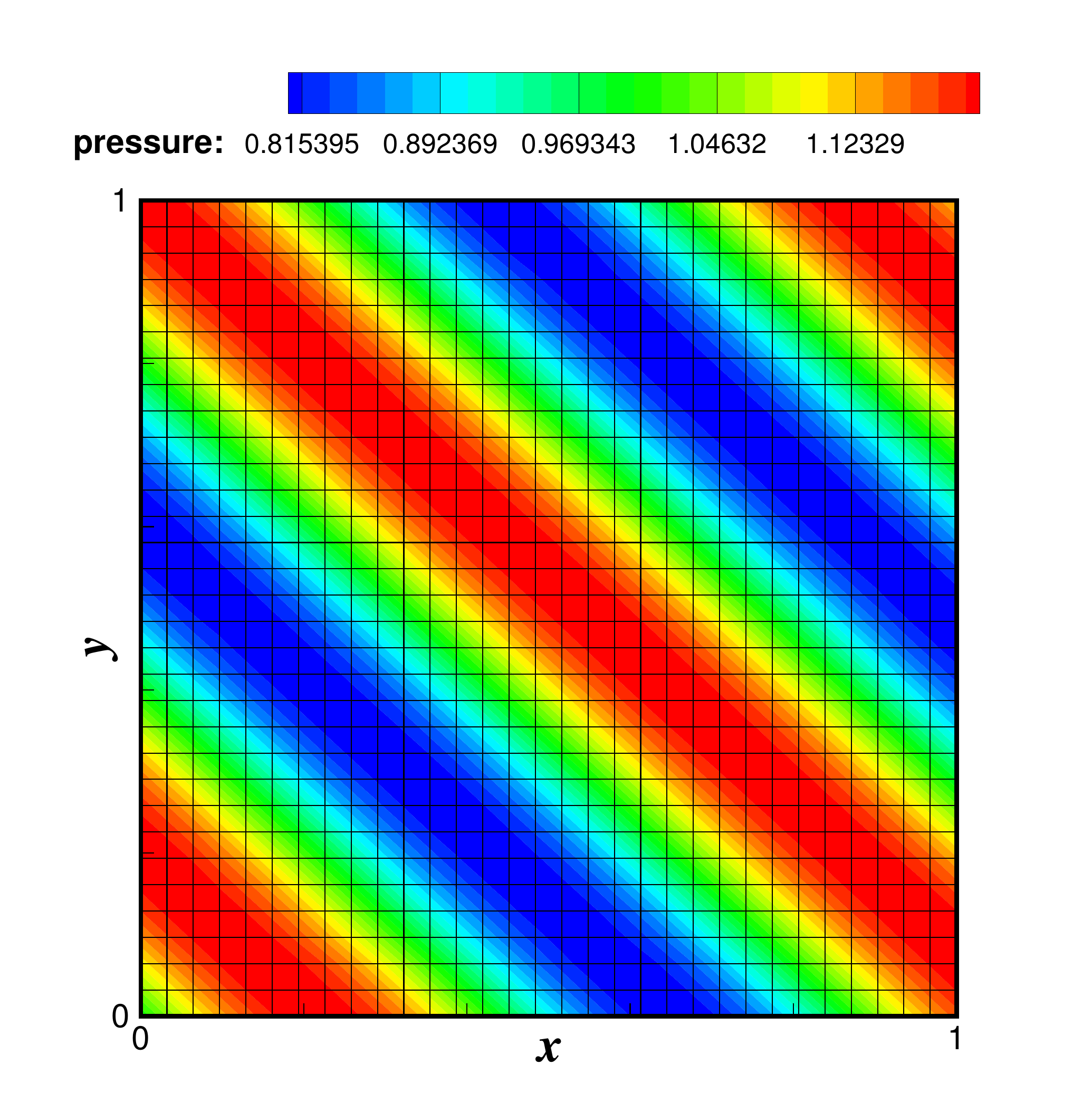}
  \caption[]{The coarsest grid and the exact solution.}
\label{fig:twod_verification_reg_quad_grid}%
      \end{subfigure}
      \begin{subfigure}[t]{0.48\textwidth}
  \includegraphics[width=0.99\textwidth,trim=0 0 0 0, clip]{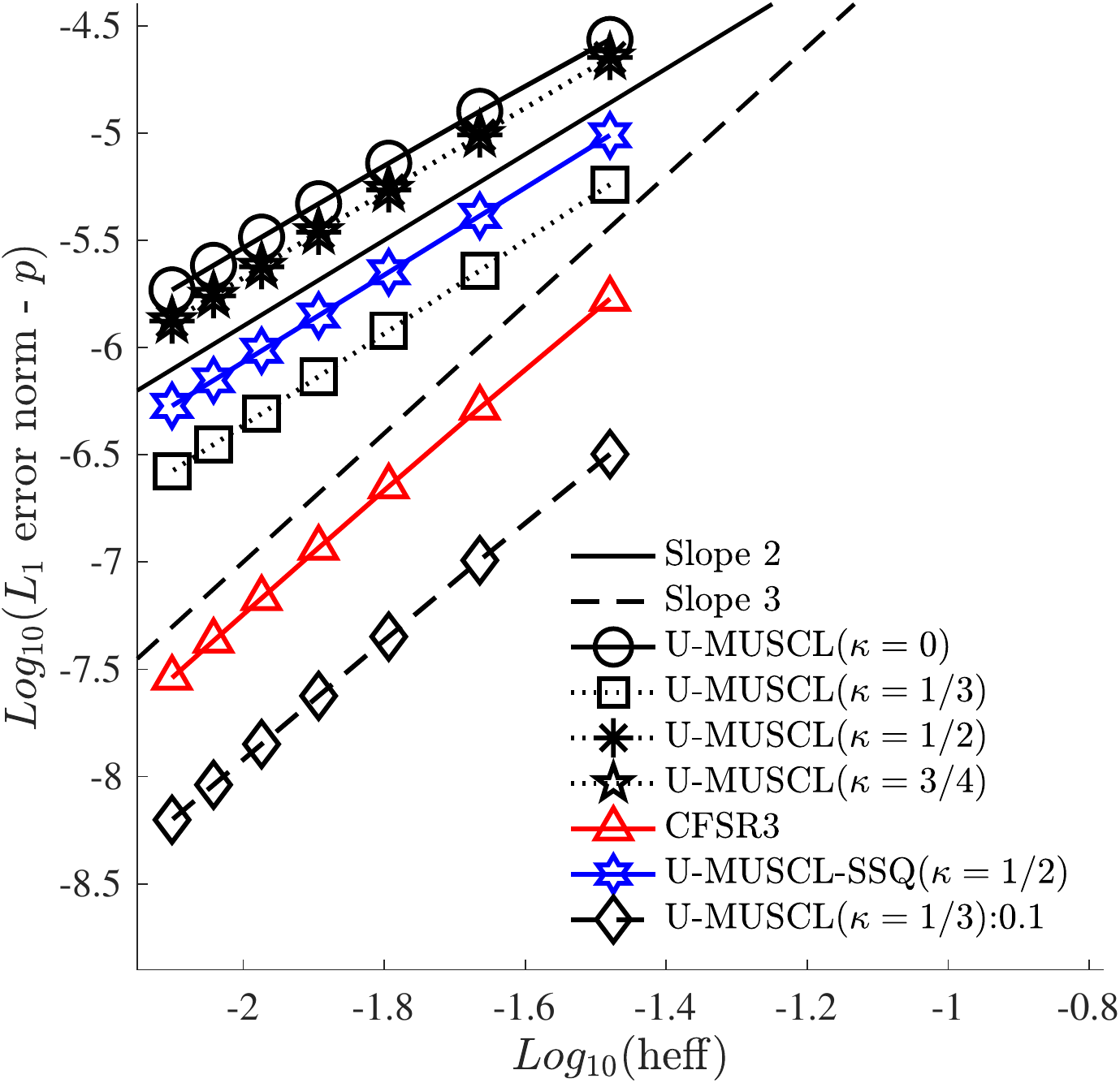}
  \caption[]{Error convergence for the pressure.}
\label{fig:twod_verification_reg_quad_err}%
      \end{subfigure}
      \caption{
\label{fig:twod_verification_reg_quad}%
Error convergence results for the Euler equations on regular quadrilateral grids.} 
\end{figure}

We first consider a series of regular quadrilateral grids with 1024, 2304, 4096, 6400, 9216, 12544, and 16384 nodes. Figure \ref{fig:twod_verification_reg_quad_grid} shows the coarsest grid and the exact solution contours. The error convergence results are shown in Figure \ref{fig:twod_verification_reg_quad_err}. {\color{black} As expected, CFSR3 is third-order accurate. All others are second-order accurate, including U-MUSCL-SSQ with $\kappa_s=1/6$.} Among them, $\kappa=1/3$ gives the lowest errors (black squares). It is close to third-order; the superior accuracy is considered due to the perturbed form of the exact solution used, which effectivey linearizes the Euler equations as shown in  Ref.[\citen{Nishikawa_FakeAccuracy:2020}], and the fact that $\kappa=1/3$ indeed gives third-order accuracy for linear equations. To verify this, we included results for the U-MUSCL($\kappa=1/3$) with a smaller amplitude $C=0.1$, as denoted by U-MUSCL($\kappa=1/3$):0.1 in Figure \ref{fig:twod_verification_reg_quad_err}. Clearly, the scheme achieves third-order accuracy. See Ref.[\citen{Nishikawa_FakeAccuracy:2020}] for further details. 


\subsubsection{Regular equilateral triangles}

\begin{figure}[htbp!]
    \centering
      \begin{subfigure}[t]{0.48\textwidth}
  \includegraphics[width=0.99\textwidth,trim=0 0 0 0 ,clip]{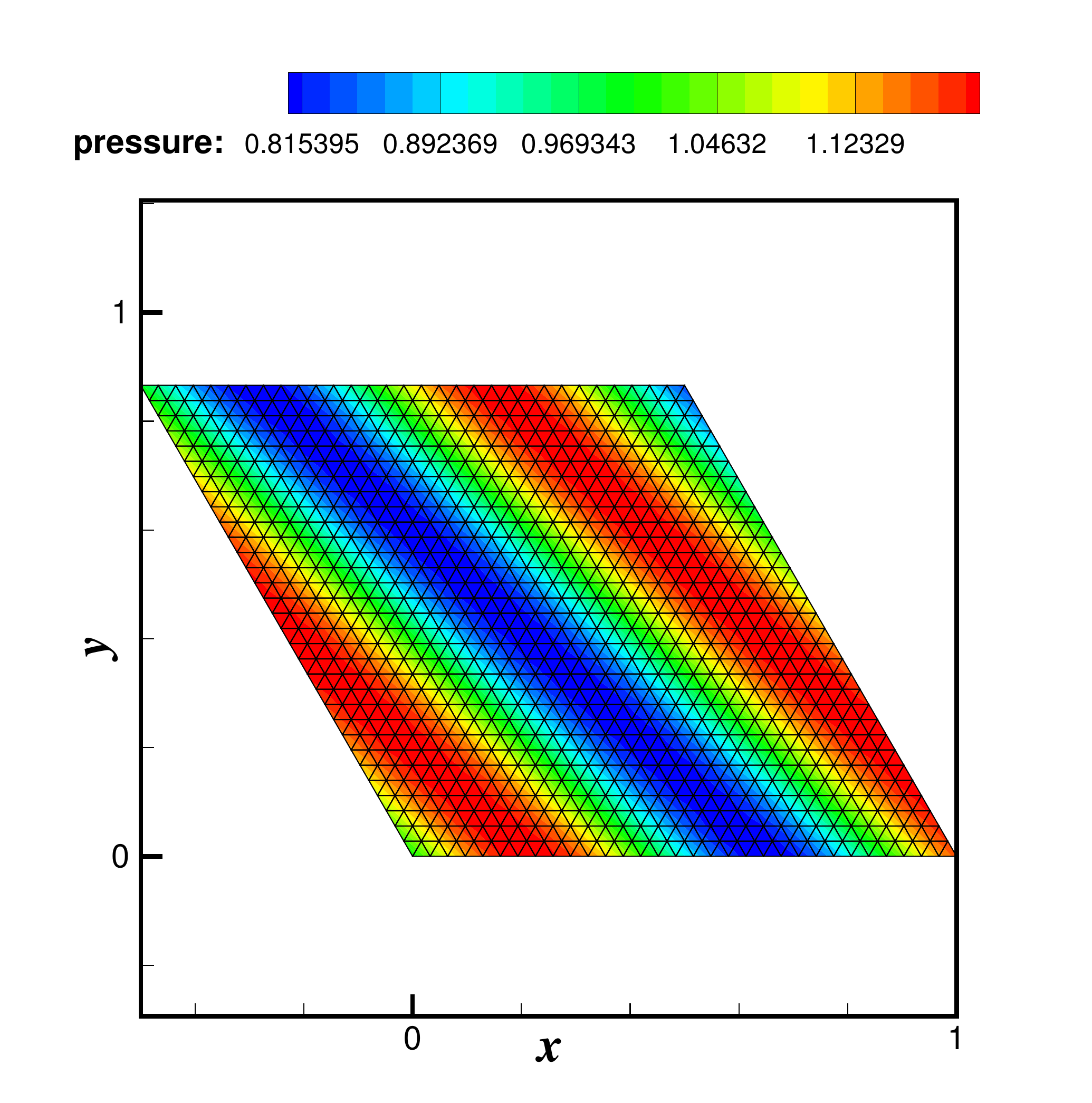}
  \caption[]{The coarsest grid and the exact solution.}
\label{fig:twod_verification_tria_equi_grid}%
      \end{subfigure}
      \begin{subfigure}[t]{0.48\textwidth}
  \includegraphics[width=0.99\textwidth,trim=0 0 0 0, clip]{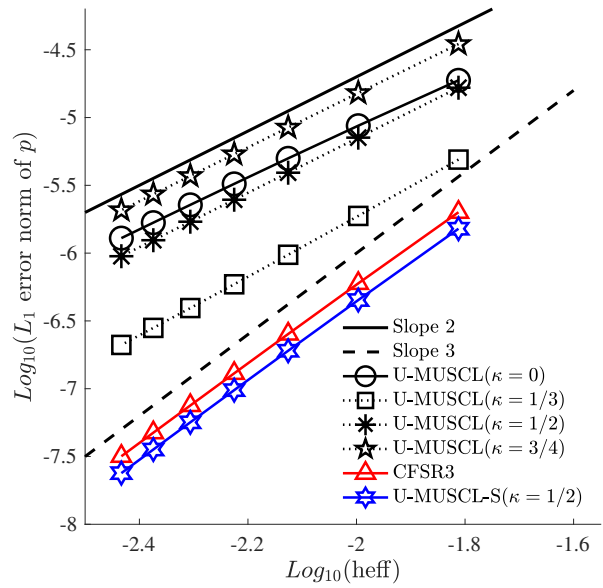}
  \caption[]{Error convergence for the pressure.}
\label{fig:twod_verification_tria_equi_err}%
      \end{subfigure}
      \caption{
\label{fig:twod_verification_tria_equi}%
Error convergence results for the Euler equations on  equilateral triangular grids.} 
\end{figure}

Next, we consider regular equilateral triangular grids with 1024, 2304, 4096, 6400, 12544, and 16384 nodes. The coarsest grid is shown in Figure \ref{fig:twod_verification_tria_equi_grid}. Error convergence results are shown in Figure \ref{fig:twod_verification_tria_equi_err}. First, we see that CFSR3 and U-MUSCL-SSQ with $\kappa_s=1/4$ achieve third-order accuracy as expected. All other schemes are second-order accurate. Again, $\kappa=1/3$ gives the lowest 
second-order error.

\subsubsection{Regular right-isosceles triangles}

\begin{figure}[htbp!]
    \centering
      \begin{subfigure}[t]{0.48\textwidth}
  \includegraphics[width=0.99\textwidth,trim=0 0 0 0 ,clip]{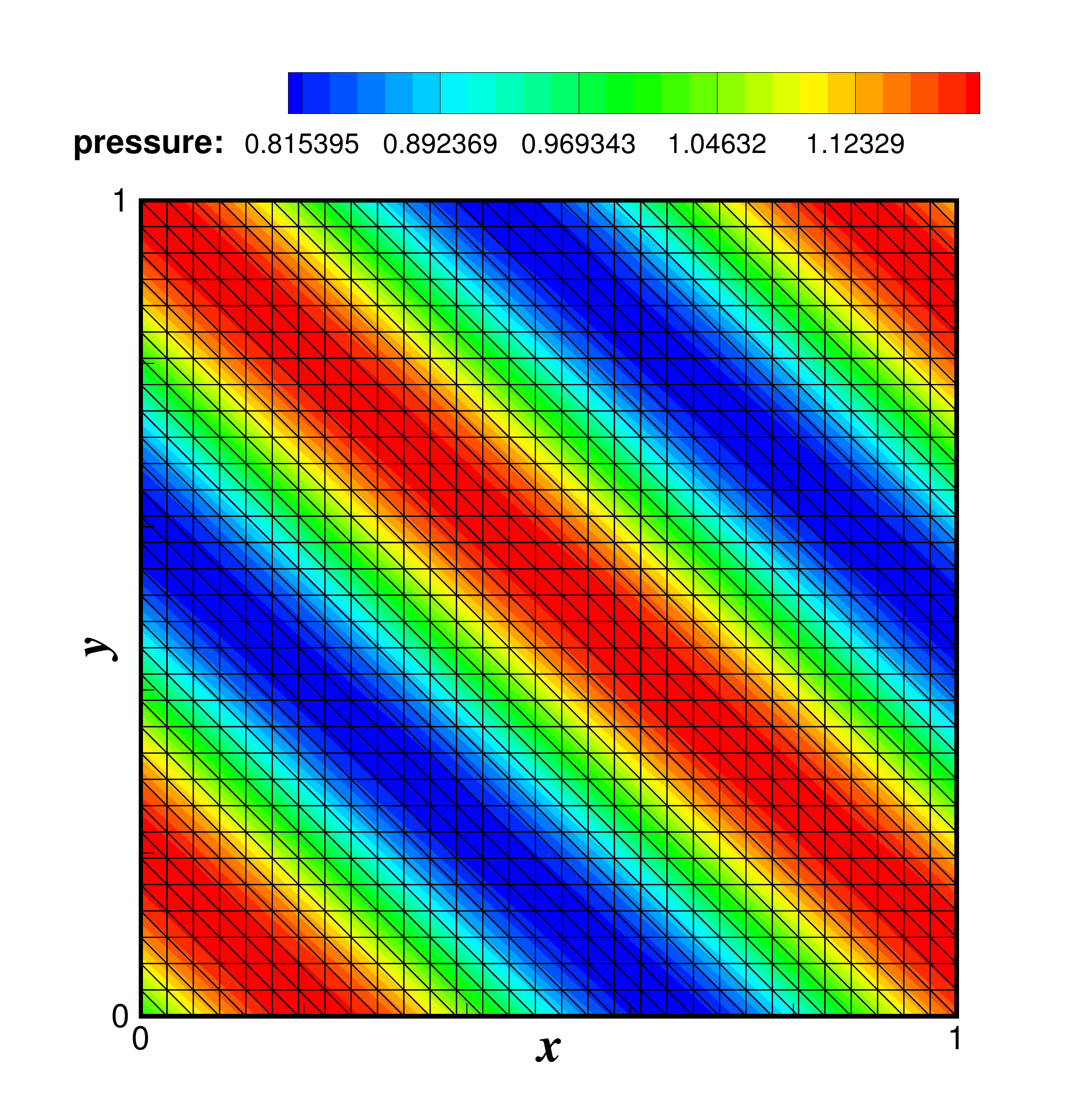}
  \caption[]{The coarsest grid and the exact solution.}
\label{fig:twod_verification_tria_right_grid}%
      \end{subfigure}
      \begin{subfigure}[t]{0.48\textwidth}
  \includegraphics[width=0.99\textwidth,trim=0 0 0 0, clip]{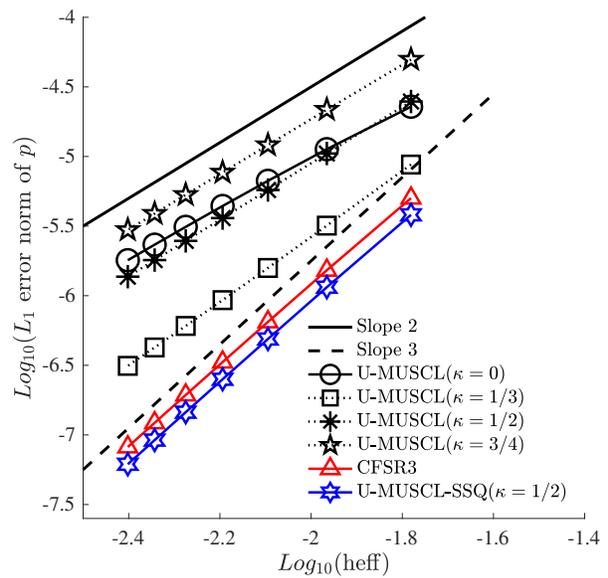}
  \caption[]{Error convergence for the pressure.}
\label{fig:twod_verification_tria_right_err}%
      \end{subfigure}
      \caption{
\label{fig:twod_verification_tria_right}%
Error convergence results for the Euler equations on right  triangular grids.} 
\end{figure}

Next, we consider regular right-isosceles triangular grids with 1024, 2304, 4096, 6400, 12544, and 16384 nodes. The coarsest grid is shown in Figure \ref{fig:twod_verification_tria_right_grid}. Error convergence results are shown in Figure \ref{fig:twod_verification_tria_right_err}. As expected, the results are very similar to those for the equilateral triangular grids in the previous section. CFSR3 and U-MUSCL-SSQ with $\kappa_s=1/4$ are both third-order accurate as long as the grid is regular and triangular.  

\subsubsection{Irregular triangles}

\begin{figure}[htbp!]
    \centering
      \begin{subfigure}[t]{0.48\textwidth}
  \includegraphics[width=0.99\textwidth,trim=0 0 0 0 ,clip]{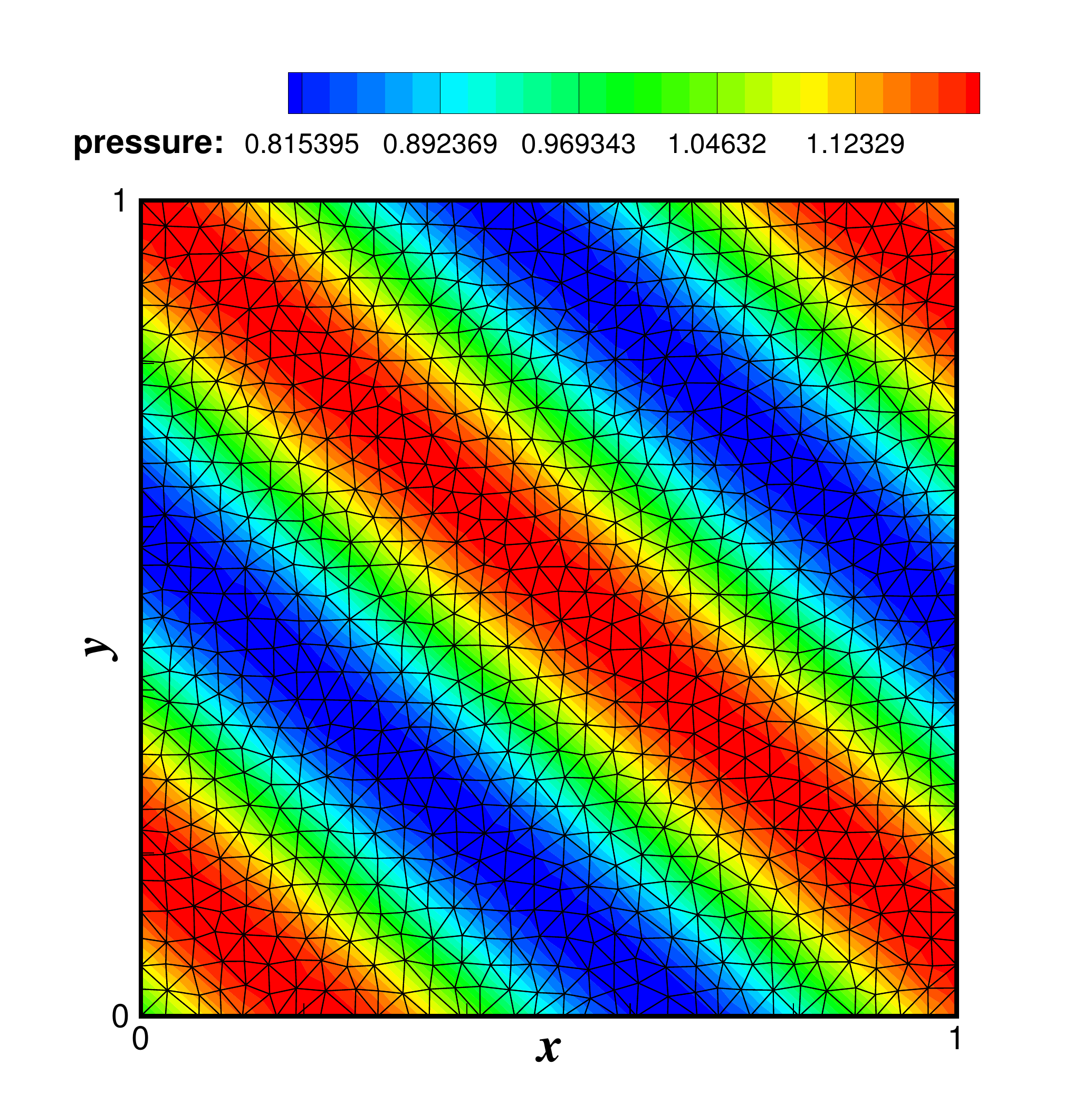}
  \caption[]{The coarsest grid and the exact solution.}
\label{fig:twod_verification_tria_irrg_grid}%
      \end{subfigure}
      \begin{subfigure}[t]{0.48\textwidth}
  \includegraphics[width=0.99\textwidth,trim=0 0 0 0, clip]{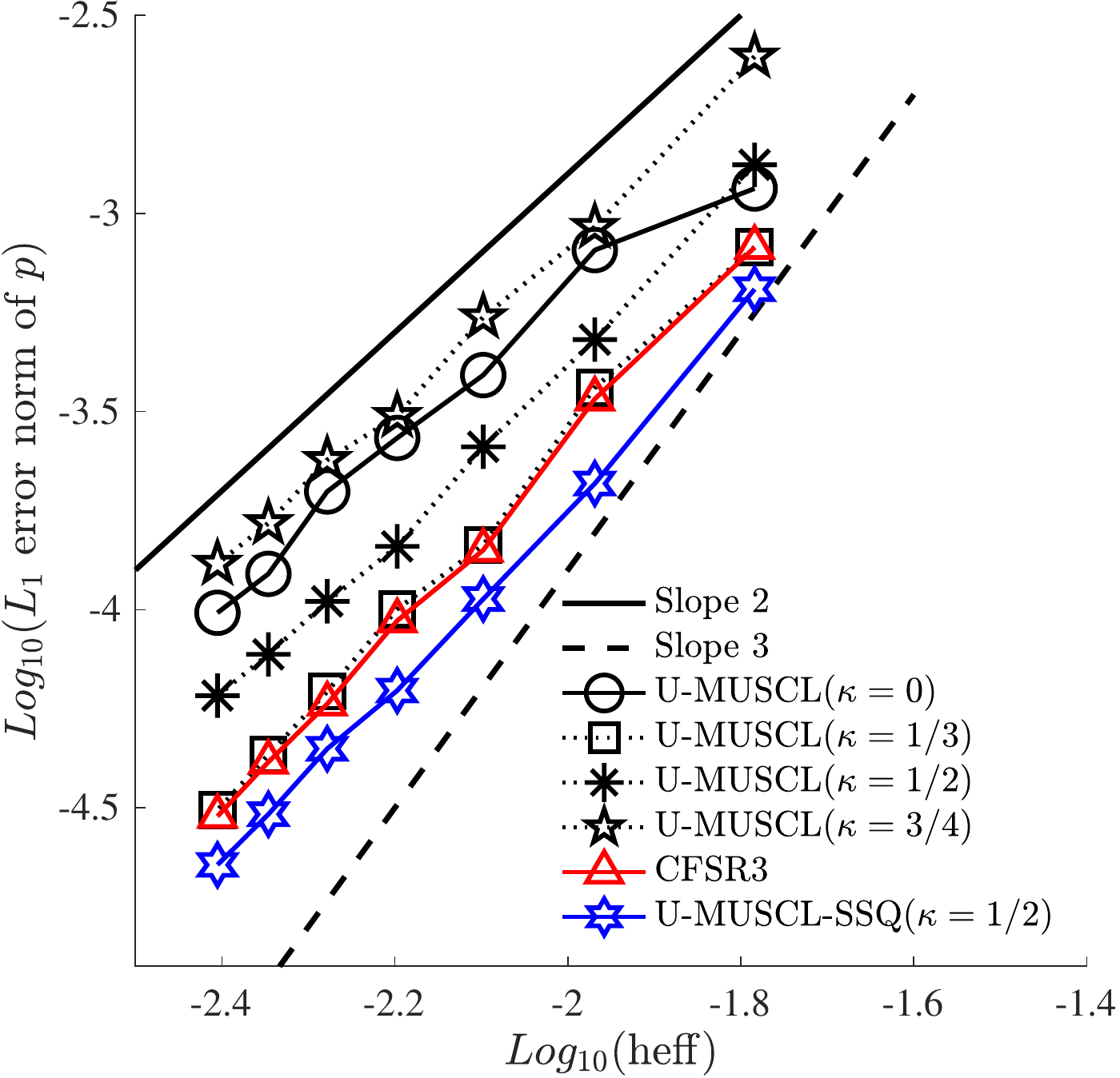}
  \caption[]{Error convergence for the pressure.}
\label{fig:twod_verification_tria_irrg_err}%
      \end{subfigure}
      \caption{
\label{fig:twod_verification_tria_irrg}%
Error convergence results for the Euler equations on irregular triangular grids.} 
\end{figure}

Next, we consider irregular triangular grids with 1024, 2304, 4096, 6400, 12544, and 16384 nodes. The coarsest grid is shown in Figure \ref{fig:twod_verification_tria_irrg_grid}. Error convergence results are shown in Figure \ref{fig:twod_verification_tria_irrg_err}. In this case, no schemes achieve third-order accuracy. However, CFSR3 and U-MUSCL-SSQ still provide accurate solutions compared with others except U-MUSCL($\kappa=1/3$), which is as accurate as CFSR3. It is interesting to note that U-MUSCL-SSQ with $\kappa_s=1/4$ gives the lowest errors for all triangular-grid cases.

\subsection{Unsteady vortex transport problem in two dimensions}
\label{results_2d_vorex}
 
\begin{figure}[htbp!]
    \centering
      \begin{subfigure}[t]{0.48\textwidth}
  \includegraphics[width=0.99\textwidth,trim=0 0 0 0 ,clip]{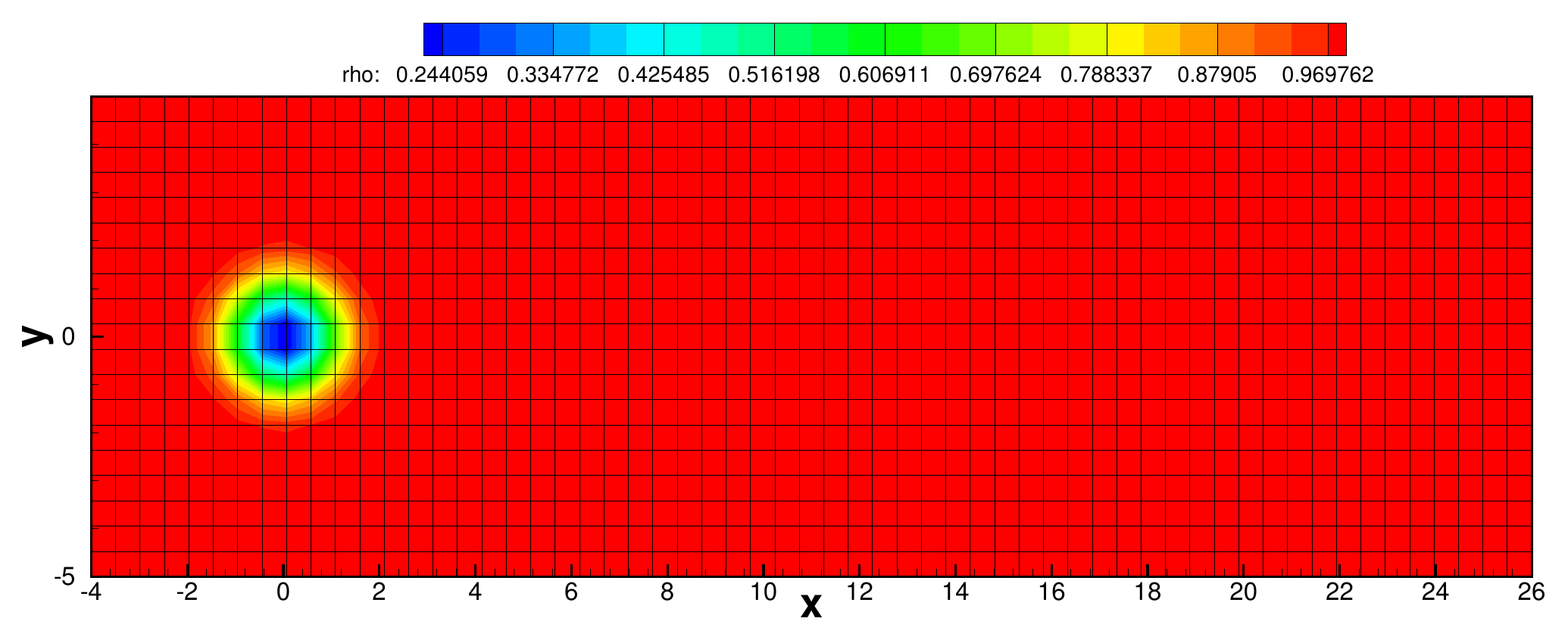}
  \caption[]{Quadrilateral grid and the initial density contours.}
\label{fig:twod_vortec_grid_quad}%
      \end{subfigure}
      \begin{subfigure}[t]{0.48\textwidth}
  \includegraphics[width=0.99\textwidth,trim=0 0 0 0 ,clip]{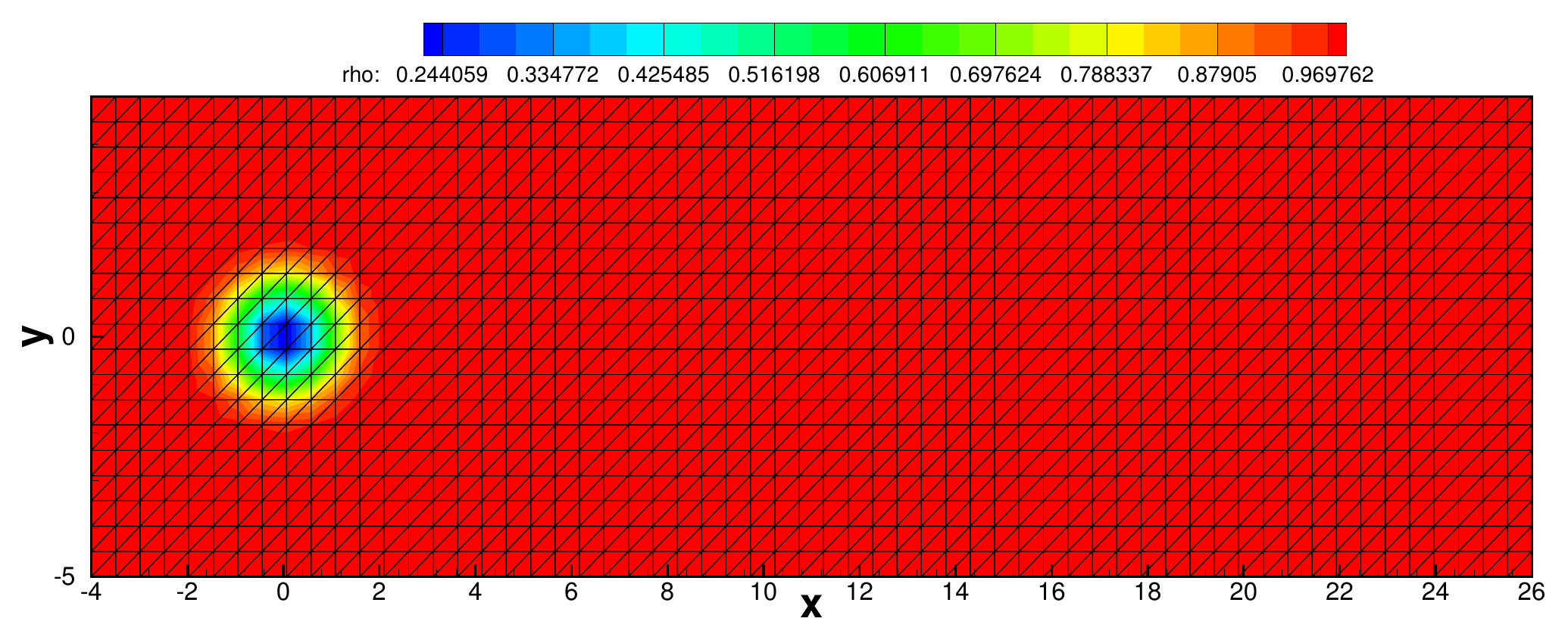}
  \caption[]{Triangular grid and the initial density contours.}
\label{fig:twod_vortec_grid_tria}%
      \end{subfigure}
      \begin{subfigure}[t]{0.48\textwidth}
  \includegraphics[width=0.99\textwidth,trim=0 0 0 0, clip]{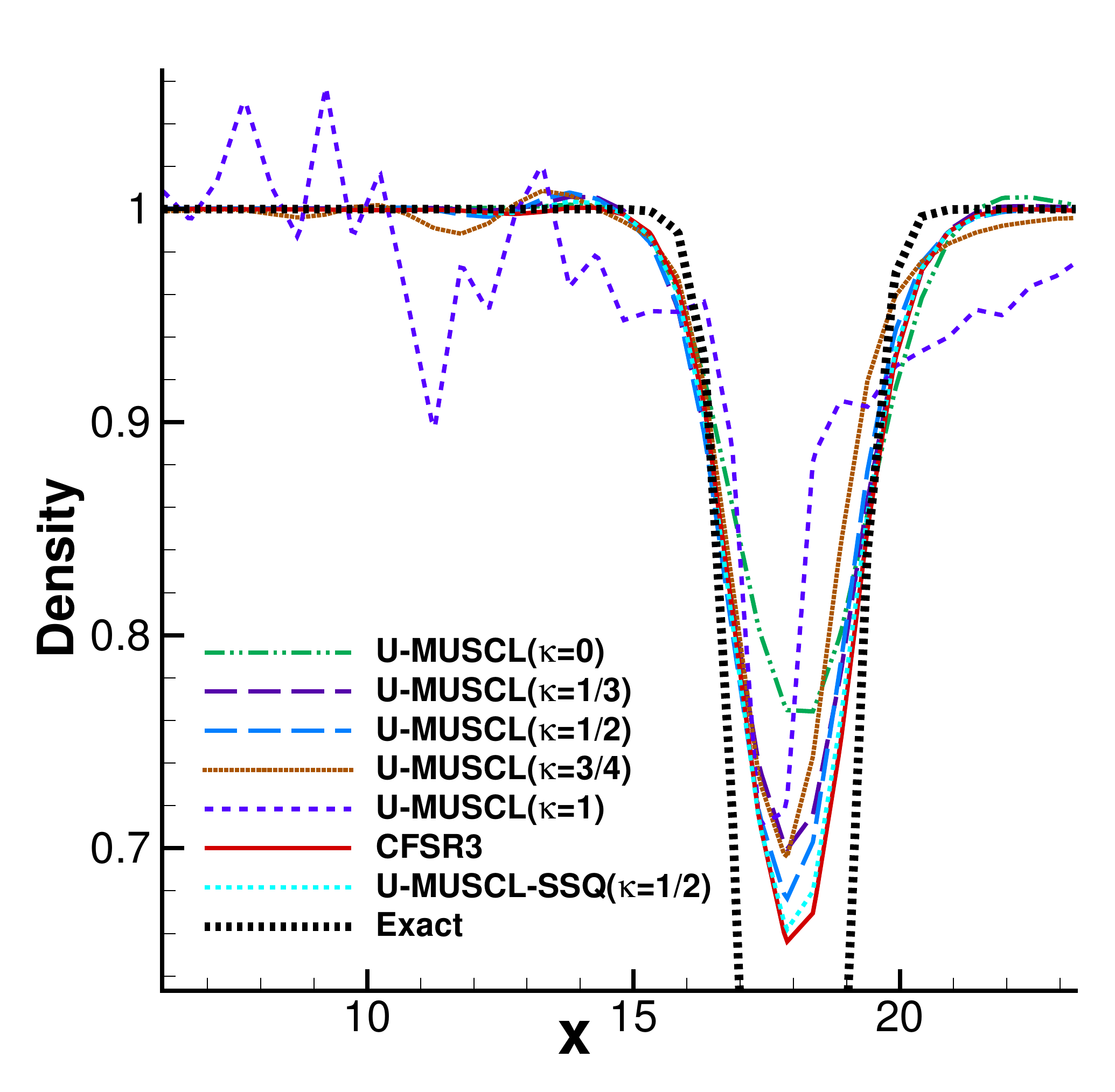}
  \caption[]{Section plot taken at $y=0$ (quadrilaterals).}
\label{fig:twod_vortex_section_quad}%
      \end{subfigure}
      \begin{subfigure}[t]{0.48\textwidth}
  \includegraphics[width=0.99\textwidth,trim=0 0 0 0, clip]{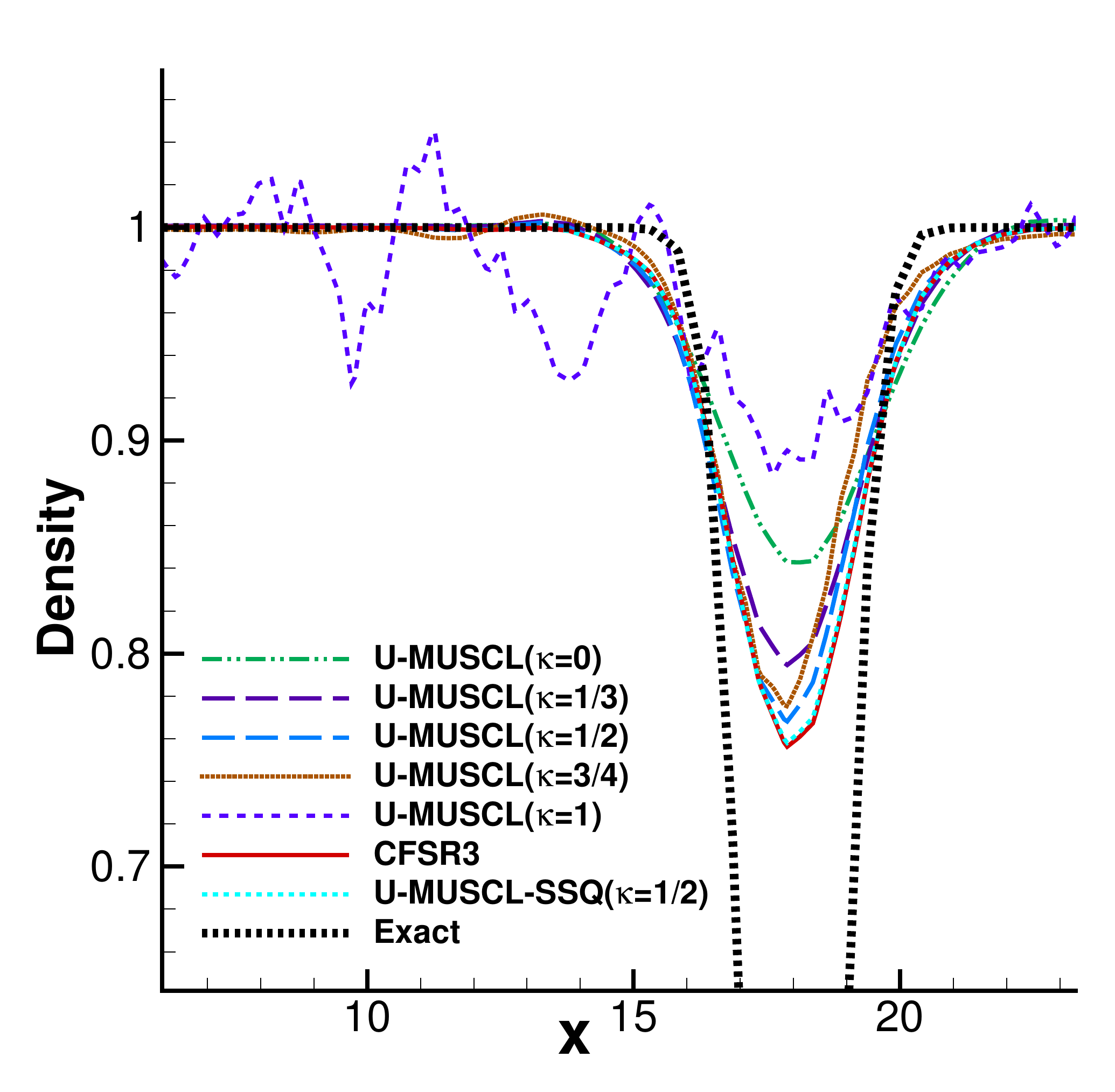}
  \caption[]{Section plot taken at $y=0$ (triangles).}
\label{fig:twod_vortex_section_tria}%
      \end{subfigure}
      \caption{
\label{fig:twod_veortex_quad}%
Vortex-transport problem on regular quadrilateral and triangular grids.} 
\end{figure}

To investigate the impact of various schemes on unsteady problems, we consider an inviscid vortex-transport problem, whose solution is an exact solution to the Euler equations without a source term:
\begin{eqnarray}
u =  u_\infty  -  \frac{ K \overline{y} }{2 \pi }   \exp \left(  \frac{1-\overline{r} ^2 }{2} \right) ,  \quad
v  = v_\infty   + \frac{ K \overline{x} }{2 \pi }   \exp \left(  \frac{1-\overline{r} ^2 }{2} \right),
\end{eqnarray}
and
\begin{eqnarray}
T =  1  -  \frac{ K^2 (\gamma-1)  }{8 \pi^2 }   \exp \left(  1-\overline{r} ^2   \right) , \quad
\rho = T^{  \frac{1}{\gamma-1}  } , \quad
p = \frac{ \rho ^{  \gamma   }  }{\gamma}    ,
\end{eqnarray}
where $\overline{x} = x  - u_\infty t$, $\overline{y} = y  - v_\infty t$, $\overline{r} ^2 =\overline{x}^2 + \overline{y}^2 $, and $(u_\infty,v_\infty)=(0.5,0.0)$. The initial solution at $t=0$ is shown in Figure \ref{fig:twod_vortec_grid_quad}. Here, $K$ is a key parameter. As discussed in Ref.[\citen{Nishikawa_FakeAccuracy:2020}], a small value of $K$ effectively linearizes the Euler equations and thus leads to false high-order accuracy of U-MUSCL schemes. To avoid this problem, we set $K=6$, which is large enough to avoid the problem. See Ref.[\citen{Nishikawa_FakeAccuracy:2020}] for further details. The solution is computed at the final time $t=36.0$ with the third-order SSP Runge-Kutta time-integration scheme with a fixed time step $\Delta t= 0.0005$. For this problem, we employ the Roe flux and the explicit SSP RK3 time-stepping scheme [\citen{SSP:SIAMReview2001}].

We consider a regular quadrilateral grid and also a regular triangular grid, both of which have 60$\times$20 nodes, as shown in Figures \ref{fig:twod_vortec_grid_quad} and \ref{fig:twod_vortec_grid_tria}, respectively. The solutions obtained with various schemes at the final time are compared in Figures \ref{fig:twod_vortex_section_quad} and \ref{fig:twod_vortex_section_tria}. For the quadrilateral grid, we observe that CFSR3 is the most accurate as expected. U-MUSCL-SSQ is also very accurate. This scheme requires the inversion of the mass matrix at each stage of the RK scheme; it took four Gauss-Seidel relaxations to reduce the linear residual by six orders of magnitude. Among the U-MUSCL schemes, we find that unlike the previous test cases, U-MUSCL($\kappa=1/2$) is more accurate than U-MUSCL($\kappa=1/3$) for this fully nonlinear problem. The slightly superior performance of $\kappa=1/2$ may be due to the less dissipative nature as discussed in Section \ref{quadratic_extrapolation}. However, the dissipation cannot be too small. The solution obtained with $\kappa=3/4$ has a larger dispersive error and the peak is not better maintained than $\kappa=1/2$. Also, the central scheme ($\kappa=1$) suffers from severe dispersive errors as can be clearly seen. Note that the third-order RK scheme is stable with the central discretization 
($\kappa=1$) and therefore the oscillation is due to the dispersive error and is not a stability problem. The U-MUSCL scheme of $\kappa=0$ is dissipative and 
dispersive, and therefore the least accurate option. 

\begin{figure}[htbp!]
    \centering
      \begin{subfigure}[t]{0.32\textwidth}
  \includegraphics[width=0.99\textwidth,trim=0 0 0 0 ,clip]{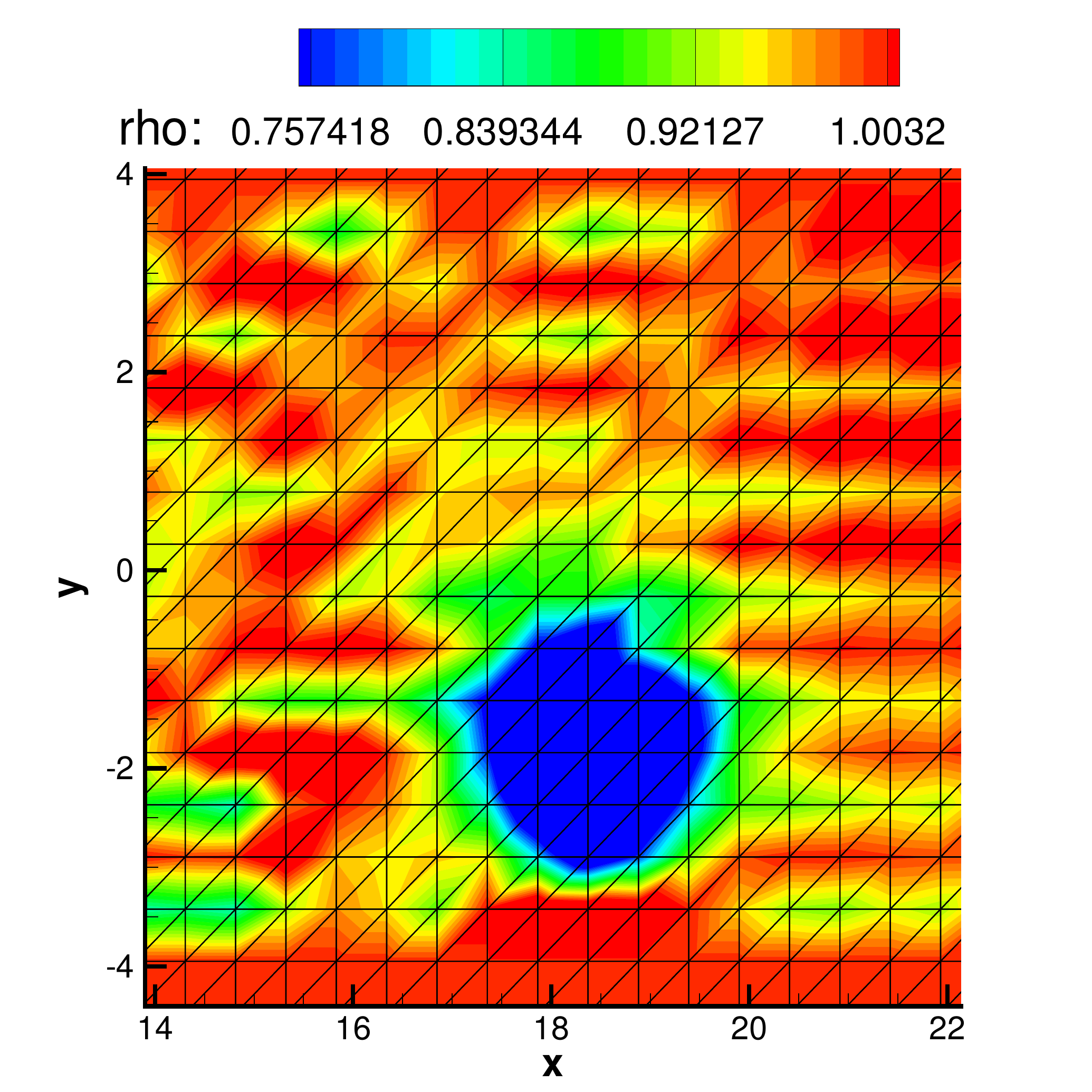}
  \caption[]{U-MUSCL ($\kappa=1$)}
\label{fig:twod_veortex_tria_top_view_1p0}%
      \end{subfigure}
      \begin{subfigure}[t]{0.32\textwidth}
  \includegraphics[width=0.99\textwidth,trim=0 0 0 0 ,clip]{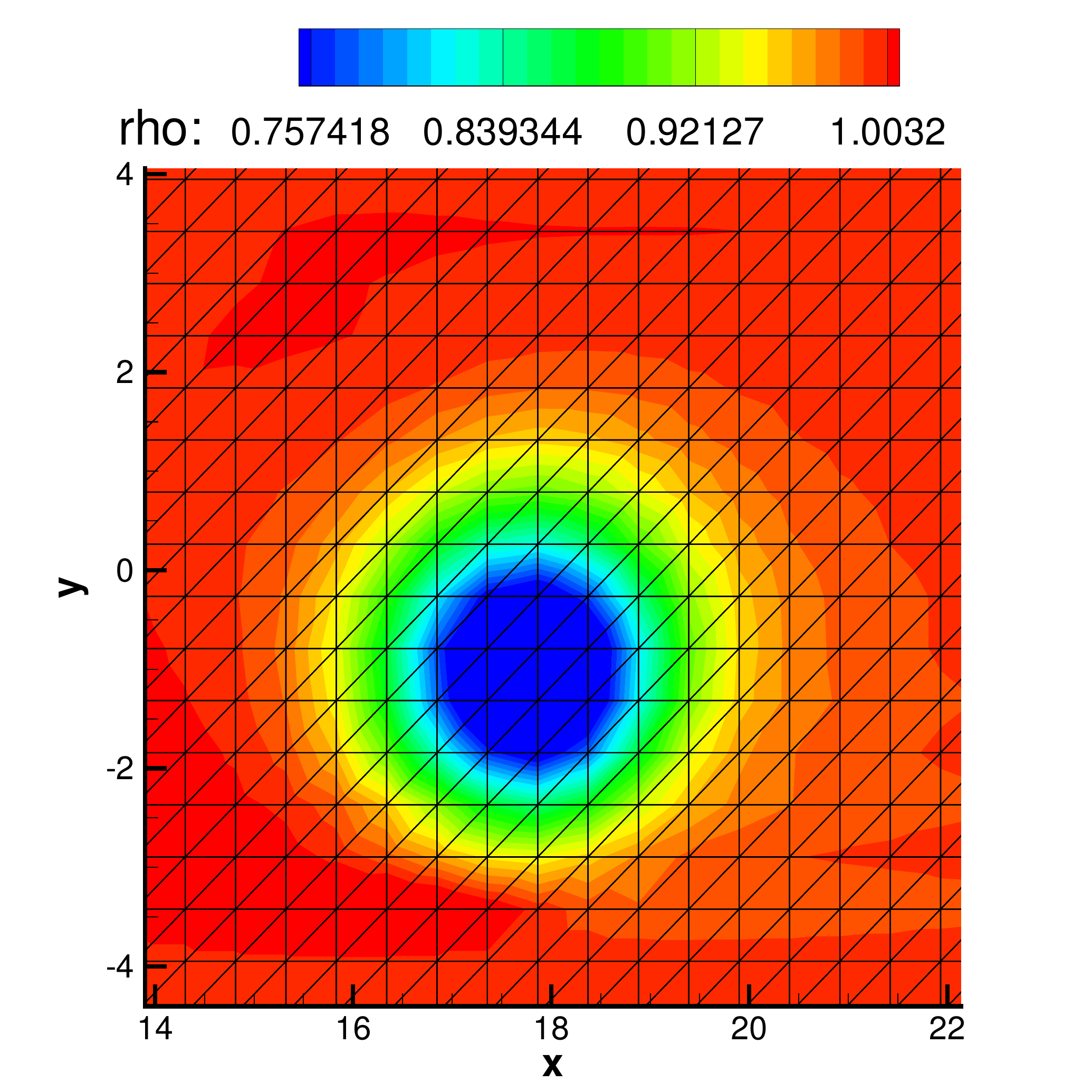}
  \caption[]{U-MUSCL ($\kappa=3/4$)}
\label{fig:twod_veortex_tria_top_view_1o3}%
      \end{subfigure}
      \begin{subfigure}[t]{0.32\textwidth}
  \includegraphics[width=0.99\textwidth,trim=0 0 0 0 ,clip]{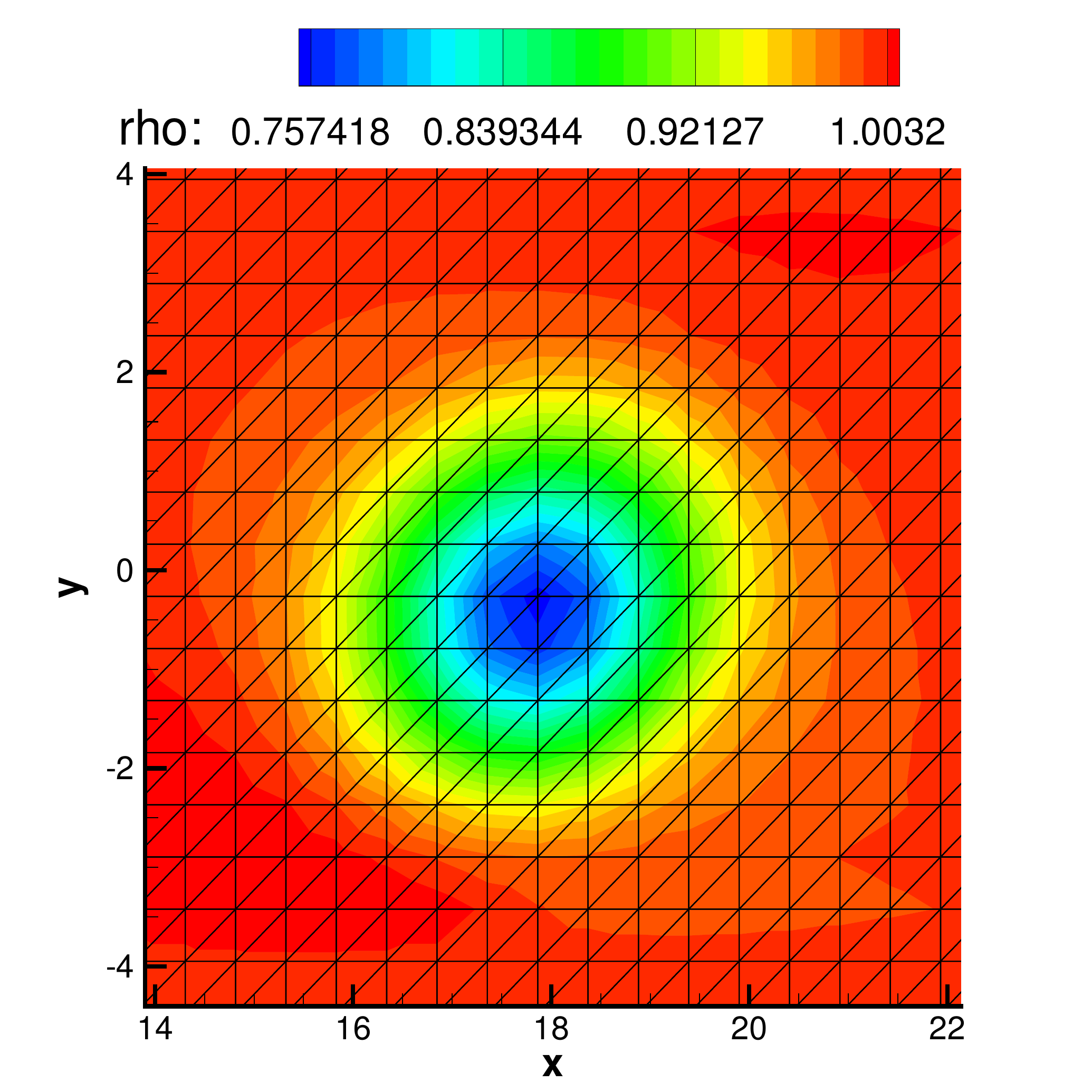}
  \caption[]{U-MUSCL ($\kappa=1/2$)}
\label{fig:twod_veortex_tria_top_view_1o2}%
      \end{subfigure}
      \begin{subfigure}[t]{0.32\textwidth}
  \includegraphics[width=0.99\textwidth,trim=0 0 0 0 ,clip]{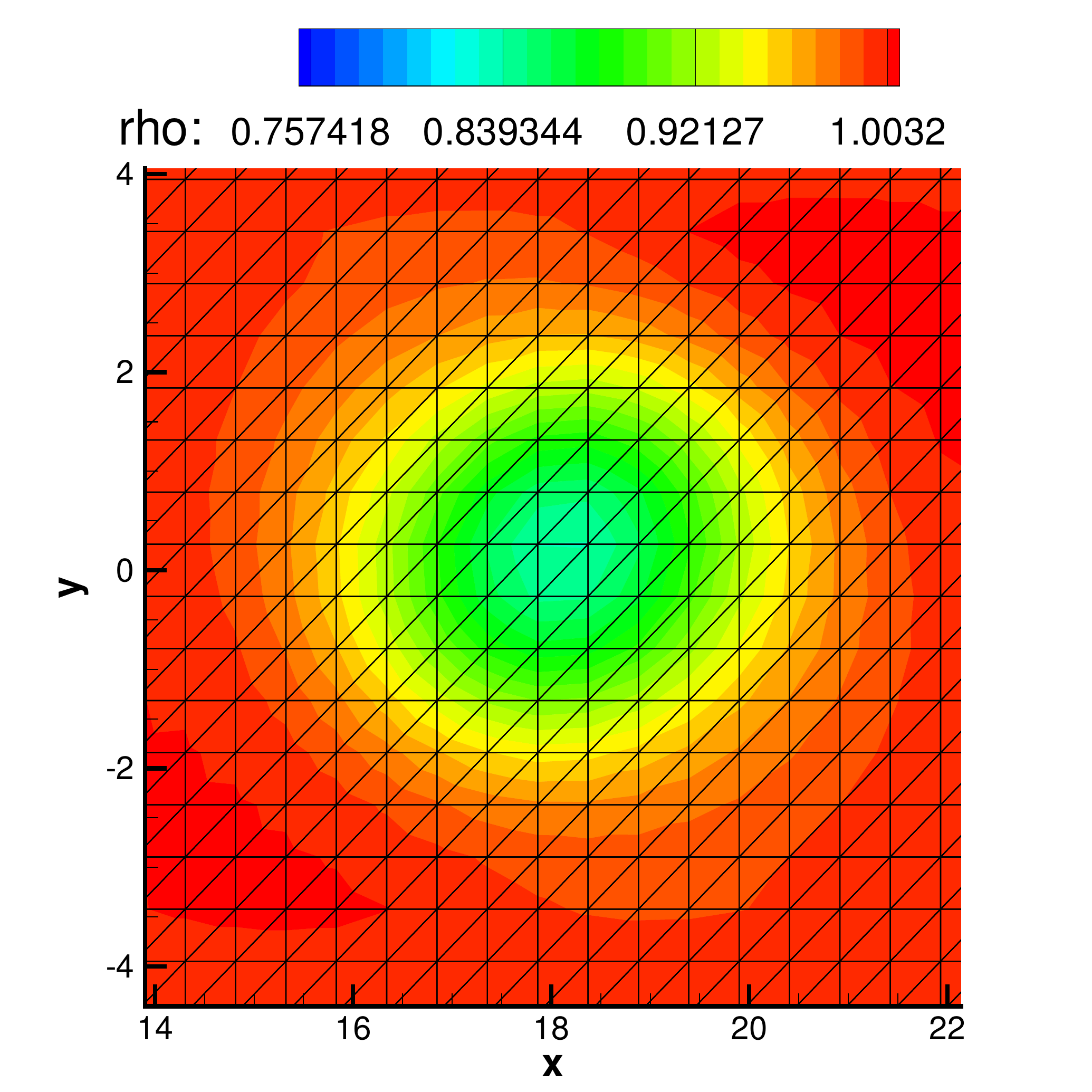}
  \caption[]{U-MUSCL ($\kappa=0$)}
\label{fig:twod_veortex_tria_top_view_0p0}%
      \end{subfigure}
      \begin{subfigure}[t]{0.32\textwidth}
  \includegraphics[width=0.99\textwidth,trim=0 0 0 0, clip]{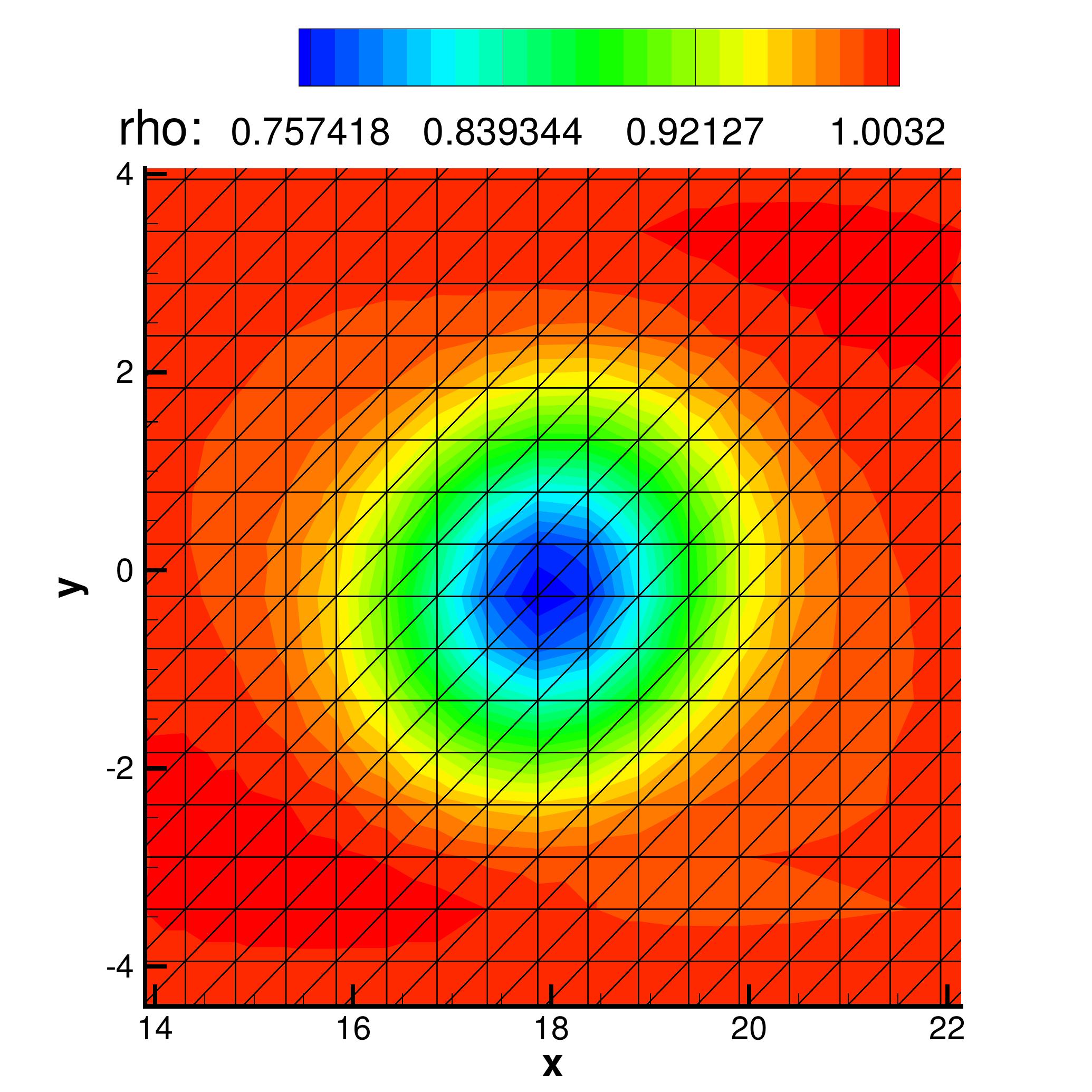}
  \caption[]{CFSR3}
\label{fig:twod_veortex_tria_top_view_cfsr3}%
      \end{subfigure}
      \begin{subfigure}[t]{0.32\textwidth}
  \includegraphics[width=0.99\textwidth,trim=0 0 0 0, clip]{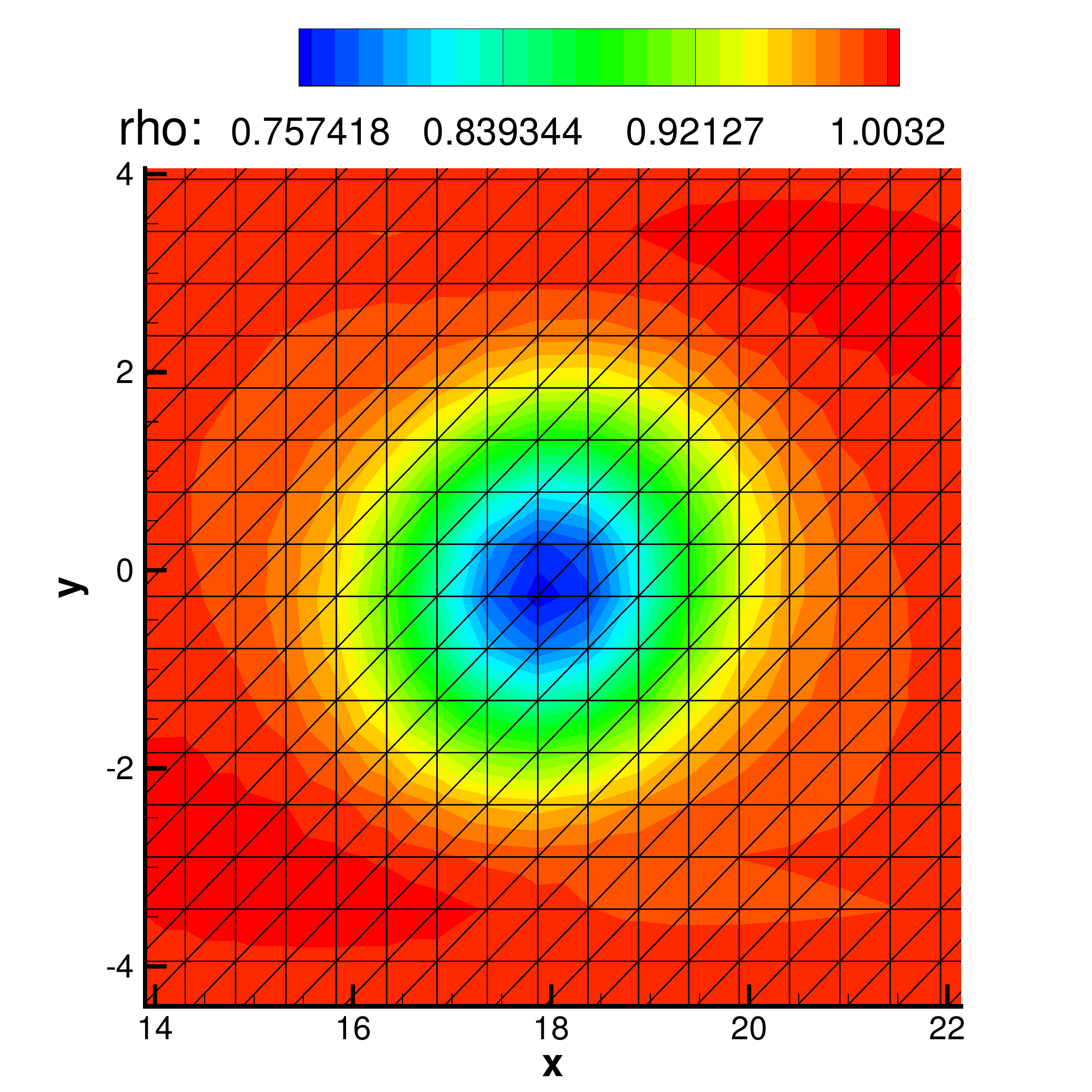}
  \caption[]{U-MUSCL-SSQ}
\label{fig:twod_veortex_tria_top_view_ssq}%
      \end{subfigure}
      \caption{
\label{fig:twod_veortex_tria_top_view}%
Contours of the vortex at the final time on triangular grids, shown with the same contour range. } 
\end{figure}

The solutions on the triangular grid exhibit the same trend, which are more dissipative than those on the regular quadrilateral grid.
See Figure \ref{fig:twod_vortex_section_tria}. 
The dispersive error is significant for U-MUSCL($\kappa=3/4$) and U-MUSCL($\kappa=1$). On the other hand, U-MUSCL($\kappa=1/2$) and 
U-MUSCL($\kappa=1/3$) have a better dispersion property (close to the third-order scheme) although they are theoretically second-order accurate. 
Compare it with CFSR3 and U-MUSCL-SSQ. Again, U-MUSCL($\kappa=1/2$) keeps the peak better than $\kappa=1/3$. 
Note that U-MUSCL($\kappa=0$) is again dissipative and U-MUSCL($\kappa=1$) may look also dissipative. However, it is actually due to a large dispersive error.
To illuminate this, we show top-views of the final solutions in Figure \ref{fig:twod_veortex_tria_top_view}. The vortex initially located at the origin
 is simply convected to the right; therefore it must stay along the centerline and be centered at $(x,y)=(18, 0)$ at the final time. 
 However, it can be seen that the vortex is 
 located significantly below the centerline in U-MUSCL($\kappa=1$). See Figure \ref{fig:twod_veortex_tria_top_view_1p0}. 
 The same is observed for U-MUSCL($\kappa=3/4$) as in Figure \ref{fig:twod_veortex_tria_top_view_1o3}. 
 The error is small with U-MUSCL($\kappa=1/2$), and slightly smaller with CFSR3 and U-MUSCL-SSQ, as shown in Figures 
 \ref{fig:twod_veortex_tria_top_view_1o2}, \ref{fig:twod_veortex_tria_top_view_cfsr3}, and \ref{fig:twod_veortex_tria_top_view_ssq}, respectively.
 Finally, we see from Figure \ref{fig:twod_veortex_tria_top_view_0p0} that U-MUSCL($\kappa=0$) is indeed rather dissipative than dispersive: 
 the vortex nearly at the right location but severely dissipated. Comparing it with Figure \ref{fig:twod_veortex_tria_top_view_1p0}, we see that the 
 apparently dissipative solution of U-MUSCL($\kappa=1$), as in Figure \ref{fig:twod_vortex_section_tria}, is due to the dispersive error whereas it is 
 due to the dissipative error in the case of U-MUSCL($\kappa=0$).


In both quadrilateral and triangular grids, U-MUSCL($\kappa=0$) and U-MUSCL($\kappa=1$) are significantly less accurate than others and thus 
U-MUSCL($\kappa=1/2$), U-MUSCL-SSQ, and CFSR3 are recommended. However, improvements of U-MUSCL-SSQ and CFSR3
over others are somewhat minor. As we will show in a subsequent paper, more significant improvements can be achieved by their fourth-order versions, which 
reduces not only dispersive errors but also dissipative errors. 

\subsection{Unsteady separated laminar flow over a cylinder}
\label{results_3d_laminar_cylinder}
 
To investigate the genuine third-order accuracy for practical three-dimensional computations, we 
implemented the CFSR3 scheme in NASA's FUN3D code [\citen{fun3d_website}]. In 3D, as a genuine third-order scheme, 
we will focus on CFSR3 rather than U-MUSCL-SSQ because $\kappa$ is a free parameter in CFSR3 (not in U-MUSCL-SSQ), and can be
adjusted to control the amount of dissipation, and also because CFSR3 is third-order accurate on both 
regular hexahedral and tetrahedral grids while U-MUSCL-SSQ is not third-order accurate for regular hexahedral grids. 

\begin{figure}[htbp!]
    \centering
      \begin{subfigure}[t]{0.30\textwidth}
  \includegraphics[width=0.99\textwidth,trim=0 0 0 0 ,clip]{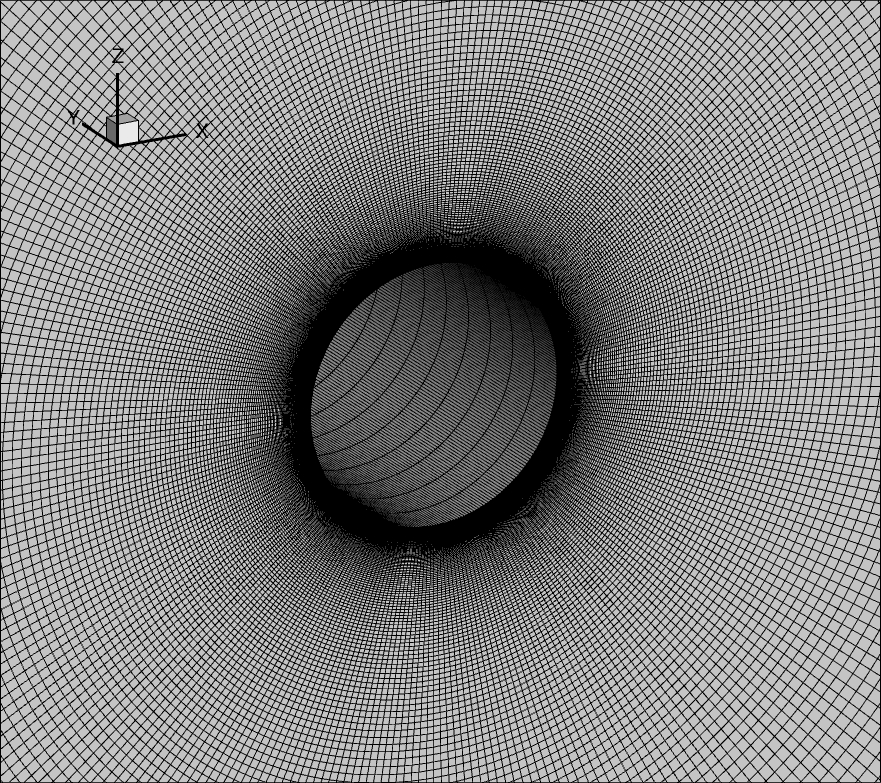}
  \caption[]{Grid close to the body}
      \end{subfigure}
      \begin{subfigure}[t]{0.30\textwidth}
  \includegraphics[width=0.99\textwidth,trim=0 0 0 0 ,clip]{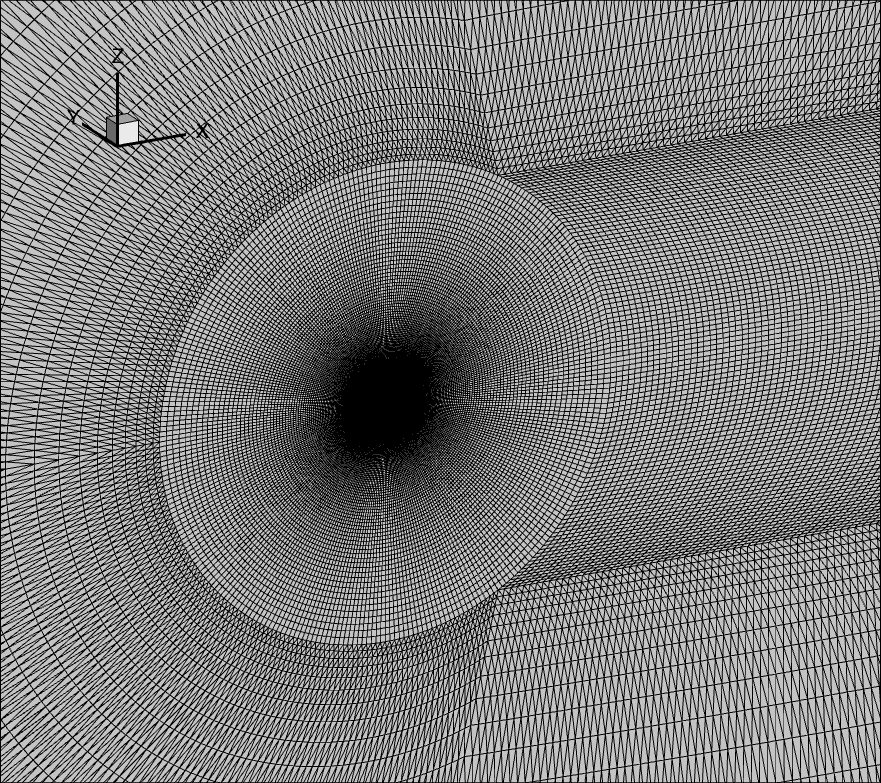}
  \caption[]{Grid in wake region}
      \end{subfigure}
      \begin{subfigure}[t]{0.30\textwidth}
  \includegraphics[width=0.99\textwidth,trim=0 0 0 0 ,clip]{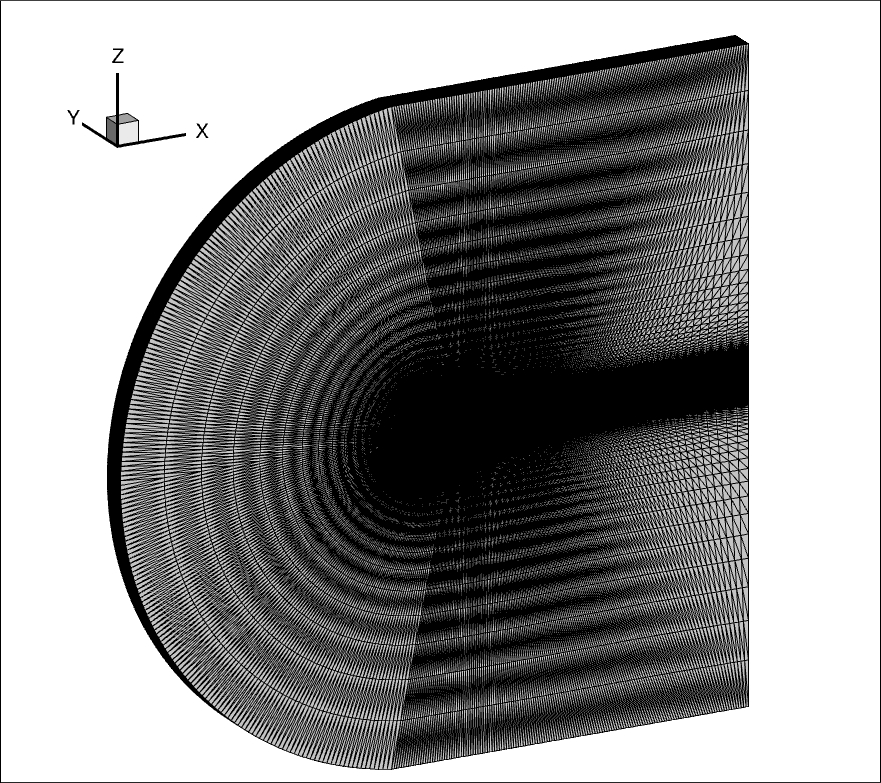}
  \caption[]{Grid in full domain}
      \end{subfigure}
      \caption{
\label{fig:cyl_tet_mesh}
A hexahedral-prismatic mixed grid for the laminar cylinder test case}
\end{figure}

We tested  the CFSR3 and U-MUSCL schemes for a laminar flow over a cylinder with the Mach number $0.1$ at zero
angle of attack, the Reynolds number $3,900$, and the free stream temperature $460.0[R]$. The cylinder of diameter 1.0
has a span-wise length of 2.0 and the outer boundary is located approximately at the distance 100 from the center of the cylinder. 
The domain is taken as periodic at two planes with the minimum and maximum $y$-coordinates. 
The grid used for these simulations is a mixed grid with prismatic and hexahedral cells as shown in Figure \ref{fig:cyl_tet_mesh}.
 It has 3,934,800 nodes, 1,175,000 prisms, and 3,243,000 hexahedra.
 Simulations were run with the Roe flux, the Galerkin viscous discretization, and the second-order backward-difference 
 (BDF2) time-stepping scheme with a non-dimensionalized time step of $0.05$ for 10,000 time steps, with 400 cores. 
 At each physical time step, we perform 20 nonlinear iterations at maximum. 
In this problem, we compare CFSR3 of $\kappa=1/3$ and $\kappa=1/2$ with U-MUSCL of $\kappa= 0$, $1/3$, and $1/2$. 
 

\begin{figure}[htbp!]
    \centering
      \begin{subfigure}[t]{0.48\textwidth}
  \includegraphics[width=0.99\textwidth,trim=0 0 0 0 ,clip]{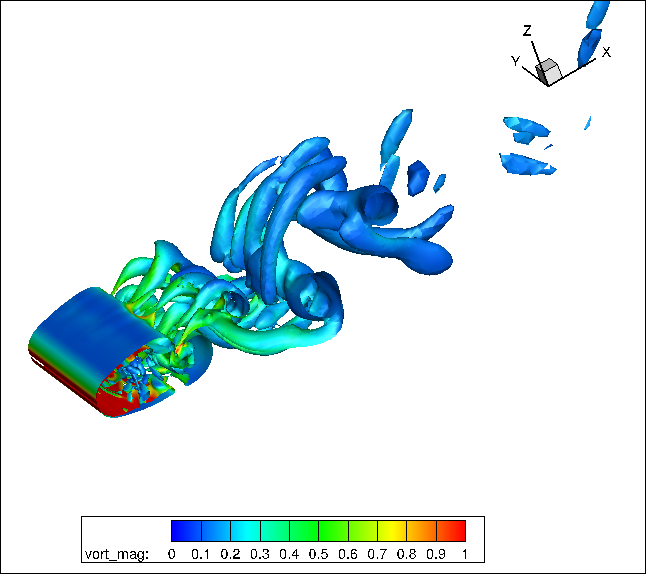}
  \caption[]{U-MUSCL ($\kappa = 0$)}
\label{fig:unsteady_kappa0}%
      \end{subfigure}
      \\
      \begin{subfigure}[t]{0.48\textwidth}
  \includegraphics[width=0.99\textwidth,trim=0 0 0 0 ,clip]{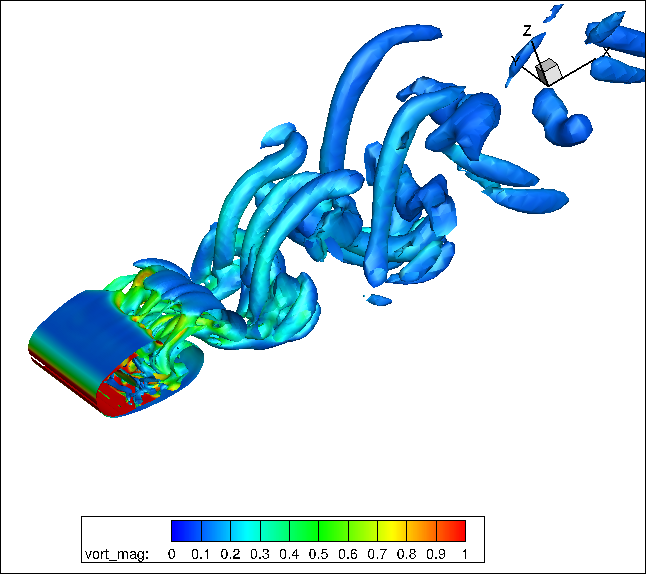}
  \caption[]{U-MUSCL ($\kappa = 1/3$)}
\label{fig:unsteady_kappa13}%
      \end{subfigure}
      \hfill
      \begin{subfigure}[t]{0.48\textwidth}
  \includegraphics[width=0.99\textwidth,trim=0 0 0 0 ,clip]{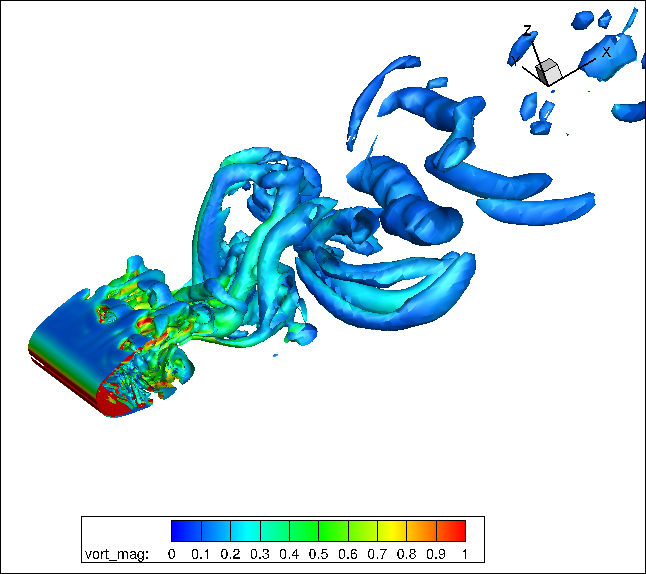}
  \caption[]{U-MUSCL ($\kappa = 1/2$)}
\label{fig:unsteady_kappa12}%
      \end{subfigure}
      \hfill
      \begin{subfigure}[t]{0.48\textwidth}
  \includegraphics[width=0.99\textwidth]{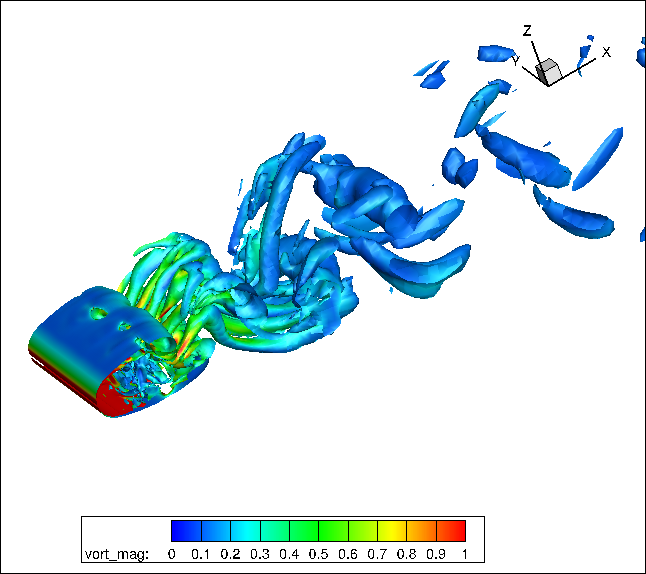}
\caption{ CFSR ($\theta = 1/3$, $\kappa=1/3$) 
\label{fig:unsteady_CFSR_k1o3}}
      \end{subfigure}  
      \hfill
\begin{subfigure}[t]{0.48\textwidth}
  \includegraphics[width=0.99\textwidth]{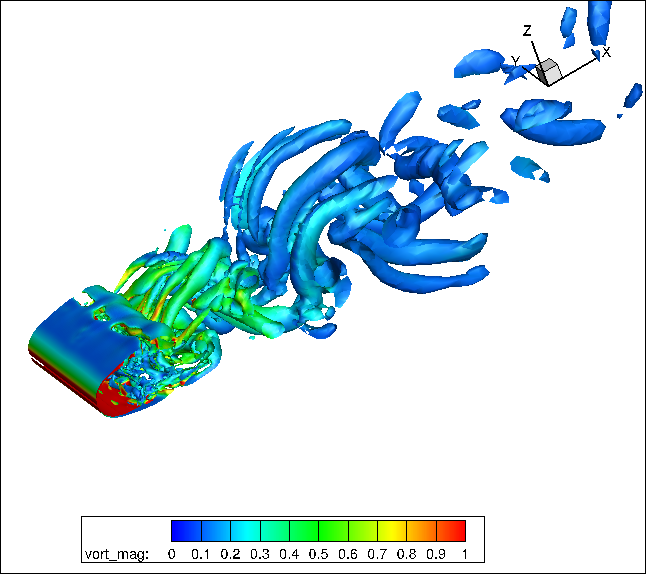}
\caption{ CFSR ($\theta = 1/3$, $\kappa=1/2$) 
\label{fig:unsteady_CFSR_k1o2}}
      \end{subfigure}
      \caption{
\label{fig:kappa_unsteady_cylinder_far}%
Unsteady separated laminar flow results.} 
\end{figure}

The Q-isosurface (colored by the vorticity magnitude) plots are compared in Figure \ref{fig:kappa_unsteady_cylinder_far}. 
Figures \ref{fig:unsteady_kappa0}, \ref{fig:unsteady_kappa13}, and \ref{fig:unsteady_kappa12} show that the resolution 
improves as the dissipation is reduced from $\kappa=0$ to $\kappa=1/2$. It is important to note that U-MUSCL of $\kappa=1/2$ is not third-order 
scheme but a low-dissipation second-order scheme. In Refs.[\citen{yang_harris:AIAAJ2016,yang_harris:CCP2018}], U-MUSCL of $\kappa=1/2$ is
considered as a third-order scheme apparently following Burg's paper [\citen{burg_umuscl:AIAA2005-4999}]. However, numerical 
results in Refs.[\citen{yang_harris:AIAAJ2016,yang_harris:CCP2018}] show second-order error convergence for U-MUSCL of $\kappa=1/2$. 
These numerical results are correct: it is not a third-order scheme as we discussed in Section \ref{1d_fd_fr}. But it provides low dissipation 
simply because it is closer to $\kappa=1$ (zero dissipation) than $\kappa=0$ and $\kappa=1/3$. 

\begin{figure}[htbp!]
    \centering
      \begin{subfigure}[t]{0.48\textwidth}
  \includegraphics[width=0.99\textwidth,trim=0 0 0 0 ,clip]{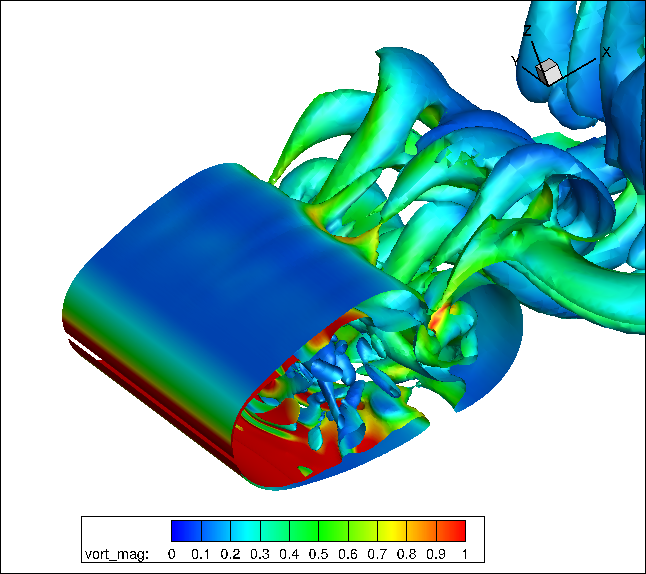}
  \caption[]{U-MUSCL ($\kappa = 0$)}
\label{fig:unsteady_kappa0_zoom}%
      \end{subfigure}
      \\
      \begin{subfigure}[t]{0.48\textwidth}
  \includegraphics[width=0.99\textwidth,trim=0 0 0 0 ,clip]{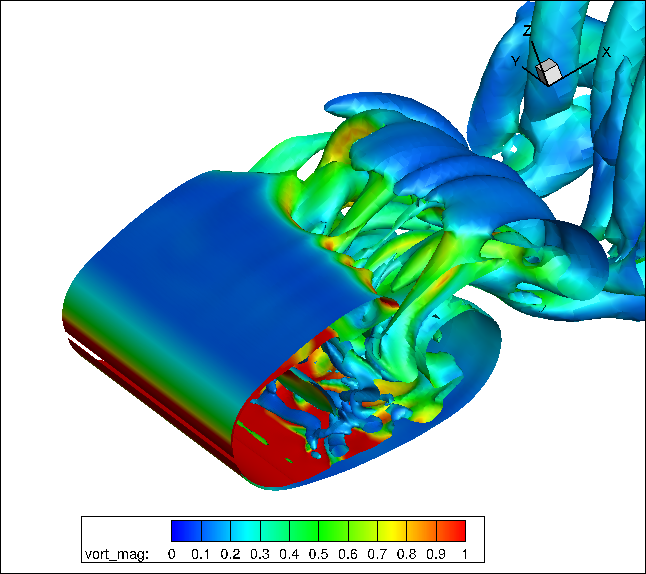}
  \caption[]{U-MUSCL ($\kappa = 1/3$)}
\label{fig:unsteady_kappa13_zoom}%
      \end{subfigure}
      \hfill
      \begin{subfigure}[t]{0.48\textwidth}
  \includegraphics[width=0.99\textwidth,trim=0 0 0 0 ,clip]{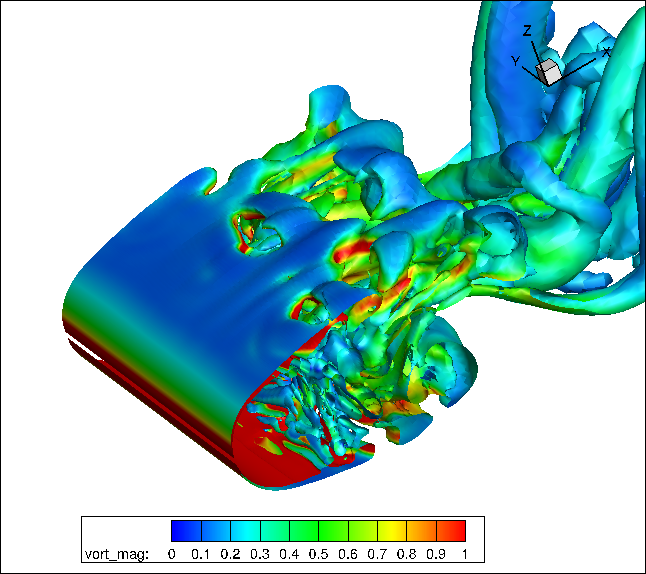}
  \caption[]{U-MUSCL ($\kappa = 1/2$)}
\label{fig:unsteady_kappa12_zoom}%
      \end{subfigure}
      \hfill
      \begin{subfigure}[t]{0.48\textwidth}
  \includegraphics[width=0.99\textwidth]{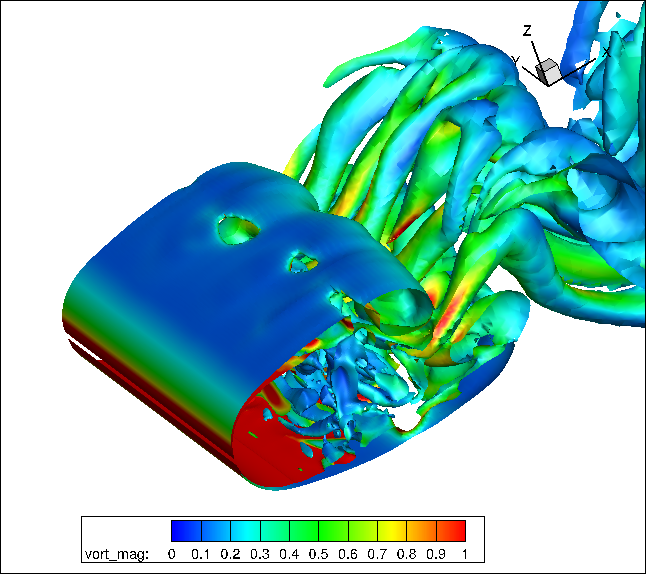}
\caption{CFSR ($\theta = 1/3$, $\kappa=1/3$) 
\label{fig:unsteady_CFSR_k1o3_zoom}}
      \end{subfigure}  
      \hfill
\begin{subfigure}[t]{0.48\textwidth}
  \includegraphics[width=0.99\textwidth]{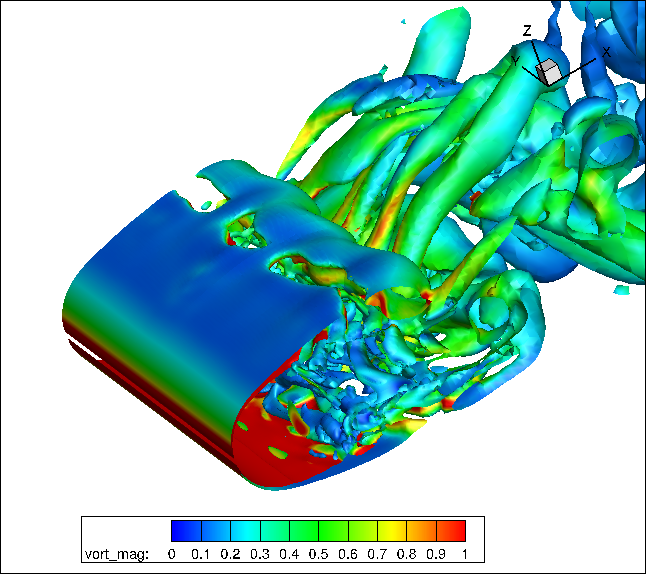}
\caption{ CFSR ($\theta = 1/3$, $\kappa=1/2$) 
\label{fig:unsteady_CFSR_k1o2_zoom}}
      \end{subfigure}
      \caption{
\label{fig:kappa_unsteady_cylinder_zoom}%
Zoomed-in view of unsteady separated laminar flow results.} 
\end{figure}

Results for CFSR3 are shown in Figures \ref{fig:unsteady_CFSR_k1o3} and \ref{fig:unsteady_CFSR_k1o2}. 
Comparing them for the same $\kappa$ value, we see that CFSR3 gives slightly better resolved results, especially for $\kappa=1/3$. 
As discussed earlier in Sections \ref{oned_diss_disp} and \ref{results_2d_vorex}, a third-order scheme eliminates the leading dispersion 
error from a second-order scheme of the same family and thus it does not directly impact on dissipation, but it can result in a more accurate solution. 
These features can be seen more clearly in the zoomed-in views provided in Figure \ref{fig:kappa_unsteady_cylinder_zoom}. 
Significantly less dissipative solutions are expected with fourth-order schemes, which will be reported in a subsequent paper. 

\section{Concluding Remarks}
\label{conclusions}
 
In this paper, we have resolved confusions over third-order accuracy of the U-MUSCL scheme in the edge-based discretization. In one dimension, the U-MUSCL scheme can be third-order in the following cases: 
\begin{enumerate}
    \item Finite-volume scheme with cell-averaged solutions and source terms with $\kappa=1/3$, 
    \item Finite-volume scheme with point-valued solutions and cell-averaged source terms with $\kappa=1/2$, 
    \item Finite-difference scheme with point-valued solutions and source terms with $\kappa=1/3$ (for linear equations only).
\end{enumerate}
With a single flux evaluation per edge (i.e., without high-order flux quadrature), the first two schemes do not extend to higher dimensions. 
The third one achieves third-order accuracy on regular grids in higher dimensions, but only for linear equations. In conclusion, the U-MUSCL 
scheme as typically implemented in practical CFD codes is second-order accurate at best even on regular grids for realistic applications. 

 In accuracy verification by the method of manufactured solutions, the U-MUSCL scheme was found more accurate with $\kappa=1/3$ than 
 with $\kappa=1/2$. This is shown to be due to a perturbed form of the exact solution, which (as verified theoretically in 
 Ref.[\citen{Nishikawa_FakeAccuracy:2020}]) effectively linearizes the Euler equations with small perturbation and the fact
 that $\kappa=1/3$ gives third-order accuracy for linear equations. On the other hand, for an inviscid vortex transport 
 problem, we obtained slightly less dissipative solutions with $\kappa=1/2$ than with $\kappa=1/3$; apparently, the superior 
 performance is due to the lower dissipation of $\kappa=1/2$ than $\kappa=1/3$. Further reducing the dissipation with $\kappa > 1/2$ did not 
 improve the solution due to large dispersive errors. These results indicate that $\kappa=1/3$ is suitable for problems with small 
 perturbations and $\kappa=1/2$ otherwise. It is also noted that $\kappa=1/2$ tends to provide more robust iterative convergence on unstructured 
 grids as it reduces the effect of 
 LSQ gradients in the U-MUSCL solution reconstruction scheme, more than $\kappa=1/3$ does. This could be another major reason 
 for the popularity of $\kappa=1/2$. Finally, the U-MUSCL scheme can achieve third-order accuracy with $\kappa=1/2$ for 
 a steady nonlinear problem without source terms, but only on regular triangular/tetrahedral grids. As we have shown, it corresponds to 
 a special third-order U-MUSCL-SSQ scheme derived in this paper. This can be yet another reason that $\kappa=1/2$ has been found 
 advantageous over other values. 
 
Finally, we discussed two techniques to achieve genuine third-order accuracy for general nonlinear equations when the grid is regular: 
\begin{enumerate}
    \item CFSR3: Conservative finite-difference scheme with point-valued solutions and the U-MUSCL flux reconstruction with $\theta=1/3$ 
    and solution reconstruction with an arbitrary choice of $\kappa$, as originally developed in Ref.[\citen{Nishikawa_FakeAccuracy:2020}]. 
    The parameter $\kappa$ can be chosen to reduce the dissipation; we employed $\kappa=1/2$ to keep dissipation relatively low.
    This scheme is called FSR and an efficient version based on the chain rule is called CFSR3. Third-order accuracy is obtained on regular 
    quadrilateral and triangular grids.
    \item U-MUSCL-SSQ: Conservative finite-difference scheme with point-valued solutions with no flux reconstruction but with a special source discretization at $\kappa=1/2$. This scheme is called U-MUSCL-SSQ. Third-order accuracy is obtained on regular triangular/tetrahedral grids.
This is a useful option if one's target application allows the use of regular triangular/tetrahedral grids. 
\end{enumerate}
The more flexible scheme among the two, CFSR3, has been implemented in NASA's FUN3D code and superior resolution has been demonstrated for 
realistic three-dimensional flow problems. 


In future work, we will extend the CFSR3 scheme to fourth-order accuracy in FUN3D by extending the U-MUSCL solution reconstruction 
as originally proposed in Ref.[\citen{yang_harris:AIAAJ2016}] and also extending the efficient flux reconstruction. As the fourth-order scheme 
eliminates the leading dissipative error
of CFSR3, significant improvements in resolution are expected. It will be compared also with a second-order scheme with a 
low-dissipation numerical flux scheme [\citen{nishikawa_liu_aiaa2018-4166}]. 


\addcontentsline{toc}{section}{Acknowledgments}
\section*{Acknowledgments}

 The authors gratefully acknowledges support from the U.S. Army Research Office 
under the contract/grant number W911NF-19-1-0429 with Dr. Matthew Munson as the program manager. 
 The second author acknowledges support also from Software CRADLE, part of Hexagon and 
 the Hypersonic Technology Project, through the Hypersonic Airbreathing Propulsion Branch of the NASA Langley Research Center, under Contract No. 80LARC17C0004. The authors are grateful to Dr. Yi Liu (National Institute of Aerospace) for his assistance in the setup of three-dimensional test cases.

\addcontentsline{toc}{section}{References}
\bibliography{./bibtex_nishikawa_database}
\bibliographystyle{aiaa} 


\clearpage
\newpage

\section*{Appendix A: Truncation Error Analysis for U-MUSCL-SSQ on Regular Grids}
\label{appendixA}
\renewcommand{\thefigure}{A.\arabic{figure}}\setcounter{figure}{0}
\renewcommand{\thetable}{A.\arabic{table}}\setcounter{table}{0}
\renewcommand{\theequation}{A.\arabic{equation}}\setcounter{equation}{0}
\renewcommand{\thesubsection}{A.\arabic{subsection}}\setcounter{subsection}{0}

We follow Ref.[\citen{nishikawa_liu_source_quadrature:jcp2017}] and derive the truncation error of the U-MUSCL-SSQ scheme on regular simplex-element (triangular and tetrahedral) grids. As discussed in Ref.[\citen{nishikawa_liu_source_quadrature:jcp2017}], a regular grid is defined as a grid having an identical compact stencil (i.e., the stencil composed of edge-connected neighbors) at all nodes. The stencil is symmetric with respect to the central node and invariant under translation, i.e., for any edge-connected nodes $j$ and $k$, the stencil at $j$ translated along the edge $[jk]$ coincides with the stencil centered at $k$. Examples are a regular right triangular grid in Figure \ref{fig:grid_reg_tria} (or Figure \ref{fig:twod_verification_tria_right_grid}), a regular equilateral triangular grid in Figure \ref{fig:twod_verification_tria_equi_grid}, and a regular tetrahedral grid in Figure \ref{fig:grid_reg_tet}. Note that these grids are regular only away from the boundary. Regular simplex grids remain regular under any linear transformation (scaling, rotation, translation), and therefore there can be highly anisotropic and/or highly skewed regular grids. 

For simplicity but without loss of generality, we consider a scalar conservation law, $\mbox{div} \, {\bf f} = 0$, where the flux vector ${\bf f}$ is a function of a solution variable ${u}$. Note that ${\bf f}$ is a flux vector here, not a projected flux. As discussed in Ref.[\citen{nishikawa_liu_source_quadrature:jcp2017}], the 
dissipation term in the numerical flux,
\begin{eqnarray}
 {\Phi}_{jk} =  \frac{1}{2} \left[   {\bf f}(u_L)   +  {\bf f}(u_R)  \right] \cdot \hat{\bf n}_{jk}  - \frac{1}{2} \hat{\bf D}_n (  u_R -  u_L ),
\label{cfd_numerical_flux_2d_euler_appendix}
\end{eqnarray}
does not contribute to the second-order truncation error. Therefore, we will focus on the averaged flux term. To expand it, we write
\begin{eqnarray}
 {\bf f}(u_L)  &=& {\bf f}(u_j) +{\bf f}_{u} (u_L- u_j) +  \frac{1}{2} {\bf f}_{uu} (u_L- u_j)^2 + \frac{1}{6}  {\bf f}_{uuu} (u_L- u_j)^3 + O(h^4),
 \label{fL_expanded_du}
 \\[2ex]
 {\bf f}(u_R)  &=&  {\bf f}(u_j) + {\bf f}_{u} (u_R- u_j) +  \frac{1}{2} {\bf f}_{uu} (u_R- u_j)^2 + \frac{1}{6} {\bf f}_{uuu} (u_R- u_j)^3 + O(h^4),
  \label{fR_expanded_du}
  \end{eqnarray}
where ${\bf f}_{u} =\partial {\bf f}/ \partial u$, ${\bf f}_{uu} =\partial^2 {\bf f}/ \partial u^2$, and ${\bf f}_{uuu} =\partial^3 {\bf f}/ \partial u^3$. Consider the Taylor expansions around a node $j$ of an exact solution at a neighbor node $k$ and LSQ gradients at $k$ and $j$. 
\begin{eqnarray}
                           u_k &=& u_j + \djk u_j +  \frac{1}{2} \,  \ddjk u_j   +  \frac{1}{6} \,  \dddjk u_j  + O(h^4), \\ [2ex]
\overline{\nabla} u_k &=& \nabla u_k + {\cal E}(\overline\nabla u_k) 
   =  \nabla u_j + \djk (\nabla u_j) + \frac{1}{2} \,  \ddjk  u_j + {\cal E}(\overline{\nabla} u_j) +  O(h^3), 
   \label{expanded_fR} \\ [2ex]
\overline{\nabla} u_j &=& \nabla u_j + {\cal E}(\overline{\nabla} u_j) + O(h^3),
\end{eqnarray}
where the operator $\partial_{jk}$ is defined as
\begin{eqnarray}
\djk  =   ( {\bf x}_{k} - {\bf x}_j ) \cdot \nabla,
\end{eqnarray}
and ${\cal E}(\overline{\nabla} u_k) $ and ${\cal E}(\overline{\nabla} u_j)$ denote the vectors of 
the local truncation errors of the nodal gradients at $k$ and $j$, respectively. 
Substitute these expansions into the U-MUSCL solution reconstruction schemes, we obtain
 \begin{eqnarray}
u_L  
&=&    \kappa \frac{ u_{j}  + u_{k} }{2} + (1-\kappa) \left[  u_{j} +  \frac{1}{2} \overline{\nabla} u_j \cdot ( {\bf x}_{k} - {\bf x}_j )   \right]   \\[2ex] 
&=& u_j + \frac{1}{2} \djk u + \frac{\kappa}{4} \ddjk u  + \frac{\kappa}{12}  \dddjk u  + \frac{1-\kappa}{2} {\cal E}(\overline{\nabla} u_j)  \cdot ( {\bf x}_{k} - {\bf x}_j )  
+ O(h^4),
\end{eqnarray}
and
 \begin{eqnarray}
u_R  
&=&  \kappa \frac{ u_{k}  + u_{j} }{2} + (1-\kappa) \left[  u_{k} -  \frac{1}{2} \overline{\nabla} u_k \cdot ( {\bf x}_{k} - {\bf x}_j )   \right],   \\[2ex] 
&=& u_j + \frac{1}{2} \djk u + \frac{\kappa}{4} \ddjk u  + \frac{2 \kappa -1}{12}  \dddjk u  - \frac{1-\kappa}{2} {\cal E}(\overline{\nabla} u_j)  \cdot ( {\bf x}_{k} - {\bf x}_j )  
+ O(h^4).
\end{eqnarray}
Thus, we have, up to third-order order (which will contribute to 
the second-order error eventually), 
\begin{eqnarray}
u_L  - u_j &=&  \frac{1}{2} \djk u + \frac{\kappa}{4} \ddjk u  + \frac{\kappa }{12}  \dddjk u
  + \frac{1-\kappa}{2} {\cal E}(\overline{\nabla} u_j)  \cdot ( {\bf x}_{k} - {\bf x}_j )  + O(h^4), \\ [2ex]
(u_L  - u_j)^2 &=&  \frac{1}{4} \left( \djk u  \right)^2 + \frac{\kappa}{4}  ( \djk u) (  \ddjk u )
+ \frac{1-\kappa}{2}   ( \djk u ) {\cal E}(\overline{\nabla} u_j)  \cdot ( {\bf x}_{k} - {\bf x}_j ) + O(h^4), \\ [2ex]
(u_L  - u_j)^3 &=& \frac{1}{8} \left(  \djk u  \right)^3+ O(h^4).
\end{eqnarray}
and
\begin{eqnarray}
u_R  - u_j &=&  \frac{1}{2} \djk u + \frac{\kappa}{4}  \ddjk u  + \frac{2 \kappa -1}{12} \dddjk u
  - \frac{1-\kappa}{2} {\cal E}(\overline{\nabla} u_j)  \cdot ( {\bf x}_{k} - {\bf x}_j ) + O(h^4),  \\ [2ex]
(u_R  - u_j)^2 &=&  \frac{1}{4} \left(   \partial_{jk} u  \right)^2 + \frac{\kappa}{4}  ( \djk u) ( \ddjk u )
- \frac{1-\kappa}{2}   ( \djk u ) {\cal E}(\overline{\nabla} u_j)  \cdot ( {\bf x}_{k} - {\bf x}_j ) + O(h^4),  \\ [2ex]
(u_R  - u_j)^3 &=& \frac{1}{8} \left(  \djk u  \right)^3+ O(h^4).
\end{eqnarray}
Then, substituting these into Equations (\ref{fL_expanded_du}) and (\ref{fR_expanded_du}), we find 
\begin{eqnarray}
 \frac{1}{2} \left[   {\bf f}(u_L)   +  {\bf f}(u_R)  \right] \cdot \hat{\bf n}_{jk} 
 = \left[   
  {\bf f}(u_j) + \frac{1}{2}  {\bf f}_{u} \djk u 
  + {\cal E}_2 
  + {\cal E}_3
 \right] \cdot \hat{\bf n}_{jk}, 
\end{eqnarray}
or, because ${\bf f}_{u} \djk u  =  \djk \, {\bf f}$,
\begin{eqnarray}
 \frac{1}{2} \left[   {\bf f}(u_L)   +  {\bf f}(u_R)  \right] \cdot \hat{\bf n}_{jk} 
 = \left[   
  {\bf f}(u_j) + \frac{1}{2}  \djk  \, {\bf f}
  + {\cal E}_2 
  + {\cal E}_3
 \right] \cdot \hat{\bf n}_{jk}, 
\end{eqnarray}
where
\begin{eqnarray} 
{\cal E}_2 &=& 
\frac{1}{4} \left\{
\kappa {\bf f}_{u}  ( \ddjk  u )
+ {\bf f}_{uu}  (  \djk u )^2
  \right\}, 
  \\ [2ex]
{\cal E}_3 &=& 
\frac{1}{48} \left\{
    {\bf f}_{uuu}  (  \djk  u )^3 
    + 6 \kappa  {\bf f}_{uu} ( \djk  u )( \ddjk  u ) 
    + (6 \kappa-2)  {\bf f}_{u} ( \dddjk  u )
  \right\}.
\end{eqnarray}
For $\kappa=1/2$, the term ${\cal E}_3$ can be written as 
\begin{eqnarray} 
{\cal E}_3 =
\frac{1}{48} \left\{
    {\bf f}_{uuu}  ( \djk  u )^3 
    +3  {\bf f}_{uu} ( \djk  u )( \ddjk  u ) 
    +   {\bf f}_{u} ( \dddjk  u )
  \right\}
  = 
  \frac{1}{48}  \dddjk  \, {\bf f} 
  .
\end{eqnarray}
Therefore, the truncation error ${\cal T}$ is given by
\begin{eqnarray} 
{\cal T} = 
\frac{1}{V_j} \sum_{k \in \{ k_j\} }  {\Phi}_{jk} A_{jk}
= 
\frac{1}{V_j}
\sum_{k \in \{ k_j\} } 
 \left[   
  {\bf f}(u_j) + \frac{1}{2}  \djk \, {\bf f}
  + {\cal E}_2 
  + \frac{1}{48 V_j} \dddjk \,  {\bf f} 
 \right] \cdot {\bf n}_{jk} + O(h^3).
\end{eqnarray}
First, note that we have the following identities for simplex-element grids, 
\begin{eqnarray} 
\sum_{k \in \{ k_j\} }  {\bf n}_{jk} = 0 , \quad
\frac{1}{V_j} \sum_{k \in \{ k_j\} } \frac{1}{2} \Delta {x}_{jk} {\bf n}_{jk} =1,
\quad 
\frac{1}{V_j} \sum_{k \in \{ k_j\} } \frac{1}{2} \Delta {y}_{jk} {\bf n}_{jk} =1,
\end{eqnarray}
and therefore
\begin{eqnarray} 
{\cal T} =  
\mbox{div} \, {\bf f}_j + 
\frac{1}{V_j}
\sum_{k \in \{ k_j\} } 
  {\cal E}_2  
  \cdot {\bf n}_{jk}
 +
  \frac{1}{48 V_j}
 \sum_{k \in \{ k_j\} } \dddjk  \, {\bf f} 
 \cdot {\bf n}_{jk}  + O(h^3),
\end{eqnarray}
or, since $\mbox{div} \, {\bf f}_j  = 0$ for any exact solution,
 \begin{eqnarray} 
{\cal T} =  
\frac{1}{V_j}
\sum_{k \in \{ k_j\} } 
  {\cal E}_2  
  \cdot {\bf n}_{jk}
 +
  \frac{1}{48 V_j}
 \sum_{k \in \{ k_j\} } \dddjk  \, {\bf f} 
 \cdot {\bf n}_{jk}  + O(h^3),
\end{eqnarray}
Second, we note that the second term involves the following quantities,
which are all zero on regular simplex-element grids,
\begin{eqnarray} 
\sum_{k \in \{ k_j\} } \frac{1}{2} \Delta {x}^2 {\bf n}_{jk}
=
\sum_{k \in \{ k_j\} } \frac{1}{2} \Delta {y}^2 {\bf n}_{jk}
=
\sum_{k \in \{ k_j\} } \frac{1}{2} \Delta {z}^2 {\bf n}_{jk}
= 0, \\ [2ex]
\sum_{k \in \{ k_j\} } \frac{1}{2} \Delta {x}  \Delta {y}  {\bf n}_{jk}
=
\sum_{k \in \{ k_j\} } \frac{1}{2} \Delta {y} \Delta {z}  {\bf n}_{jk}
=
\sum_{k \in \{ k_j\} } \frac{1}{2} \Delta {z}  \Delta {x}  {\bf n}_{jk}
= 0, 
\end{eqnarray}
where $\Delta x = x_k - x_j$, $\Delta y = y_k- y_j$, and $\Delta z = z_k- z_j$, 
and thus it vanishes. Therefore, we are left with
\begin{eqnarray} 
{\cal T} =  
  \frac{1}{48 V_j}
 \sum_{k \in \{ k_j\} } \dddjk \,  {\bf f} 
 \cdot {\bf n}_{jk}  + O(h^3).
\end{eqnarray}
This is identical to the second-order truncation error of the third-order edge-based 
discretization on a regular simplex-element grid except the factor 
$1/48$, which is $-1/24$ for the third-order edge-based discretization 
Ref.[\citen{nishikawa_liu_source_quadrature:jcp2017}]. 
The rest of the proof for third-order accuracy without source terms 
can be found in Ref.[\citen{nishikawa_liu_source_quadrature:jcp2017}]; the second-order 
error term in the above can be written as a sum of second derivatives of $\mbox{div} \, {\bf f}$ (therefore it vanishes).
As proved in Ref.[\citen{nishikawa_liu_source_quadrature:jcp2017}], the factored form is given, in two dimensions,
\begin{eqnarray} 
{\cal T} 
=  
   \frac{1}{48 V_j} \left[ 
   Q_{xx} 
  + Q_{yy}
  + Q_{xy}   \right] ( \mbox{div} \, {\bf f} )_j,
  \label{factored_second_order_error_2d_2}
 \end{eqnarray}
where $Q_{xx}$, $Q_{yy}$, and $Q_{xy}$ are defined as
\begin{eqnarray}
Q_{xx}  =
\sum_{k \in \{ k_j\} } \!\!\!  n_x \Delta x^3 \partial_{xx} 
,    \quad
Q_{yy} =
\sum_{k \in \{ k_j\} } \!\!\! n_y \Delta y^3  \partial_{yy}
,    \quad
Q_{xy}  =
\sum_{k \in \{ k_j\} }\!\!\! 3 n_x \Delta x^2 \Delta y \, \partial_{xy}.
\end{eqnarray}
and in three dimensions,
\begin{eqnarray}
{\cal T} 
 =
 \! \frac{1}{48 V_j} \left[ 
         Q_{xx}   
  \!  + Q_{yy}   
  \!  + Q_{zz} 
  \!  + Q_{xy} 
  \!  + Q_{yz}   
  \!  + Q_{zx}      \right] ( \mbox{div} \, {\bf f} )_j   + O(h^3),
  \label{factored_second_order_error_3d_2}
 \end{eqnarray}
 where
 \begin{eqnarray}
Q_{xx}  =
\sum_{k \in \{ k_j\} } \!\!\!  n_x \Delta x^3 \partial_{xx} 
,    \quad
Q_{yy}  =
\sum_{k \in \{ k_j\} } \!\!\!  n_y \Delta y^3 \partial_{yy} 
,    \quad
Q_{zz} =
\sum_{k \in \{ k_j\} }  \!\!\!  n_z \Delta z^3 \partial_{zz},
\end{eqnarray}
 \begin{eqnarray}
Q_{xy}  =
 \sum_{k \in \{ k_j\} } \!\!\! 6 n_z \Delta x \Delta y \Delta z \, \partial_{xy} 
,    \quad
Q_{yz}  =
\sum_{k \in \{ k_j\} }\!\!\! 6 n_x \Delta x \Delta y \Delta z   \, \partial_{yz} 
,    \quad
Q_{zx} =
\sum_{k \in \{ k_j\} } \!\!\! 6 n_y \Delta x \Delta y \Delta z    \partial_{zx}.
\end{eqnarray} 
Here, we have omitted the subscript in the directed-area vector components for brevity: $ {\bf n}_{jk} = (n_x,n_y,n_z)$.

To achieve third-order accuracy for equations with source terms, we must have a compatible second-order error from the source discretization. Let us write the source discretization (\ref{residual_imuscl_s}) in the general form with $D=2$ for two dimensions and $D=3$ for three dimensions: 
\begin{eqnarray} 
\tilde{s}_j = \frac{1}{V_j} \sum_{k \in \{ k_j\} }  {\psi}_{jk} V_{jk},
\end{eqnarray}
 \begin{eqnarray}
 V_{jk} = \frac{1}{2 D} (  {\bf x}_k - {\bf x}_j  )\cdot {\bf n}_{jk} 
 ,
 \quad
  {\psi}_{jk} = 
   \kappa_s \frac{ {s}_j + {s}_k }{2} 
 + (1 - \kappa_s) \left[ {s}_j +
 \frac{1}{2} \overline{\nabla} {s}_j 
 \cdot (  {\bf x}_k - {\bf x}_j  ) \right]
.
\end{eqnarray}
Expanding the source term $s$ in the same way as we have done for $u$, and we obtain
the truncation error ${\cal T}_s $ of the source discretization as
\begin{eqnarray} 
{\cal T}_s 
&=& \frac{1}{V_j} \sum_{k \in \{ k_j\} } 
\left[
s_j +  \frac{1}{2}\partial_{jk} s + \frac{\kappa_s}{4} \partial_{jk}^2 s  +
O(h^3)  
\right] V_{jk} \\ [2ex]
&=&
s_j  
+ \frac{1}{2 V_j}  \sum_{k \in \{ k_j\} }  ( \partial_{jk} s )  V_{jk}
+  \frac{\kappa_s}{4 V_j  }  \sum_{k \in \{ k_j\} }   ( \partial_{jk}^2 s ) V_{jk} 
+ O(h^3).
\end{eqnarray}
As shown in Ref.[\citen{nishikawa_liu_source_quadrature:jcp2017}], the second term 
vanishes. Also shown in Ref.[\citen{nishikawa_liu_source_quadrature:jcp2017}], the third term, which is the second-order error, can be factored as 
\begin{eqnarray}
{\cal T}_s =  \frac{\kappa_s}{4 V_j} \left( \frac{1}{3} \right)  \sum_{k \in \{ k_j\} }
\left[
   Q_{xx} 
  + Q_{yy}
  + Q_{xy} 
\right]
s_j,
\label{te_s_expression_2d}
\end{eqnarray}
in two dimensions, and
\begin{eqnarray}
{\cal T}_s =  \frac{\kappa_s}{4 V_j} \left( \frac{5}{18} \right)  \sum_{k \in \{ k_j\} }
\left[
         Q_{xx}   
  \!  + Q_{yy}   
  \!  + Q_{zz} 
  \!  + Q_{xy} 
  \!  + Q_{yz}   
  \!  + Q_{zx}  
\right]
s_j, 
\label{te_s_expression_3d}
\end{eqnarray}
 in three dimensions. 
To achieve third-order accuracy, we must have, in two dimensions,
\begin{eqnarray} 
{\cal T} 
=
\frac{1}{V_j} \sum_{k \in \{ k_j\} }  {\Phi}_{jk} A_{jk} -  \frac{1}{V_j} \sum_{k \in \{ k_j\} }  {\psi}_{jk} V_{jk}
=  
  \frac{1}{48 V_j} \left[ 
   Q_{xx} 
  + Q_{yy}
  + Q_{xy}   \right] ( \mbox{div} \, {\bf f}_j  - s_j),
  \label{factored_second_order_error_2d_2_source}
 \end{eqnarray}
which is achieved if 
\begin{eqnarray} 
\frac{\kappa_s}{4} \left( \frac{1}{3} \right) = \frac{1}{48} ,
 \end{eqnarray}
thus
\begin{eqnarray} 
\kappa_s = \frac{1}{4}.
 \end{eqnarray}
Similarly, in three dimensions, we set 
\begin{eqnarray} 
 \frac{\kappa_s}{4 } \left( \frac{5}{18} \right)  = \frac{1}{48} ,
 \end{eqnarray}
 and thus
 \begin{eqnarray} 
\kappa_s = \frac{3}{10}.
 \end{eqnarray}


\end{document}